\documentclass[a4paper,12pt,twoside]{article}
\usepackage{amssymb,amsmath,latexsym,geometry,fancyhdr,lineno,hyperref,titletoc}

\usepackage{ulem}

\usepackage{import}
\usepackage{xifthen}
\usepackage{pdfpages}
\usepackage{transparent}
\usepackage{amsmath}
\usepackage{mathrsfs}
\contentsmargin{0pt}

\dottedcontents{section}[2em]{\vspace{-1mm}\small}{2em}{0pt}
\dottedcontents{subsection}[5em]{\vspace{-1mm}\small}{3em}{5pt}

\hypersetup{colorlinks,linkcolor= blue,citecolor=blue}

\newtheorem{theorem}{Theorem}[section]

\newtheorem{lemma}{Lemma}[section]

\newtheorem{proposition}{Proposition}[section]

\newtheorem{remark}{Remark}[section]

\newtheorem{claim}{Claim}[section]

\newtheorem{corollary}{Corollary}[section]

\newtheorem{definition}{Definition}[section]

\newcommand{\hdot}{^\text{\r{}}\hspace{-.33cm}H}

\numberwithin{equation}{section}
\allowdisplaybreaks[0]

\def\Vs{\vskip8pt}\def\vs{\vskip4pt}

\allowdisplaybreaks[4]

\begin{document}

\begin{center}
{\bf\Large Local uniqueness and non-degeneracy of blowup solutions for regular Liouville systems}
\end{center}

\vs\centerline{Zetao Cheng}  
\begin{center}
{\footnotesize
{Department of Mathematics and Research Institute for Natural Sciences,\\ College of Natural Sciences, Hanyang University,\\ 222 Wangsimni-ro Seongdong-gu, Seoul 04763, Republic of Korea\\
{\em E-mail}: chengzt20@mails.tsinghua.edu.cn}}
\end{center}

\vs\centerline{Haoyu Li}  
\begin{center}
{\footnotesize
{ Department of Mathematics, Faculty of Science and Technology,\\ University of Macau, Taipa, Macau, China\\
{\em E-mail}:  hyli1994@hotmail.com}}
\end{center}

\vs\centerline{Lei Zhang}
\begin{center}
{\footnotesize
{Department of Mathematics, University of Florida\\
1400 Stadium Rd, Gainesville FL 32611\\
{\em E-mail}:  leizhang@ufl.edu}}
\end{center}

\noindent{\bf\normalsize Abstract.} {\small
We study the following Liouville system defined on a compact Riemann surface $M$, 
\begin{equation}
     -\Delta u_i=\sum_{j=1}^n a_{ij}\rho_j\Big(\frac{h_j e^{u_j}}{\int_M h_j e^{u_j}}-1\Big)\mbox{ in }M\mbox{ for }i=1,\cdots,n,\nonumber
\end{equation}
where $M$ is a compact Riemann surface, the coefficient matrix $A=(a_{ij})_{n\times n}$ is nonnegative, $h_1, \ldots, h_n$ are positive smooth functions, and $\rho_1, \ldots, \rho_n$ are positive constants. For the blowup solutions, we establish their uniqueness and non-degeneracy based on natural assumptions.

    The main results significantly generalize corresponding results for single Liouville equations \cite{BartJevLeeYang2019,BartYangZhang20241,BartYangZhang20242}. To overcome several substantial difficulties, we develop certain tools and extend them into a more general framework applicable to similar situations. Notably, to address the considerable challenge of a continuum of standard bubbles, we refine the techniques from Huang-Zhang \cite{HuangZhang2022} and Zhang \cite{Zhang2006,Zhang2009} to achieve extremely precise pointwise estimates. Additionally, to address the limited information provided by the Pohozaev identity, we develop a useful Fredholm theory to discern the exact role that the Pohozaev identity plays for systems. The considerable difference between systems and a single equation is also reflected in the location of blowup points, where the uncertainty of the energy type of the blowup point makes it difficult to determine the sufficiency of pointwise estimates. In this regard, we extend our highly precise pointwise estimates to any finite order. This aspect is drastically distinct from analyses of single equations.

}



\tableofcontents

\section{Introduction and main results}
Let $(M,g)$ be a Riemann surface with metric $g$ and no boundary. 
We study bubbling solutions to the following system:
\begin{equation}\label{e:001}
-\Delta u_i=\sum_{j=1}^n a_{ij}\rho_j\Big(\frac{h_j e^{u_j}}{\int_M h_j e^{u_j}}-1\Big)\mbox{ in }M\mbox{ for }i=1,\cdots,n.
\end{equation}
Here, $M$ is a compact Riemann surface, $h_1,...h_n$ are positive smooth functions, $\rho_1,..,\rho_n$ are positive constants, $(a_{ij})_{n\times n}$ is a matrix satisfying 
\begin{itemize}
    \item [$(\mathcal{H}_1)$] $A$ is non-negative, invertible, symmetric and irreducible;
    \item [$(\mathcal{H}_2)$] $a^{ii}\leq0$, $a^{ij}\geq0$ and $\sum_{j=1}^na^{ij}\geq 0$  for any $i,j=1,\cdots,n$ with $i\neq j$. 
\end{itemize}
Here, $a^{ij}$ denotes the $(i,j)$-element of the inverse of $A$.
Since $A$ is a nonnegative matrix, the system (\ref{e:001}) is called a Liouville system. 

     Liouville systems have numerous applications across a wide range of areas. In the field of geometry, when simplified to a single equation ($n=1$), they extend the famous Nirenberg problem, which has been the subject of extensive study over the last few decades (see \cite{barto2,barto3,bart-taran-jde,bart-taran-jde-2,ccl,chen-li-duke,kuo-lin-jdg,li-cmp,li-shafrir,wei-zhang-adv,wei-zhang-plms,Zhang2006,Zhang2009}). In the realm of physics, Liouville systems arise from the mean field limit of point vortices in Euler flow (see \cite{biler,wolansky1,wolansky2,yang}) and are closely associated with self-dual condensate solutions in the Abelian Chern-Simons model with $N$ Higgs particles \cite{kimleelee,Wil}. In biological contexts, they appear in the stationary solutions of the multi-species Patlak-Keller-Segel system \cite{wolansky3} and play a vital role in the study of chemotaxis \cite{childress}. Grasping the concept of bubbling solutions, a major challenge in Liouville systems, is essential for making progress in related fields.

In comparison Toda systems are completely integrable systems used in various fields including including solid-state physics, mathematical physics, and in the study of integrable systems and cluster algebras, etc (see \cite{takasaki,lin-wei-ye}). Here we mention three main features of Liouville systems that make them drastically different from Toda systems. One is Liouville systems are enormously inclusive. Various models in different fields can be described as Liouville systems. The second is the total integrations on all components of a Liouville system form a hypersurface, in other words, the energy is a continuum. The third main feature is that the method of moving planes can be emplyed on Liouville systems to prove the radial symmetry of solutions. In comparison solutions of Toda systems generally do not have symmetry but the integrals of global Toda systems form discrete points.

Using $I=\{1,\cdots,n\}$ we identity an important quantity
\[\Lambda_{I,N}=\Lambda_{I,N}(\rho)=4\sum_{i=1}^n\frac{\rho_i}{2\pi N}-\sum_{i=1}^n\sum_{j=1}^na_{ij}\frac{\rho_i}{2\pi N}\frac{\rho_j}{2\pi N}.\] 
and we use $\Lambda_{I,N}$ to define
\begin{equation*}
\Gamma_{N}=
\left\{
\rho\,\big|\, \rho_i>0, i\in I;\,\,  \Lambda_{I,N}(\rho)=0.
\right\}.
\end{equation*}
Clearly, for any solution $(u_1,\cdots,u_n)$ to
 (\ref{e:001}),  $(u_1+C_1,\cdots,u_n+C_n)$ is also a solution for any constants $C_i$ ($i=1,\cdots,n$). Therefore we assume that all the solutions of (\ref{e:001}) belong to the following space: 
 \[\hdot^{1,n}(M)= \hdot^{1}(M)\times....\times  \hdot^{1}(M)\] where 
$$ \hdot^{1}(M):=\{v\in L^2(M);\quad \nabla v\in L^2(M), \mbox{and }\,\, \int_M v dV_g=0\}. $$

In this paper, we study blowup solutions defined as follows.
\begin{definition}
A sequence of solutions $u^k=(u^k_1,\cdots,u_N^k)\in \hdot^{1}(M)$ to
\begin{equation}\label{e:00k}
     -\Delta u^k_i=\sum_{j=1}^n a_{ij}\rho^k_j\Big(\frac{h_j e^{u^k_j}}{\int_M h_j e^{u^k_j}}-1\Big)\mbox{ in }M\mbox{ and }\mbox{ and }k=1,2,\cdots
\end{equation}
is called \emph{blowup solutions} if there exist $N$ disjoint points $q^*_1,\cdots,q^*_N\in M$ such that
\begin{itemize}
    \item [$(1)$] $u_i^k\to-\infty$ on any compact subset of $M\backslash\{q^*_1,\cdots,q^*_N\}$ as $k\to\infty$;
    \item [$(2)$] for $i=1,\cdots,n$, $\frac{h_i e^{u_i^k}}{\int_\Omega h_i e^{u_i^k}}\to\sum_{t=1}^N \alpha_t\delta_{q^*_t}$ as $k\to\infty$. Here, the constant $\alpha_t>0$ and $\delta_{q^*_t}$ is the Dirac measure supported by $\{q^*_t\}$ for $t=1,\cdots,N$.
\end{itemize}
\end{definition}
It is proved in \cite{LinZhang2013,HuangZhang2022,GuZhang2025} that if there exists a sequence of blowup solutions $u^k=(u^k_1,\cdots,u^k_n)$ to (\ref{e:00k})
then $\rho^k:=(\rho^k_1,\cdots,\rho^k_n)\to\rho^*\in\Gamma_{N}$ as $k\to\infty$ for some $N$. Here, $N$ denotes the number of blowup points.
For any blowup solution $u^k=(u_1^k,\cdots,u_N^k)$ of (\ref{e:001}) we denote
\begin{align}
M^{k}_{t}:=\max_{i\in I}\max_{x\in B_{\tau}(q_t^*)}\Big\{u_i^k(x)-\ln\int_M h_i e^{u_i^k(x)}dx\Big\}\quad \mbox{ and }\quad
\varepsilon^{k}_{t}:=e^{-\frac{M^{k}_{t}}{2}}.\nonumber
\end{align}
where $\delta>0$ is small.  It is well known that blow up points have to be the critical point of the following function:
\begin{equation}\label{c-f}f(x_1,...,x_N):=\sum_{t=1}^N\sum_{i=1}^n\rho_i\Big[\ln h_i(x_{t})+\pi m_i^* \gamma(x_t,x_t)+\sum_{s,s\neq t}\pi m_i^* G(x_t,x_s)\bigg ]. 
\end{equation}
So we assume:
\[(G):\quad (q_1^*,\cdots,q_N^*) \mbox{  is a non-degenerate critical point of $f(x_1,...,x_N)$ in (\ref{c-f})}.
\]
In the following for simplicity, we omit $k$ in the notation and use $\varepsilon$ to denote $\varepsilon_1,\cdots,\varepsilon_N$ and use $\varepsilon_t$, $u_{i,\varepsilon}$, $u_\varepsilon$, $\rho_{i,\varepsilon}$, $\rho_\varepsilon$ to denote $\varepsilon_{k,t}$, $u_i^k$, $u^k$, $\rho_i^k$, $\rho^k$, respectively.
Let $\rho^*\in\Gamma_N$ be the limit of $\rho_\varepsilon$ with
$\Lambda_{I,N}(\rho^*)=0$.
 and $q^k=(q_{1}^{k},\cdots,q_{N}^{k})$. Let
\[m_i^*=\sum_{j=1}^na_{ij}\frac{\rho_j^*}{2\pi N}
\quad\mbox{for}\quad  i=1,\cdots,n,
\quad \mbox{and}\quad
m^*=\min\{m^*_i|i=1,\cdots,n\}.\]

The profiles of blowup solutions are significantly different for $m^*<4$ and $m^*=4$, respectively.

For any $p\in M$, let $G(x,p)$ be the Green's function at $p$:
\begin{equation}\label{e:GreensFunction}
     -\Delta_x G(x,q)=\delta_{q} -1\mbox{, }\int_M G(x,q)dx=0,
\end{equation}
where $\gamma$ denotes the regular part of $G$, represented as
\[\gamma(x,q)=G(x,q)+\frac{1}{2\pi} \ln|x-q|, \quad \mbox{for $x$ close to $q$}.\]
 We set for any $t=1,\cdots,N$ and $q=(q_1,\cdots,q_N)\in M^N$,
\begin{align}\label{def:Gt}
G^*(x;q_t)=\gamma(x,q_t)+\sum_{s\neq t}G(x,q_s).
\end{align}

For any $q=(q_1,\cdots,q_N)\in M^N$ with $q_t\neq q_s$ if $t\neq s$ and $t,s=1,\cdots,N$, we define
\begin{align}\label{def:H}
H_{it}(q)=2\pi m_i^*G^*(q_t;q_t) +\ln h_i(q_t),
\end{align}
and let $\Omega_t$ ($t=1,\cdots,N$) be $N$ sub-domains of $M$ such that $q_t^*\in\Omega_t$, $\overline{\cup_{t=1}^N\Omega_t}=M$ and $\Omega_{t_1}\cap\Omega_{t_2}=\emptyset$ for $t_1,t_2=1,\cdots,N$ with $t_1\neq t_2$. 
To go to the main theorem, we need the following further notions. For any $i\in I$ and any $t=1,\cdots,N$, we denote
\begin{align}
D_{it}=\lim_{\tau\to0+}\Bigg[-\frac{2\pi e^{I_i}}{m^*_i -2}\tau^{2-m^*_i}+e^{I_i}\int_{\Omega_t\backslash B_\tau^{M}(q_t^*)} \frac{h_i(x)}{h_i(q^*_t)}e^{2\pi m_i^*[G^*(x;q_t^*)-G^*(q_t^*;q_t^*)]}dV_g\Bigg]\nonumber
\end{align}
and
\begin{align}\label{def:Lit}
L_{it}=\Delta(\ln h_i (q_t^*))-2K(q_t^*)+8\pi N+\Big|\nabla (\ln h_i(q_t^*))+8\pi  \nabla_1 G^*(q_t^*;q_t^*)\Big|^2.
\end{align}
where $\nabla_1$ means the derivatives with respect to the first $q_t^k$, $I_i$ are constants in the expansion of a global solution (see (\ref{e:LimitSystemAsy})).
Using $D_{it}$,$L_{i,t}$ and $m^*$ we postulate the following natural assumption:
\begin{align}\label{nature-a}
\left\{
\begin{aligned}
& A: \, m^*<4 \, \mbox{and}\, \sum_{i=1}^n \sum_{t=1}^N D_{it}\frac{e^{H_{i,t}(q^*)}}{e^{H_{i,1}(q^*)}}\neq0;\\
 &   B: \, m^*=4 \,\mbox{and}\, \sum_{i=1}^n \sum_{t=1}^N L_{it}\frac{e^{H_{i,t}(q^*)}}{e^{H_{i,1}(q^*)}}\neq0;\\
  &  C:\, m^*=4 \,\mbox{and} \,\sum_{i=1}^n \sum_{t=1}^N L_{it}\frac{e^{H_{i,t}(q^*)}}{e^{H_{i,1}(q^*)}}=0 \, \mbox{ and } \,\sum_{i=1}^n \sum_{t=1}^N  D_{it}\frac{e^{H_{i,t}(q^*)}}{e^{H_{i,1}(q^*)}}\neq0. 
    \end{aligned}
    \right.
\end{align}

\begin{theorem}\label{t:MAIN} Let $u^{(1),k}=(u_1^{(1),k},..,u_n^{(1),k})$ and $u^{(2),k}=(u_1^{(2),k},...,u_n^{(2),k})$ be any two blowup solutions of (\ref{e:00k}) with the same parameters $\rho^k=(\rho_1^k,..,\rho_n^k)$ and $q^*=(q^*_1,\cdots,q^*_N)$.
Suppose $(\mathcal{H}_1)$, $(\mathcal{H}_2)$ and $(G)$  hold and one of the three options in 
(\ref{nature-a}) occurs, then there exists $k_0>1$ such that $u_i^{(1),k}=u_i^{(2),k}$ for all $k\ge k_0$ and $i=1,...,n$.
\end{theorem}
The second main result is the following non-degeneracy theorem.

\begin{theorem}\label{non-deg}
Let $u^k=(u^k_1,..,u^k_n)$ be a sequence of blowup solutions of (\ref{e:00k}) and $\phi^k=(\phi_1^k,...,\phi_n^k)$ be a sequence of solutions of 
\begin{align}
\Delta\Big(\sum_j a^{ij}\phi^k_j\Big)+\rho^k_i\frac{h_ie^{u_i^k}}{\int_M h_ie^{u_i^k}}\Bigg(\phi^k_i-\int_M\frac{h_i e^{u^k_i}\phi^k_i}{\int_M h_i e^{u^k_i}}\Bigg)=0\mbox{ in }M\mbox{ for }i=1,\cdots,n,
\end{align}
where  $\phi^k_i\in \mathring{H}^1(M)=\{\phi\in H^1(M)|\int_M \phi=0\}$ for  $i=1,\cdots,n$. Then, under the same assumptions as in Theorem \ref{t:MAIN}, there exists $k_0>1$ such that for all $k\ge k_0$, $\phi_i^k\equiv 0$ for $i=1,...,n$.
\end{theorem}
\begin{remark}
       The term $L_{it}$ pertains only to the curvature information of the blowup points, whereas $D_{it}$ involves an integration over the entire manifold. According to the main theorems, when $m^* < 4$, only $D_{it}$ is significant. However, when $m^* = 4$, $L_{it}$ becomes predominant. If the term involving $L_{it}$ is zero, then the global integration $D_{it}$ assumes the leading role.
\end{remark}

       Theorems \ref{t:MAIN} and \ref{non-deg} serve as comprehensive results encompassing all blowup cases for regular Liouville systems. Achieving such generality involves overcoming at least four significant challenges. Firstly, the continuous nature of energy in Liouville systems complicates blowup analysis, particularly when multiple bubbles occur simultaneously. Our work uncovers strong similarities between profiles around different blowup points, providing a complete solution to this issue. Secondly, while the Pohozaev identity is crucial in many curvature prescribing equations, its role is limited in systems due to the difficulty of using a single Pohozaev identity to govern multiple equations. We propose viewing the Pohozaev identity as an orthogonal condition within a Fredholm theory for systems, developing a Fredholm theory for specific linearized operators to fully understand the Pohozaev identity's role.  From this perspective, we identify a non-vanishing blowup phenomenon in Liouville systems, as detailed in Subsection \ref{subsectionb:NonvanishBlowup}, which is markedly distinct from the single equation scenario. The third challenge arises from the weak decay rate of blowup functions for systems. Unlike single Liouville equations where solutions rapidly decay to negative infinity, Liouville systems do not behave similarly. To prove the main theorems and address all possible scenarios, the weak decay rate necessitates a highly precise analysis of the asymptotic behavior of blowup solutions. The fourth major difficulty comes from the procedure in the proof of the main theorems. One needs to prove three parameters $b_0,b_1,b_2$ to be zero. For the Liouville equation with non singular source, the main task is to prove $b_0=0$ and the other two parameters can be proved equal zero easily. However in this article, we have to prove a precise vanishing rate of $b_0^k$ in order to prove that other parameters vanish. 
       We anticipate that the methods introduced in this article will apply to a variety of similar problems.

  The paper is organized as follows: Section 2 lists some preliminary results from previous works. Sections 3 and 4 conduct a precise asymptotic analysis of the profiles of bubbling solutions. It is crucial to reduce the error to be below a sufficiently small threshold to ensure a good match of bubbles at different blowup points. The proofs of the main theorems are provided in Section 5, and Section 6 studies the Pohozaev identity within the framework of a Fredholm alternative and a related construction of blowup solutions with non-vanishing first derivatives of coefficient functions.

\section{Preliminary Results}

In this section we mention some previous results and study the limit system in more detail.  First we define $p_{t,\varepsilon}\in B_\tau(q_t^*)$ to be the point where the following maximum is attained: 
\[\max_{i\in I}\max_{x\in B_{\tau}(q_t^*)}\Big\{u_i^k(x)-\ln\int_M h_i e^{u_i^k(x)}dx\Big\}.\]

Then we recall an important result on the location of blowup points:
\begin{align}\label{SharpEstimate1}
\sum_{i=1}^n\Big[\nabla \ln h_i(q_{t}^{k})+2\pi m_i^* \nabla_1 G^*(q_t^k,q_t^k)\Big]\rho_i^*=\left\{
\begin{aligned}
&O(\varepsilon_1^{m^*}) &\mbox{ if } m^*<4;\\
&O(\varepsilon_1^2\ln\frac{1}{\varepsilon_1})&\mbox{ if } m^*= 4,
\end{aligned}
\right.
\end{align}
 (\ref{SharpEstimate1}) was proved by Lin-Zhang \cite{LinZhang2013} for one-bubble case and by Huang-Zhang \cite{HuangZhang2022} and Gu-Zhang \cite{GuZhang2025} for multi-bubble case.

The second result of the \emph{sharp estimates} states that
\begin{align}\label{SharpEstimate2}
\Lambda_{I,N}(\rho_\varepsilon)=\left\{
\begin{aligned}
&-\sum_{i\in\hat{I}}\frac{m^*-2}{N}\Big(D_{i,1}+\sum_{t>1}\frac{e^{H_{i,t}(q^*)}}{e^{H_{i,1}(q^*)}} D_{i,t}+o_\tau(1)\Big)\varepsilon_1^{m^*-2} &\mbox{ if } m^*<4;\\
&-\sum_{i=1}^n\frac{2}{N}\Big(L_{i,1}+\sum_{t>1}\frac{e^{H_{i,t}(q^*)}}{e^{H_{i,1}(q^*)}} L_{i,t}+o_\tau(1)\Big)\varepsilon_1^2\ln\frac{1}{\varepsilon_1} &\mbox{ if } m^*= 4.
\end{aligned}
\right.
\end{align}
See \cite{LinZhang2013,HuangZhang2022,GuZhang2025}. The functions $H_{i,t}$ are defined in (\ref{def:H}).
The corresponding limit problem is \begin{equation}\label{e:LimitSystem}
    \left\{
   \begin{array}{lr}
     -\Delta v_i=\sum_{j=1}^n a_{ij}e^{v_j}\mbox{ for }x\in\mathbb{R}^2,\\     \frac{1}{2\pi}\int_{\mathbb{R}^2}e^{v_i}=\sigma_i \mbox{ for }i\in I=\{1,\cdots,n\}.
   \end{array}
   \right.
\end{equation}
The entire solution $V=(V_1,\cdots,V_n)$ satisfies
\begin{equation}\label{e:LimitSystemAsy}
    \left\{
   \begin{array}{lr}
     V_i(x)=-m_i^* \ln|x|+I_i-\sum_{j=1}^n\frac{a_{ij}e^{I_j}}{(m_j-2)^2}|x|^{2-m_j^*}+\sum_{j,l=1}^n\frac{a_{ij}a_{jl}e^{I_j+I_l}}{(m_l-2)^2(m_j+m_l-4)^2}|x|^{4-m_j-m_l}+O(|x|^{6-3m^*}),\\
     DV_i(x)=-m_i^*\frac{x}{|x|^2}+\sum_{j=1}^n\frac{a_{ij}e^{I_j}}{m_j-2}|x|^{-m_j^*}x+O(|x|^{3-2m^*})\mbox{ for }i\in I=\{1,\cdots,n\}
   \end{array}
   \right.
\end{equation}
for $|x|$ large, and
\begin{equation}\label{e:LimitSystemAsy1}
    \left\{
   \begin{array}{lr}
     V_i(x)=V_i(0)-\frac{1}{4}\sum_{j=1}^n a_{ij}e^{V_j(0)}|x|^2+\frac{1}{196}\sum_{j=1}^n\sum_{l=1}^n a_{ij}a_{jl}e^{V_j(0)+V_{l}(0)}|x|^4+O(|x|^6),\\
     DV_i(x)=-\frac{1}{2}\sum_{j=1}^n a_{ij}e^{V_j(0)}x+\frac{1}{49}\sum_{j=1}^n\sum_{l=1}^n a_{ij}a_{jl}e^{V_j(0)+V_{l}(0)}|x|^2x+O(|x|^5)
   \end{array}
   \right.
\end{equation}
for $i\in I=\{1,..,n\}$ and $|x|$ small. Here we remark that the error term $O(|x|^{6-3m^*})$ can be specified and the order of the error can be made smaller if necessary. 

A classical result for $(V_1,\cdots,V_n)$ is that
\begin{proposition}\label{p:Kernel}
If function $\phi=(\phi_1,\cdots,\phi_n)$ satisfies
\begin{equation*}
    \left\{
   \begin{array}{lr}
    \Delta\phi_i+\sum_{j=1}^n a_{ij}e^{V_j}\phi_j=0\mbox{ in }\mathbb{R}^2,\\
     |\phi_i(z)|\leq C(1+|z|)^\tau\mbox{ for }i\in I=\{1,\cdots,n\}
   \end{array}
   \right.
\end{equation*}
for $\tau\in (0,1)$. Then
\begin{align}
\phi\in\mbox{span}\{(Z_{0,1},\cdots,Z_{0,n}),(Z_{1,1},\cdots,Z_{1,n}),(Z_{2,1},\cdots,Z_{2,n})\}.\nonumber
\end{align}
Here, 
\begin{equation}\label{kernels}
Z_{0,i}=|x|\dot{V}_i(|x|)+2, \quad Z_{1,j}=\partial_{x_1} V_i,\quad Z_{2,i}=\partial_{x_2}V_i,\quad i=1,...,n.
\end{equation} 
\end{proposition}

The proof of Proposition \ref{p:Kernel} can be found in \cite{ChengLiZhang2025,GuZhang2025,LinZhang2013}.


\noindent{\bf Notations.} Before we proceed we mention some notations for clarity.
\begin{itemize}
    \item For some $\tau>0$, we denote the variable in the coordinate of $B_\tau(0)$ by $x$ and points $q$;
    \item For some $\tau>0$, we denote the variable in the coordinate of $B_{\tau/\varepsilon^k_t}(0)$ by $y$ and points $p$;
    \item In this paper, we use several different scales, say $\varepsilon_t^{(l),k}$, the scale of the $t$-th bubble of the solution $u^{(l),k}$ for $l=1,2$ and $k=1,2,\cdots$. However, we prove that $\varepsilon_t^{(l),k}\sim\varepsilon_s^{(l'),k}$ for any $t,s,=1,\cdots,N$ and $l,l'=1,2$. Therefore, for the error terms in the form of $O((\varepsilon_t^{(l),k})^\alpha)$ for some $\alpha\in\mathbb{R}$, we always write $O(\varepsilon^\alpha)$ for short;
    \item In this paper, we always use the error terms in the form of $O(\varepsilon^{\alpha-\delta})$ and $O((1+|y|)^{\alpha+\delta})$ for some $\alpha\in\mathbb{R}$. We assume that the positive number $\delta$ is much smaller than $|\alpha|$;
    \item We use $M$ to denote a generically large integer which might be different from line to line;
    \item We frequently use higher-order expansions and their continuous dependence on parameters, while the specific expression of the expansion itself is irrelevant. Therefore in this paper, we always denote the sum of higher-order terms as $T_{M,k}$ with certain superscripts or subscripts. We refer, for instance, Theorems \ref{t:IntegralExpansion} and \ref{t:IntegralExpansion}, Propositions \ref{prop:OutAverageExpansion} and \ref{prop:Pohozaev} and Corollaries \ref{coro:RefinedRate} and \ref{coro:TauMk}.
\end{itemize}

\section{A higher order expansion around a blowup point}

\subsection{An initial coordinate and classical estimates}\label{subsection:ClassicalBubbling}

In what follows, we assume without loss of generality that $\int_M h_i e^{u_i^k}dV_g=1$. Otherwise, we can change $u_i^k$ to
$u_i^k-\ln\int_M h_i e^{u_i^k} dV_g$
for $i=1,\cdots,n$. Then, $u_i^k$ solves
\begin{equation}
    \left\{
\begin{array}{lr}
    -\Delta u_i^k=\sum_{j=1}^n a_{ij}\rho_j^k(h^k_j e^{u_j^k}-1)\mbox{ in }M\\
     \int_M h^k_ie^{u_i^k}=1\mbox{ for }i\in I=\{1,2,\cdots,n\}.
   \end{array}
   \right.
\end{equation}
Here we take $\tau_0>0$ sufficiently small and let $q_{0,t}^k$ be the maximum point of $u_1^k|_{B_{2\tau_0}(q_t^*)}$ and we denote $q_t^*$ as the limit of $q_{0,t}^k$: $q_t^*=\lim_{k\to\infty}q_{0,t}^k$.
For $t=1,\cdots,N$, we define the isothermal coordinate system around $q_{0,t}^k$ as follows:
\begin{equation}\label{def:psi}
    \left\{
   \begin{array}{lr}
     ds^2=e^{\psi}(dy_1^2+dy_2^2),\\
     \psi(0)=|\nabla \psi(0)|=0,\\
     \Delta\psi=-2Ke^{\psi}.
   \end{array}
   \right.
\end{equation}
where, $K$ denotes the Gauss curvature of $M$.
Define
\begin{equation}\label{def:def:fik}
    \left\{
   \begin{array}{lr}
     -\Delta f_i^k=\sum_{j=1}^n a_{ij}\rho_j^k e^{\psi}\mbox{ in }B_{\tau_0},\\
     f_i^k(0)=|\nabla f_i^k(0)|=0\mbox{ for }i\in I.
   \end{array}
   \right.
\end{equation}
Define the function
\begin{align}
v_i^{0,k}(x)=u_i^k(q_{0,t}^k+x)-f_i^k(x)+\ln\rho_i^k+\ln h_i^k(q_{0,t}^k).\nonumber
\end{align}
This solves the equation
\begin{align}
-\Delta v_i^{0,k}=\sum_{j=1}^n a_{ij} H^{0,k}_j( y)e^{v_i^{0,k}}\mbox{ in }B_{\tau_0}(0)\nonumber
\end{align}
with $H_i^{0,k}(\cdot)=h^k_i(q_{0,t}^{k}+\cdot)e^{\psi+f_i^k}$ as we define above.
Define the harmonic function to encode the oscillation of $v_i^{0,k}$ on $\partial B_{\tau_0}$:
\begin{equation}\label{def:HarmonicFunction}
    \left\{
   \begin{array}{lr}
     -\Delta \phi_{i,t}^{0,k}=0\mbox{ in }B_{\tau_0},\\
     \phi_{i,t}^{0,k}|_{\partial B_{\tau_0}}=v_i^{0,k} -\frac{1}{2\pi\tau_0}\int_{\partial B_{\tau_0}}v_i^{0,k} \mbox{ for }i\in I.
   \end{array}
   \right.
\end{equation}
Let $q_{1,k}^k:=\varepsilon_t^k p_{1,t}^k\in B_{\tau_0}(q_{0,t}^k)$ be the maximum point of the function $v_i^{0,k}-\phi_{i,t}^{0,k}$. By the classical estimate (see \cite[Lemma 4.1]{LinZhang2013}), we have $|p_{1,t}^k|=O(\varepsilon^k_t)$. Define the number
\begin{align}
\varepsilon_t^k=exp\{-\max_{i\in I}\max_{x\in B(q_{0,t}^k,\tau_0)}(v_i^{0,k}(x)-\phi_{i,t}^{0,k}(x))/2\}\nonumber
\end{align}
and
\begin{align}
v_i^{1,k}(y)=v_i^{0,k}(\varepsilon_t^k y+q_t^{1,k})+\ln\rho_i^k+\ln h_i^k(q_t^{1,k})-\phi_{i,t}^{0,k}(\varepsilon_t^k y)+2\ln\varepsilon_t^k\nonumber
\end{align}
for $|y|\leq\tau_1^k(\varepsilon_t^k)^{-1}$ with $\tau_0=\tau_1^k+\varepsilon_t^k|p_{1,t}^{k}|$. It is evident that $v^{1,k}=(v_1^{1,k},\cdots,v^{1,k}_n)$ solves
\begin{align}\label{e:Equationv1k}
-\Delta v_i^{1,k}=\sum_{j=1}^n a_{ij} H^{1,k}_j( y)e^{v_j^{1,k}}\mbox{ in }B_{\tau_1^k}(0)
\end{align}
with \[H_i^{1,k}(\cdot)=h^k_i(\varepsilon_t^k \cdot+q_{1,t}^{k})e^{\psi+f_i^k+\phi_{t,i}^{0,k}(q_{1,t}^k+\cdot)}.\]
Let the function $V^{1,k}_t=(V_{t,1}^{1,k},\cdots,V_{t,n}^{1,k})$ be the solution to
\begin{equation}\label{e:LocalStandardBubble}
    \left\{
   \begin{array}{lr}
     -\Delta V_{t,i}^{1,k}=\sum_{j=1}^n a_{ij}e^{V_{t,j}^{1,k}}\mbox{ in }\mathbb{R}^2.\\
     V_{t,i}^{1,k}(0)= v_i^{1,k}(q_{1,t}^k)+\ln\rho_i^k+\ln h_i^k(q_{1,t}^k)-\phi_{i,t}^{1,k}(q_{1,t}^k) \mbox{ for }i\in I.
   \end{array}
   \right.
\end{equation}
Based on the above, we recall some of the classical estimates.

~~

\noindent{\bf Comparison between bubbles}

~~

Here we list several pieces of information concerning the comparison between bubbles.
Define the local masses
\begin{align}
\sigma_{t,i}^{1,k}=\frac{1}{2\pi}\int_{B(p_t^{1,k},\tau_1^k)}h_i^k e^{u_i^k},\quad m_{t,i}^{1,k}=\sum_{j=1}^n a_{ij}\sigma_{t,j}^{1,k}.\nonumber
\end{align}
Passing $k$ to infinity, we get
\begin{align}
\sigma_{t,i}=\lim_{k\to\infty}\sigma_{t,i}^{1,k}\mbox{ and }m_{t,i}=\lim_{k\to\infty}m_{t,i}^{1,k}.\nonumber
\end{align}
By \cite[pp. 13 $\And$ Proposition 3.1]{HuangZhang2022}, we get
\begin{align}\label{e:LocalMass}
m_{i,t}^{1,k}=\sum_{j=1}^N\frac{a_{ij}\rho_j^k}{2\pi N}+O((\varepsilon_1^k)^{m^*-2-\delta})
\end{align}
for any $i=1,\cdots,n$ and any $t=1,\cdots,N$.  Therefore, $\sigma_{t,i}=\sum_{j=1}^N\frac{a_{ij}\rho_j^*}{2\pi N}$ for any $i=1,\cdots,n$ and any $t=1,\cdots,N$.
For the standard bubble (\ref{e:LocalStandardBubble}), we get
\begin{align}\label{e:BubbleMass0}
\int_{\mathbb{R}^2} e^{V_{t,i}^k}=\frac{\rho_i^k}{N}+O((\varepsilon_1^k)^{m^*-2-\delta})
\end{align}
for any $i=1,\cdots,n$ and any $t=1,\cdots,N$.

Let us recall the computation in \cite[pp. 16-18]{HuangZhang2022}.
We consider $V_{1,i}^k$ and $V_{2,i}^k$ for examples. By (\ref{e:LimitSystemAsy}), we get
\begin{align}\label{Expansion:BubblesComparison}
\left\{
\begin{aligned}
V_{1,i}^k(x)&=-m_{1,i}^k\ln|x|-\frac{m_{1,i}^k-2}{2}M_{1}^{k}+I_{1,i}^k
+O((\varepsilon_1^k)^{m^*-2}),\\
\overline{V}_{2,i}^k(x)&:=V_{2,i}(
\eta x)-2\ln\eta_k=-m_{2,i}^k\ln|x|-\frac{m_{2,i}^k-2}{2}M_{2}^{k}+I_{2,i}^k
+O((\overline{\varepsilon}_1^k)^{m^*-2})
\end{aligned}
\right.
\end{align}
with $\varepsilon_1^k=e^{-\frac{M_i^k}{2}}$ and $\overline{\varepsilon}_1^k=e^{-\frac{M_2^k}{2}}$. By \cite[Theorem 1.1]{HuangZhang2022}, $\varepsilon_1^k\sim\overline{\varepsilon}_1^k$.    We further note that the relationships between the profiles of different bubbles are derived from precise pointwise estimates, both near and away from the blowup points.

Denote $\eta_k:=e^{\frac{M_{1}^k -M_2^k}{2}}$ and $\eta_k=O(1)$. Then, we claim that
\begin{claim}\label{c:BubblesCompare}
Let $\alpha^k:=\sum_i|\alpha^k_i-\overline{\alpha}^k_i|$ and $w^k_i:=V_{1,i}^k-\overline{V}_{2,i}^k$. Here, $\alpha^k_i=V_{1,i}^k(0)$ and $\overline{\alpha}_i^{k}:=\overline{V}_{2,i}^k(0)$.
It holds that
\begin{equation}\label{wes}
|w_i^k(y)|\leq C((\varepsilon_1^k)^\delta\alpha^k+(\varepsilon_1^k)^{m^*-2})(1+|y|)^{2\delta}\nonumber
\end{equation}
for $|y|\leq C/\varepsilon_1^k$. Here, $C$ is a positive constant and $\delta>0$ is sufficiently small.
As a consequence $\alpha^k=O((\varepsilon_1^k)^{m^*-2})$. 
\end{claim}
\noindent{\bf Proof.}
Suppose that
\begin{align}
\Lambda_k:=\max_i\sup_{|y|\leq{C/\varepsilon_1^k}}\frac{|w_i^k(y)|}{((\varepsilon_1^k)^\delta\alpha^k+(\varepsilon_1^k)^{m^*-2})(1+|y|)^{2\delta}}\to+\infty.\nonumber
\end{align}
We denote the maximum attended at $r=r_k$.
On $\partial B_{\tau/\varepsilon_k}$, (\ref{Expansion:BubblesComparison}) implies that
\begin{align}
\max_i\frac{|w_i^k(y)|}{((\varepsilon_1^k)^\delta\alpha^k+(\varepsilon_1^k)^{m^*-2})(1+|y|)^{2\delta}}\Bigg|_{\partial B_{\tau/\varepsilon_k}}\to0.\nonumber
\end{align}
Therefore, $r_k$ is contained in $B_\frac{\tau}{\varepsilon_1^k}$. Let 
\[\hat w_i^k(y)=\frac{w_i^k(y)}{\Lambda_k((\varepsilon_1^k)^{\delta}+(\varepsilon_1^k)^{m^*-2})(1+r_k)^{2\delta}},\quad |y|\le \tau/\varepsilon_k.\]

~

It is evident that
\begin{align}
\left\{
\begin{aligned}
&\Delta \hat w_i^k+\sum_j a_{ij}e^{\xi_j^k}\hat w_j^k=0\mbox{ in }B_\frac{C}{\varepsilon},\nonumber\\
&\frac{d}{dr}\hat w^k_i(0)=0\mbox{ for }i=1,\cdots,n,\nonumber
\end{aligned}
\right.
\end{align}
and
\[\exists i_0,\, \hat w_{i_0}^k(r_k)=\pm 1,\quad |\hat w_i(y)|\le \frac{(1+|y|)^{2\delta}}{(1+r_k)^{2\delta}}.\]
Here, we remark that $\xi_j^k\to  V_j$ over any compact subset of $\mathbb R^2$ , which is due to the difference of the equations of $V_{1,i}^k$ and $\overline{V}_{2,i}^k$.
Passing $k\to\infty$, we find that
$\hat w_i^k\to cZ_{0,i}$ in $C_{loc}^2(\mathbb{R}^2)$
for certain constant $c\in\mathbb{R}$. From the definition of $\eta_k$ we see that at least for one index $i_1$ we have $\hat w_{i_1}(0)=0$, thus $c=0$ and $\hat w_i=o(1)$ over any compact subset of $\mathbb R^2$. This rules out the possibility that $|r_k|<C$ for $C$ independent of $k$. But if $r_k\to \infty$, we use the Green's representation of $\hat w^k_{i_0}(r_k)$ to see that $\hat w^k_{i_0}(r_k)=o(1)$, a contradiction to $\hat w^k_{i_0}(r_k)=\pm 1$. 
The readers may read \cite[pp.402-403]{Zhang2009} for more details in this argument. (\ref{wes}) is established. Claim \ref{c:BubblesCompare} is established.

\noindent $\Box$

\medskip

From Claim \ref{c:BubblesCompare} we immediately obtain $|\alpha^k |=O(\varepsilon_k^{m^*-2})$ by evaluating $w_i^k(0)$.

Due to \cite{HuangZhang2022}, it is known that for any $t,s=1,\cdots,N$, $\varepsilon_t^k\sim\varepsilon_s^k$. Therefore, as we agree in above, we omit the footnotes and the superscripts in $\varepsilon_t^k$. To be precise, we denote $O(\varepsilon^\alpha)$ for $O((\varepsilon_t^k)^\alpha)$.

\begin{remark}
In order to obtain refined bound for $|\alpha^k|$, we use an analogue to (\ref{Expansion:BubblesComparison}) of $\varepsilon^{2m^*-4}$-order. Via a similar method as Claim \ref{c:BubblesCompare}, we prove that $|\alpha^k|=O(\varepsilon^{2m^*-4})$. Repeating this argument, we can prove that $|\alpha^k|=|\alpha_i^k-\overline{\alpha}_i^k|=O(\varepsilon^{M(m^*-2)})$ for arbitrary large integer $M$.
\end{remark}

By a classical result \cite[Lemma 3.2]{LinZhang2010}, we get
\begin{align}\label{ineq:eVDataK}
\sum_{i=1}^n\Big|\int_{\mathbb{R}^2} e^{V_{1,i}^k}-\int_{\mathbb{R}^2} e^{V_{2,i}^k}\Big|=O(\varepsilon^{M(m^*-2)-\delta})
\end{align}
and
\begin{align}\label{ineq:BubblesDataK}
|I_t^k-I_s^k|,\mbox{ }|V_{t,i}^k(0)-V_{s,i}^k(0)|=O(\varepsilon^{M(m^*-2)-\delta})
\end{align}
for any positive integer $M$.

For $t,s=1,\cdots,N$ with $t\neq s$, the scales $\varepsilon_t^k$ and $\varepsilon_s^k$ satisfy
\begin{align}\label{e:VarepsilonTVarespilonS}
\Bigg(\frac{\varepsilon^k_t}{\varepsilon_s^k}\Bigg)^{m_i^k-2}=\frac{e^{H_{i,t}(q^k)}}{e^{H_{i,s}(q^k)}} +O(\varepsilon^{m^*-2-\delta}).
\end{align}
Moreover, the points $p_t^k$ satisfy that
\begin{align}\tag{$SE_1$}\label{SE1}
\sum_{i=1}^n\Big[\nabla \ln h_i(q_{t}^{k})+2\pi m_i^* \nabla_1 G^*(q_t^k,q_t^k)\Big]\rho_i^*=\left\{
\begin{aligned}
&O(\varepsilon^{m^*-2}) &\mbox{ if } m^*<4;\\
&O\Big(\varepsilon^2\ln\frac{1}{\varepsilon}\Big)&\mbox{ if }m^*= 4.
\end{aligned}
\right.
\end{align}
and
\begin{align}\label{e:LHitHis}
\frac{H^k_{it}(q^k)-H^k_{is}(q^k)}{m_i^*-2}=\frac{H^k_{jt}(q^k)-H^k_{js}(q^k)}{m_j^*-2}+O(\varepsilon^{m^*-2-\delta})
\end{align}
with
$H^k_{it}$ defined in (\ref{def:H})
for any $i,j=1,\cdots,n$ and any $t,s=1,\cdots,N$ with $i\neq j$ and $t\neq s$. This is due to \cite[Theorem 2.6]{GuZhang2025}.

~~

\noindent{\bf Information outside of the bubbling disks}

~~

In follows, we list several information outside of the bubbling disks.
Recalling \cite[pp.20]{HuangZhang2022} and \cite[(3.37)]{HuangZhang2022}, we get
\begin{align}\label{e:ExpansionOutDisks}
u_i^k(x)=\int_M u_i^k+2\pi m_i^k \sum_{s=1}^N G(x,q_s^k)+O(\varepsilon^{m^*-2-\delta})
\end{align}
for $x\in M\backslash \cup_{t=1}^N B_\tau(q_t^k)$ and
\begin{align}\label{e:Averageuik}
\int_M u_i^k=(1-\frac{m_i^k}{2})M_{k,t}-\ln(\rho_i^k h_i^k(q_t^k))-2\pi m_i^k  G^*(q_t^k,q_t^k)+I_t^k+O(\varepsilon^{m^*-2-\delta})
\end{align}
for any $t=1,\cdots,N$. For the harmonic function defined as in (\ref{def:HarmonicFunction}), an immediate corollary is
\begin{align}\label{e:HarmonicLimit}
\phi_{t,i}^{0,k}(y)=2\pi m_i^*(G^*(y;q_t^*)-G^*(q_t^*;q_t^*))-f_i^k(y)+O(\varepsilon^{m^*-2-\delta})
\end{align}
for any $t=1,\cdots,N$, $i=1,\cdots,n$ and any $y\in B_{\tau}(q_t^k)$.

\subsection{Choosing a refined coordinate}\label{Subsection:RefinedCoordinateAndNotations}

In order to obtain a higher order expansion, we choose a suitable coordinate to make the boundary oscillation as small as possible. To do this, we proceed a reduction argument. Based on the results in the last subsection, we have the following preliminary estimate.
First, we locate the point $p_{1,t}^{k}\in B_{\tau_0/\varepsilon_1^k}(0)$.

\begin{lemma}\label{l:pik}
For any $i=1,\cdots,n$, $p_{1,t}^{k}=-\Big(\frac{2\partial_1\phi_{1,t}^{0,k}(0)}{\sum_{j=1}^n a_{1j} h_j^k(0) e^{V^k_j(0)}},\frac{2\partial_1\phi_{1,t}^{0,k}(0)}{\sum_{j=1}^n a_{1j} h_j^k(0) e^{V^k_j(0)}}\Big)\varepsilon_1^k+ O(\varepsilon^2)\in B_{\tau_0/\varepsilon_1^k}(0)$.
\end{lemma}
\noindent{\bf Proof.}
By \cite{LinZhang2013}, we know that $p_{1,t}^{k}=O(\varepsilon_k)$. Hence, we only need to search for $p_{1,t}^{k}$ in $B_0=\{y\in\mathbb{R}^2||y|\leq\varepsilon_1^k\ln\frac{1}{\varepsilon_1^k}\}$. Denoting $w_1^{0,k}(y)=v_1^{0,k}(y)-V_{1,t}^{0,k}(y)-\phi_{1,t}^{0,k}(\varepsilon_1^k y)$, we get
\begin{align}\label{e:pikmax}
0=\nabla_y(v_1^{0,k}(y)-\phi_{1,t}^{0,k}(\varepsilon_1^k y))|_{y=p^{k}_{1,t}}=\nabla_y(V_{1,t}^{0,k}(y)+w_1^{0,k}(y))|_{y=p_{1,t}^{k}}.
\end{align}
By \cite[Theorem 3.1]{LinZhang2013} and that $\nabla v_i(0)=0$, we know that
\begin{align}
\nabla_y w_1^{0,k}(y)=-\varepsilon_1^k\nabla\phi_{i,t}^{0,k}(0)+\varepsilon_1^kO(|y|)\,\mbox{ for }\,|y|\leq\varepsilon_1^k \ln\frac{1}{\varepsilon_1^k}.\nonumber
\end{align}
This is due to Taylor expansion of $\nabla w_1^{1,k}$ at $y=0$. Plugging this into (\ref{e:pikmax}), we get
\begin{align}
-\nabla_y w_1^{1,k}(p_{1,t}^{k})=\varepsilon_k\nabla\phi_{1,t}^{0,k}(0) +O(\varepsilon_k^2\ln\frac{1}{\varepsilon_k})=-\frac{1}{2}\sum_{j=1}^n a_{ij}h_j^k(0)e^{V_{j,t}^{0,k}(0)} p_{1,t}^{k}+O(|p_{1,t}^{k}|^3)=\nabla_yV_{1,t}^{0,k}(p_{i,k}).\nonumber
\end{align}
Here, we use (\ref{e:LimitSystemAsy1}). A direct computation implies that
\begin{align}
p_{1,t}^k=-\Big(\frac{2\partial_1\phi_{1,t}^{0,k}(0)}{\sum_{j=1}^n a_{ij} h_j^k(0) e^{V^{0,k}_{j,1}(0)}},\frac{2\partial_1\phi_{1,t}^{0,k}(0)}{\sum_{j=1}^n a_{ij} h_j^k(0) e^{V^{0,k}_{j,1}(0)}}\Big)\varepsilon_1^k+ O\Big(\varepsilon^2\ln\frac{1}{\varepsilon}\Big).\nonumber
\end{align}

Repeating the argument with $B_0$ changed into \[B'_0=\Bigg\{y\in\mathbb{R}^2\Bigg||y|\leq100\sum_{i=1}^n\Bigg(
\frac{2|\partial_1\phi_{1,t}^{0,k}(0)|}{\sum_{j=1}^n a_{ij} h_j^k(0) e^{V^{0,k}_{j,1}(0)}}+\frac{2|\partial_2\phi_{1,t}^{0,k}(0)|}{\sum_{j=1}^n a_{ij} h_j^k(0) e^{V^{0,k}_{j,1}(0)}}+1\Bigg)\varepsilon_1^k\Bigg\},\] we prove the result.
$\Box$

Denote $\Omega_{1,k}=B_{\tau_1^k/\varepsilon_1^k}(p_{1,t}^k)$. Then, we get
\begin{proposition}\label{prop:Oscv1}
It holds that $\sup_{y_1,y_2\in\partial\Omega_{1,k}}|v_i^{1,k}(y_1)-v_i^{1,k}(y_2)|=O(\varepsilon^2)$.
\end{proposition}
\noindent{\bf Proof.}
We use the Green's representation.
Let $G_k(y,\eta)$ be the Green's function on $\Omega_{1,k}$, i.e.
$G_k(y,\eta)=\frac{1}{2\pi}\ln\frac{|y|}{(\varepsilon_1^k)^{-1}\tau_1^k}+\frac{1}{2\pi}\ln\frac{|y^*-\eta|}{|y-\eta|}$
with $y^*=\frac{(\varepsilon_1^k)^{-2}(\tau_1^k)^2 y}{|y|^2}$. By (\ref{e:Equationv1k}), we get
\begin{align}
v_i^{1,k}(y)=\int_{\Omega_{1,k}}G_k(y,\eta)\Big(\sum_{j=1}^n a_{ij}H_j^{1,k}(\eta)e^{v_i^{1,k}(\eta)}\Big)d\eta-\int_{\partial\Omega_{1,k}}\partial_\upsilon G_k v_i^{1,k}d\mathcal{H}^1_\eta.\nonumber
\end{align}
We omit the second term since it is a constant. We only need to deal with the integral $$\int_{\Omega_{1,k}}G_k(y,\eta)\Big(\sum_{j=1}^n a_{ij}H_j^{1,k}(\eta)e^{v_i^{1,k}(\eta)}\Big)d\eta.$$ To do this, we rewrite it into the form:
\begin{align}
\int_{\Omega_{1,k}}\frac{1}{2\pi}\ln\frac{|y|}{(\varepsilon_1^k)^{-1}\tau_1^k}\Big(\sum_{j=1}^n a_{ij}H_j^{1,k}(\eta)e^{v_i^{1,k}(\eta)}\Big)d\eta+\int_{\Omega_{1,k}}\frac{1}{2\pi}\ln\frac{|y^*-\eta|}{|y-\eta|}\Big(\sum_{j=1}^n a_{ij}H_j^{1,k}(\eta)e^{v_i^{1,k}(\eta)}\Big)d\eta \nonumber\\
=:I_1(y)+I_2(y).\nonumber
\end{align}
Notice that $y=p_{1,k}+(\varepsilon_1^k)^{-1}\tau_1^k(\cos\theta,\sin\theta)$. Then, we get
\begin{align}
\frac{1}{4\pi}\ln\frac{|y|^2}{(\tau_1^k)^2/(\varepsilon_1^k)^2}= \frac{1}{4\pi}\ln\frac{(\varepsilon_1^k)^2|p_{1,1,t}^k|^2+(\tau_1^k)^2+2\varepsilon_1^k\tau_k p_{1,1,t}^{k}\cos\theta+2\varepsilon_1^k\tau_1^k p_{2,1,t}^k\sin\theta}{\tau^2}.\nonumber
\end{align}
By this, we get
\begin{align}
I_1(y)=\frac{1}{2}\int_{\Omega_{1,k}}\frac{1}{4\pi}\ln\frac{(\varepsilon_1^k)^2|p_{1,k}|^2+ (\tau_1^k)^2+2\varepsilon_1^k\tau_1^k p_{1,1,t}^k\cos\theta+2\varepsilon_1^k\tau_1^k p_{2,1,t}^k\sin\theta}{\tau^2} \Big(\sum_{j=1}^n a_{ij}H_j^{1,k}(\eta)e^{v_i^{1,k}(\eta)}\Big)d\eta.\nonumber
\end{align}
Then, we get
\begin{align}\label{ineq:I1Osc}
\sup_{y_1,y_2\in\partial\Omega_{1,k}}|I_1(y_1)-I_1(y_2)|=O(\varepsilon^2).
\end{align}
In order to analyze $I_2(y)$, we cut the domain $\Omega_{1,k}$ into three parts: $B_{\tau_1^k/(2\varepsilon_1^k)}(y)$, $\Omega_{1,k}\backslash (B_{\tau^k_1/(2\varepsilon_1^k)}{2}(y)\cup B_{\tau_1^k/(\varepsilon_1^k)^{1-\delta_0}}(0))$ and $B_{\tau_1^k/(\varepsilon_1^k)^{1-\delta_0}}(0)$. Here, $y\in\partial\Omega_{1,k}$ and $\delta_0$ is a small positive constant.

As for the part of $I_2(y)$ on $B_{\tau_1^k/(2\varepsilon_1^{k})}(y)$, we notice that
$\frac{1}{2\pi}\ln\frac{|y^*-\eta|}{|y-\eta|}\leq\frac{C\varepsilon_k}{|y^{**}-\eta|}$
for certain $y^{**}$ lies between $y$ and $y^*$. Since on this region, $\sum_j a_{ij}H_j^k(\eta)e^{v_i^{1,k}}=O(\varepsilon^{m^*})$. Therefore, on this region, the integral is $O(\varepsilon^{m^*})$.

As for the part of $II_2$ on $\Omega_{1,k}\backslash (B_{\tau_1^k/2\varepsilon_1^k}(y)\cup B_{\tau_1^k/(\varepsilon_1^k)^{1-\theta}}(0))$, we notice that $|y^{**}-\eta|\sim\varepsilon_k^{-1}$.
Then, the integral on this region is $O(\varepsilon^{(1-\delta_0)m^*+2\theta})$. For sufficiently small $\theta$, this is $O(\varepsilon^{m^*-\delta})$.

As for the part of $I_2(y)$ on $B_{\tau_1^k/(\varepsilon_1^k)^{1-\delta}}(0)$, we can see by a direct computation that
\begin{align}
\frac{1}{4\pi}\ln\frac{|y^*-\eta|^2}{|y-\eta|^2}&=\frac{1}{4\pi}\ln\Big(1+\frac{|y^*|^2-|y|^2+2(y^*-\eta)\cdot\eta}{|y-\eta|^2}\Big).\nonumber
\end{align}
Noticing that
$\frac{1}{|y-\eta|^2}-\frac{1}{|y|^2}=O(\varepsilon^3)|\eta|$,
we get
\begin{align}
\frac{1}{4\pi}\ln\frac{|y^*-\eta|^2}{|y-\eta|^2}&=\frac{1}{4\pi}\ln\Big(1+\frac{|y^*|^2-|y|^2}{|y|^2}+O(\varepsilon^3)|\eta|\Big)=\frac{1}{4\pi}\ln\Big(\frac{\tau^4\varepsilon_k^{-4}}{|y|^4}+O(\varepsilon^3)|\eta|\Big)\nonumber\\
&=-\frac{1}{4\pi}\ln\Big(\frac{|y|^4}{\tau_0^4(\varepsilon_1^k)^{-4}}+O(\varepsilon^3)|\eta|\Big).\nonumber
\end{align}
By $y=p_{1,t}^k\varepsilon_1^k+\tau_1^k(\cos\theta,\sin\theta)$, we get
\begin{align}\label{ineq:I2Osc}
\sup_{y_1,y_2\in\partial\Omega_{1,k}}|I_2(y_1)-I_2(y_2)|=O(\varepsilon^2).
\end{align}
(\ref{ineq:I1Osc}) and (\ref{ineq:I2Osc}) imply the result.
$\Box$

In follow, we iterate the above method to obtain a suitable coordinate with small oscillation on the boundary. 
Define the harmonic function to encode the oscillation of $v_i^{0,k}$ on $\partial B_{\tau'}$:
\begin{equation}\label{def:HarmonicFunction1}
    \left\{
   \begin{array}{lr}
     -\Delta \phi_{i,t}^{1,k}=0\mbox{ in }\Omega_{1,k},\\
     \phi_{i,t}^{1,k}|_{\partial \Omega_{1,k}}=v_i^{1,k} -\frac{1}{2\pi\tau_1^k}\int_{\partial \Omega_{1,k}}v_i^{1,k} \mbox{ for }i\in I.
   \end{array}
   \right.
\end{equation}
Let $p_{2,t}^k\in B_{\tau_1^k/\varepsilon_t^k}(p_{1,t}^k)$ be the maximum point of the function $v_i^{1,k}-\phi_{i,t}^{1,k}$. Then, $v_1^{2,k}:=v_i^{1,k}-\phi_{i,t}^{1,k}$ satisfies that
\begin{align}\label{e:Equationv2k}
-\Delta v_i^{2,k}=\sum_{j=1}^n a_{ij} H^{2,k}_j( y)e^{v_j^{2,k}}\mbox{ in }B_{\tau_2^k}(0)
\end{align}
with $H_i^{2,k}(\cdot)=\frac{h^k_i(\varepsilon_t^k\cdot +q_{2,t}^{k})}{h_i^k(q_{2,t}^{k})}e^{\psi+f_i^k+\phi_{t,i}^{0,k}(\varepsilon_t^k\cdot+q_{2,t}^k) +\phi_{t,i}^{1,k}(\varepsilon_t^k\cdot+q_{2,t}^k)}$ and $\tau_2^k=\tau_1^k-\varepsilon_1^k|p_{2,k}^k|$. Denote $\Omega_{2,k}=B_{\tau_2^k/\varepsilon_t^{k}}(p_{2,t}^k)$. By a direct computation and the method as in Lemma \ref{l:pik} and Proposition \ref{prop:Oscv1}, we get
\begin{align}
|p_{2,t}^k-p_{1,t}^k|=O(\varepsilon^3)\mbox{ and }\sup_{y_1,y_2\in\partial\Omega_{2,k}}|v_i^{2,k}(y_1)-v_i^{2,k}(y_2)|=O(\varepsilon^4).\nonumber
\end{align}

In general, if we iterate the above procedure $M$ times, we get
\begin{itemize}
  \item a sequence $p_{1,t}^k,\,p_{2,t}^k,\cdots,p_{M,t}^k$ with $|p_{l+1,t}^k -p_{l,t}^k|=O(\varepsilon^{2l+1})$ for $l=1,\cdots,M-1$;
  \item $\tau_{l+1}^k=\tau_l^k-\varepsilon_t|p_{l+1,t}^k|$ and $\Omega_{l,k}=B_{\tau_l^k/\varepsilon_t^{k}}(p_l^k)$. Notice that $\Omega_{l+1,k}$ is tangent to $\Omega_{l,k}$;
  \item $v_i^{l,k}(y)=v_i^{0,k}(y-p_{l,1}^k)-\sum_{s=0}^1\phi_{i,t}^{s,k}$, $\phi_{i,t}^{s,k}(\varepsilon_t^k y)\leq O((\varepsilon_t^k)^{2l})|y|$ and $\sup_{y_1,y_2\in\partial\Omega_{l,k}}|v_i^{l,k}(y_1)-v_i^{l,k}(y_2)|=O(\varepsilon^{2l})$ for $l=1,\cdots,M$. Here, the harmonic functions $\phi_{i,t}^{s,k}$ are defined to cancel the boundary oscillation of $v_i^{s,k}$ on $\Omega_{s,k}$;
  \item the functions $v_i^{l,k}$ satisfy the equation
  \begin{align}\label{e:Equationvlk}
-\Delta v_i^{l,k}=\sum_{j=1}^n a_{ij} H^{l,k}_j( y)e^{v_j^{l,k}}\mbox{ in }\Omega_{l,k}
\end{align}
with $H_i^{l,k}(\cdot)=h^k_i(\varepsilon_t^k\cdot+ q_{l,t}^{k})e^{\psi+f_i^k+\sum_{s=0}^{l-1} \phi_{t,i}^{0,k}(\varepsilon_t^k\cdot+ q_{l,t}^k)}$;
\item For simplicity, we introduce the function $h_i^{l,k}(x)=H_i^{l,k}(x/\varepsilon_1^k)$ for $x\in\{\varepsilon_1^k y|y\in\Omega_{l,k}\}$.
\end{itemize}
In follows, we will use the functions $v_i^{M,k}$ for sufficiently large $M$.

\subsection{A refined expansion inside the bubbling disks}\label{subsection:BubblingInside}

In this subsection, we prove a refined expansion inside the bubbling disks as follows.
\begin{theorem}\label{t:IntegralExpansion}
Under the above assumptions, we get for $m^*<4$
\begin{align}
\int_{\Omega_{M,k}}H_i^{M,k}( y)e^{v_i^{M,k}}=\int_{\Omega_{M,k}}H_i^{M,k}(0)e^{V_j^{M,k}}+T^{in}_{M',k}+O(\varepsilon^{(M+2)(m^*-2)})\nonumber
\end{align}
and for $m^*=4$
\begin{align}
\int_{\Omega_{M,k}}H_i^{M,k}( y)e^{v_i^{M,k}}=\int_{\Omega_{M,k}}H_i^{M,k}(0)e^{V_i^{M,k}}+b_i^k(\varepsilon_1^k)+T^{in,*}_{M,k}+O(\varepsilon^{(M'+2)(m^*-2)}).\nonumber
\end{align}
Here, the terms 
\begin{align}\label{def:bik}
b_i^k(\varepsilon_1^k)=-\sum_j a_{ij}m_j^k H_j^{M,k}(0)\Big[\Delta\ln h_i^{M,k}(0)+|\nabla_x h_i^{M,k}(0)|^2\Big]\cdot(\varepsilon_1^k)^2\int_{\Omega_{M,k}}e^{V_j^{M,k}}Z_{0,j}|y|^2dy=O(\varepsilon^2\ln\frac{1}{\varepsilon}),
\end{align}
$T^{in}_{M,k},T^{in,*}_{M,k}=o_\tau(1)\varepsilon^{m^*-2}$ and is $C^2$ depends on $p_{i}^{M,k}$, $h_i$ and its derivatives up to $M+3$ order. The functions $H_i^{M,k}$ and $h_i^{M,k}$ are defined as in  Subsection \ref{Subsection:RefinedCoordinateAndNotations}.
\end{theorem}

\begin{remark}\label{r:LeadingTermL}
By an elementary observation, we get
\begin{align}
b_i^k(\varepsilon_1^k)=8\pi\sum_j a_{ij}L_{it}(\varepsilon_1^k)^2\ln\frac{1}{\varepsilon_1^k}+o(\varepsilon^{m^*-2}).\nonumber
\end{align}
Here, the number $L_{it}$ is defined as in (\ref{def:Lit}).
\end{remark}

\noindent{\bf Proof of Theorem \ref{t:IntegralExpansion}.}
We establish this expansion in several steps. 

~~

\noindent{\bf Step 1. First order terms}

~~

We begin by the following equation for the first order terms of the solution.
\begin{align}
\Delta c_{1,1,i,k}+\sum_j a_{ij}H_j^{M,k}(0)e^{V_{j}^{M,k}}c_{1,1,j,k}=-\varepsilon_1^k\sum_{j=1}^n (y\cdot\nabla h_j^{M,k}(0))e^{V_j^{M,k}}
\end{align}
with $c_{1,1,i,k}=0$ on $\partial\Omega_{M,k}$ for $i=1,\cdots,n$. Such existence is implied by Corollary \ref{coro:1-FreqSolvable}.


By this, we get
\begin{claim}\label{c:1FreqMainSize}
It holds that
\begin{equation}
|c_{1,1,i,k}(y)|=   
\left\{
\begin{array}{lr}
O(\varepsilon_1^k(1+|y|)^\delta)\mbox{ if }m^*>3;\nonumber\\
O(\varepsilon_1^k(1+|y|)^{3-m^*+\delta})\mbox{ if }m^*\le 3.\nonumber
\end{array}
\right.
\end{equation}
\end{claim}
\noindent{\bf Proof.}
We only prove the case $m^*\le 3$ since the other one is similar.
A first estimate claims that
\begin{equation}\label{ineq:c11Linfty}
\|c_{1,1,i,k}\|_{L^\infty(\Omega_{M,k})}=   
\left\{
\begin{array}{lr}
O(\varepsilon_1^k)\mbox{ if }m^*>3;\\
O((\varepsilon_1^k)^{m^*-2-\delta})\mbox{ if }m^*\le 3.
\end{array}
\right.
\end{equation}
This is due to Proposition \ref{prop:invertibility} and the direct computation
\begin{equation}
\Big\|\Big((y\cdot\nabla h_1^{M,k}(0))e^{V_1^{M,k}},\cdots,(y\cdot\nabla h_n^{M,k}(0))e^{V_n^{M,k}}\Big)\Big\|_{\mathbb{Y}_\varepsilon}=   
\left\{
\begin{array}{lr}
O(1)\mbox{ if }m^*> 3;\nonumber\\
O((\varepsilon_1^k)^{m^*-3-\delta})\mbox{ if }m^*\le 3\nonumber
\end{array}
\right.
\end{equation}
for sufficiently small positive constants $\beta$ and $\delta$.
Now we proceed a Green's function argument as in \cite[Theorem 3.1]{LinZhang2013} and \cite[Proposition 3.1]{BartYangZhang20241}.

We argue by contradiction. Assume that
\begin{align}
\sup_{y\in\Omega_{M,k}}\frac{|c_{1,1,i,k}(y)|}{\varepsilon_1^k(1+|y|)^{3-m^*+\delta}}=:\Lambda_k\to+\infty.\nonumber
\end{align}
We denote the point attending $\Lambda_k$ by $y_k\in\Omega_{M,k}$. Let us define
$$\overline{c}_i(y)=\frac{c_{1,1,i,k}(y)}{\Lambda_k\varepsilon_1^k(1+|y_k|)^{3-m^*+\delta}}$$
for $i=1,\cdots,n$.
By a direct computation, we get
\begin{align}
|\overline{c}_i(y)|\leq\Bigg|\frac{c_{1,1,i,k}(y)}{\Lambda_k\varepsilon_1^k(1+|y|)^{3-m^*-\delta}}\cdot\frac{(1+|y|)^{3-m^*+\delta}}{(1+|y_k|)^{3-m^*+\delta}}\Bigg|\leq \frac{(1+|y|)^{3-m^*+\delta}}{(1+|y_k|)^{3-m^*+\delta}}.\nonumber
\end{align}

~~

\noindent{\bf Case 1. $\{|y_k|\}_k$ is bounded.}

~~

Observe that in this case, we get
\begin{align}
\Delta\overline{c}_i+\sum_j a_{ij}H_j^{M,k}(0)e^{V_j^{M,k}}\overline{c}_j=o_k(1)\sum_ja_{ij}\frac{(y\cdot\nabla_x h_j^{M,k}(0)e^{V_j^{M,k}})}{(1+|y_k|)^{3-m^*+\delta}}\mbox{ in }\Omega_{M,k}.
\end{align}
On the other hand,
by (\ref{ineq:c11Linfty}) and the construction in Appendix \ref{App:Solvability}, we get
\begin{equation}
\sum_i\int_{B_{\tau_M^k/\varepsilon_1^k}}c_{1,1,i,k}\cdot e^{V_i^{M,k}}Z_{1,i}=  
\left\{
\begin{array}{lr}
O\Big(\varepsilon_1^k\int_{1/\varepsilon_1^k}^{3/\varepsilon_1^k}r^{-m^*}dr\Big)=o(\varepsilon_1^k)\mbox{ if }m^*\geq3;\nonumber\\
O\Big((\varepsilon_1^k)^{m^*-2-\delta} \int_{1/\varepsilon_1^k}^{3/\varepsilon_1^k}r^{-m^*}dr\Big)=o(\varepsilon_1^k)\mbox{ if }m^*< 3.\nonumber
\end{array}
\right.
\end{equation}
Then, we get a smooth function $\overline{c}_{\infty,i}$ such that
\begin{align}
\left\{
\begin{aligned}
&\overline{c}_{i}\to\overline{c}_{\infty,i}\mbox{ in }C^2_{loc}(\mathbb{R}^2)\mbox{ and }|\overline{c}_i(y)|\leq C(1+|y|)^\tau\mbox{ for }\tau\in(0,1),\nonumber\\
&\Delta\overline{c}_{\infty,i}+\sum_j a_{ij}H_j^{M,\infty}(0)e^{V_j^{M,k}}\overline{c}_{\infty,j}=0\mbox{ in }\mathbb{R}^2,\nonumber\\
&\sum_i\int_{\mathbb{R}^2}\overline{c}_ie^{V_i^{M,\infty}}Z_{1,i}=\sum_i\int_{\mathbb{R}^2}\overline{c}_ie^{V_i^{M,\infty}}Z_{2,i}=\sum_i\int_{\mathbb{R}^2}\overline{c}_ie^{V_i^{M,\infty}}Z_{0,i}=0.\nonumber
\end{aligned}
\right.
\end{align}
where $Z_{0,i}$, $Z_{1,i}$ and $Z_{2,i}$ are defined in (\ref{kernels}).   This implies $\overline{c}_{\infty,i}=0$, which is a contradiction to its definition and Proposition \ref{p:Kernel}.

~~

\noindent{\bf Case 2. $\{|y_k|\}_k$ is unbounded and $|y_k|=o(1)\varepsilon_1^k$.}

~~

In this case, we use a Green's function argument to obtain the contradiction.
By {\bf Case 1}, we get $\overline{c}_i(0)\to0$ as $k\to\infty$. By Green's representation, we get
\begin{align}
\frac{1}{2}&\leq|\overline{c}_i(0)-\overline{c}_i(y_k)|\nonumber\\
&\leq\Bigg\{\int_{B_\frac{3\tau}{\varepsilon}}|G_k(y_k,\eta)-G_k(0,\eta)|\cdot\big|\sum_j a_{ij} e^{V_j^{M,k}}\overline{c}_j\big|d\eta\Bigg\}\nonumber\\
&\quad+o(1)\cdot\Bigg\{\int_{B_\frac{3\tau}{\varepsilon}}|G_k(y_k,\eta)-G_k(0,\eta)|\cdot\Bigg|\sum_j a_{ij}\frac{(\eta\cdot\nabla h_j^{M,k}(0)e^{V_j^{M,k}(\eta)})}{(1+|y_k|)^{3-m^*+\delta}}\Bigg|d\eta\Bigg\}=:A+B.\nonumber
\end{align}
Here, $G_k$ denotes the Green's function with Dirichlet boundary condition on $B_{\tau_k^M/\varepsilon_1^k}$. We will draw a contradiction by proving $A+B=o(1)$. First, we get
\begin{align}\label{ineq:GreenAPointwise}
\big|\sum_j a_{ij} e^{V_j^{M,k}}\overline{c}_j\big|\leq\frac{C}{(1+|y_k|)^{3-m^*-\delta}(1+|\eta|)^{2m^*-3-\delta}}
\end{align}
and
\begin{align}\label{ineq:GreenBPointwise}
\Bigg|\sum_j a_{ij}\frac{(\eta\cdot\nabla H_j^{M,k}(0)e^{V_j^{M,k}(\eta)})}{(1+|y_k|)^{3-m^*+\delta}}\Bigg|\leq \frac{C}{(1+|y_k|)^{3-m^*-\delta}(1+|\eta|)^{m^*-1}}.
\end{align}
By \cite[Theorem 3.2]{LinZhang2013}, we get
\begin{equation}\label{ineq:GreenCase1}
|G_k(y_k,
    \eta)-G_\varepsilon(0,\eta)| \leq \left\{
\begin{aligned}
& C(\ln|y_k| +|\ln|\eta||)&\mbox{ if }\eta\in\Sigma_1:=\{\eta\in B_{{\tau_k^M/\varepsilon_1^k}}(0)||\eta|<\frac{|y_k|}{2}\};\\
&C(\ln|y_k|+|\ln|y-\eta||)&\mbox{ if }\eta\in\Sigma_2:=\{\eta\in B_{{\tau_k^M/\varepsilon_1^k}}(0)||y_k-\eta|<\frac{|y_k|}{2}\};\\
&\frac{C|y_k|}{|\eta|}&\mbox{ if }\eta\in\Sigma_3:= B_{{\tau_k^M/\varepsilon_1^k}(0)}\backslash(\Sigma_1\cup\Sigma_2).
    \end{aligned}
    \right.
\end{equation}
Assume that $m^*<3$.
Combining (\ref{ineq:GreenAPointwise}), (\ref{ineq:GreenBPointwise}) and (\ref{ineq:GreenCase1}), we get
\begin{align}
A&\leq C\int_{\Sigma_1}\frac{\ln|y_k|+|\ln|\eta||}{(1+|y_k|)^{3-m^*}(1+|\eta|)^{2m^*-3}}d\eta+ C\int_{\Sigma_2}\frac{\ln|y_k|+|\ln|y-\eta||}{(1+|y_k|)^{3-m^*}(1+|\eta|)^{2m^*-3}}d\eta \nonumber\\
&\quad\quad+C\int_{\Sigma_3}\frac{|y_k|/|\eta|}{(1+|y_k|)^{3-m^*}(1+|\eta|)^{2m^*-3}}d\eta=o_k(1)\nonumber
\end{align}
and
\begin{align}
B&\leq C\int_{\Sigma_1}\frac{\ln|y_k|+|\ln|\eta||}{(1+|y_k|)^{3-m^*}(1+|\eta|)^{m^*-1}}d\eta+ C\int_{\Sigma_2}\frac{\ln|y_k|+|\ln|y-\eta||}{(1+|y_k|)^{3-m^*}(1+|\eta|)^{m^*-1}}d\eta \nonumber\\
&\quad\quad+C\int_{\Sigma_3}\frac{|y_k|/|\eta|}{(1+|y_k|)^{3-m^*}(1+|\eta|)^{m^*-1}}d\eta=o_k(1).\nonumber
\end{align}

~~

\noindent{\bf Case 3. $\{|y_k|\}_k$ is unbounded and $|y_k|\sim\varepsilon_1^k$.}

~~

This case is similar to {\bf Case 2}. The only difference is the estimate on the Green's function. In this case, we get
\begin{align}
|G_k(y_k,\eta)-G_k(0,\eta)|\leq C\Big(|\ln|\eta||+|\ln|y_k-\eta||+\ln\frac{1}{\varepsilon_1^k}\Big).\nonumber
\end{align}
\noindent{$\Box$}

Before we go to the next step, we remark that in some case it is convenient to write $c_{1,1,i,k}$ in polar coordinate:
\begin{align}
\left\{
\begin{aligned}
g_{1,1,i,k,1}(r)=\frac{1}{2\pi}\int_0^{2\pi} c_{1,1,i,k}(r\cos\theta,r\sin\theta)\cos\theta d\theta,\nonumber\\
g_{1,1,i,k,2}(r)=\frac{1}{2\pi}\int_0^{2\pi} c_{1,1,i,k}(r\cos\theta,r\sin\theta)\sin\theta d\theta.\nonumber
\end{aligned}
\right.
\end{align}
Then, the functions $g_{1,1,i,k,i}$ and $g_{1,1,i,k,2}$ satisfy that
\begin{equation}\label{e:1Feqg1}
    \left\{
   \begin{array}{lr}
    g_{1,1,i,k,1}''+\frac{1}{r} g_{1,1,i,k,1}'-\frac{1}{r^2} g_{1,1,i,k,1}+\sum_j a_{ij}H_j^{M,k}(0)e^{V_j^{M,k}}g_{1,1,i,k,1}=-\varepsilon_1^k\sum_ja_{ij}\partial_{x_1} h_j^{M,k}(0)re^{V_j^{M,k}},\\
    g_{1,1,i,k,1}(\tau^M_k/\varepsilon_1^{k})=0\mbox{ for }i\in I
   \end{array}
   \right.
\end{equation}
and
\begin{equation}\label{e:1Feqg2}
    \left\{
   \begin{array}{lr}
    g_{1,1,i,k,2}''+\frac{1}{r} g_{1,1,i,k,2}'-\frac{1}{r^2} g_{1,1,i,k,2}+\sum_j a_{ij}H_j^{M,k}(0)e^{V_j^{M,k}}g_{1,1,i,k,2}=-\varepsilon_1^k\sum_ja_{ij}\partial_{x_2} h_j^{M,k}(0)re^{V_j^{M,k}},\\
    g_{1,1,i,k,2}(\tau^M_k/\varepsilon_1^{k})=0\mbox{ for }i\in I.
   \end{array}
   \right.
\end{equation}

In follows we estimate the error part of the 1-frequency by using an ODE method.

\begin{proposition}\label{prop:ciError}
It holds that
\begin{align}
|w_{k,i}(y)-c_{1,1,k,i}(y)|=O((\varepsilon_1^k)^2(1+|y|)^{4-m^*+\delta})\nonumber
\end{align}
for $y\in\Omega_{M,k}$ and $\delta>0$ independent of $k$.  
\end{proposition}
\noindent{\bf Proof.}
First, we decompose $w_{k,i}$ by its frequency and denote
\begin{align}
w_{k,i}:=\sum_{l=0}^3 c_{l,k,i}\nonumber
\end{align}
with $c_{l,k,i}$ has its frequency equal to $l$ for $l=0,1,2$ while $c_{3,k,i}$ has the frequency no less than $3$. By classical result \cite{LinZhang2013}, we know that
\begin{align}
\Big|\sum_{l\neq 1}c_{l,k,i}(y)\Big|=O((\varepsilon_1^k)^2(1+|y|)^{4-m^*}).\nonumber
\end{align}
Now we estimate $c_{1,k,i}-c_{1,1,k,i}$ for $i=1,\cdots,n$.
By the construction of $w_{i,k}$, we know that
\begin{align}
\sup_{y_1,y_2\in\partial\Omega_{M,k}}|c_{1,k,i}(y_1)-c_{1,k,i}(y_2)|=O(\varepsilon^M)\nonumber
\end{align}
for $i=1,\cdots,n$ and $M$ sufficiently large. Taking a harmonic function $\phi_{1,i}$ encoding the boundary oscillation of $c_{1,k,i}$, i.e.
\begin{equation}
    \left\{
   \begin{array}{lr}
     -\Delta \phi_{1,i}=0\mbox{ in }\Omega_{M,k},\nonumber\\
     \phi_{1,i}|_{\partial \Omega_{1,k}}=c_{1,k,i} -\frac{1}{2\pi\tau_1^k}\int_{\partial \Omega_{1,k}}c_{1,k,i} \mbox{ for }i\in I.\nonumber
   \end{array}
   \right.
\end{equation}
An elementary observation gives $|\phi_{1,i}(y)|=O((\varepsilon_1^k)^{M+1}|y|)$ in $\Omega_{M,k}$. Denoting
\begin{align}
c_{1,i}^*(y)=c_{1,k,i}(y)-c_{1,1,k,i}(y)-\phi_{1,i}(y),\nonumber
\end{align}
we get
\begin{equation}
    \left\{
   \begin{array}{lr}
\Delta c_{1,i}^*+\sum_j a_{ij} H_j^{M,j}(0)e^{V_j^{M,k}}c_{1,j}^*=E_{1,i}^* \mbox{ in }\Omega_{M,k},\nonumber\\
c_{1,i}^*|_{\partial\Omega_{M,k}} =0 \mbox{ for }i\in I\nonumber
   \end{array}
   \right.
\end{equation}
with $E_{1,i}^*=O((\varepsilon_1^k)^2|y|^2(1+|y|)^{-m^*})$.
Using the idea in Remark \ref{r:MultiOrthognal} and a Green's function argument as in Claim \ref{c:1FreqMainSize}, we get 
\begin{align}
|c_{1,i}^*(y)|=O((\varepsilon_1^k)^2(1+|y|)^{4-m^*}).\nonumber
\end{align}
This completes the proof.
Notice that in the reduction scheme here, in order to make $(E_{1,1}^*,\cdots,E_{1,n}^*)\in F_\varepsilon^{odd}(T_1^\varepsilon,T_2^\varepsilon)$, we need to assume that
\begin{align}
&(y\cdot\nabla_x h_1^{M,k}(0)e^{V_1^{M,k}},\cdots,y\cdot\nabla_x h_n^{M,k}(0)e^{V_n^{M,k}}),\nonumber\\
&\quad(y\cdot\nabla_x h_1^{M,k}(0)e^{V_1^{M,k}}+E_{1,1}^*,\cdots,y\cdot\nabla_x h_n^{M,k}(0)e^{V_n^{M,k}}+E_{1,n}^*)\in F_\varepsilon^{odd}(T_1^\varepsilon,T_2^\varepsilon)\nonumber
\end{align}
for certain $T_1^\varepsilon,T_2^\varepsilon\in \mathbb{R}^2$. Here, the space $F_\varepsilon^{odd}(T_1^\varepsilon,T_2^\varepsilon)$ is defined in Remark \ref{r:MultiOrthognal}. $\Box$

~~

\noindent{\bf Step 2. Second order terms by Talyor expansion}

~~

Now we derive the second order terms via Taylor expansions.
We first write down the equation for $V_i^{M,k}+c_{1,1,i,k}$:
\begin{align}\label{ExpansionSecondOrderTotal}
\Delta(V_i^{M,k}+&c_{1,1,i,k})+\sum_j a_{ij}H_j^{M,k}( y)e^{V_j^{M,k}+c_{1,1,j,k}}\nonumber\\
&=\Delta(V_i^{M,k}+c_{1,1,i,k})+\sum_j a_{ij}H_j^{M,k}(0)e^{V_j^{M,k}+c_{1,1,j,k}+\ln\frac{H_j^{M,k}( y)}{H_j^{M,k}(0)}}\nonumber\\
&=\Delta(V_i^{M,k}+c_{1,1,i,k})+\sum_j a_{ij}H_j^{M,k}(0)e^{V_j^{M,k}}\Bigg(1+c_{1,1,j,k}+\ln\frac{H_j^{M,k}(y)}{H_j^{M,k}(0)}\nonumber\\
&\quad+\frac{1}{2}\Bigg(c_{1,1,j,k}+\ln\frac{H_j^{M,k}( y)}{H_j^{M,k}(0)}\Bigg)^2+O((\varepsilon_1^k)^3 |y|^3)\Bigg).
\end{align}
Here, in the last step we use Claim \ref{c:1FreqMainSize} and a standard regularity argument.
Now we write the Taylor expansion of the logarithmic term in polar coordinate:
\begin{align}
\ln\frac{H_j^{M,k}(\varepsilon_1^k y)}{H_j^{M,k}(0)}=\varepsilon_1^k\nabla_x\ln h_j^{M,k}(0)\cdot(\theta r)+(\varepsilon_1^k)^2_k(r^2\Theta_{2,j}+\mathcal{R}_{2,j})+O((\varepsilon_1^k)^3 |y|^3)\nonumber 
\end{align}
with
\begin{align}
\Theta_{2,j}
&=\frac{1}{4}\Bigg(\cos2\theta\Big(\partial_{x_1}^2\ln h_j^{M,k}(0)-\partial_2^2\ln h_j^{M,k}(0)\Big)+2\sin2\theta \cdot\partial_{x_1}\partial_{x_2}\ln h_j^{M,k}(0)\Bigg)\nonumber
\end{align}
and
$\mathcal{R}_{2,j}=\frac{1}{4}r^2\Delta\ln h_j^{M,k}(0)$.
Therefore, we get
\begin{align}
c_{1,1,i,k}+\ln\frac{H_j^{M,k}(\varepsilon_1^ky)}{H_j^{M,k}(0)}=c_{1,1,i,k}+\varepsilon_1^k\nabla_x\ln h_j^{M,k}(0)\cdot(\theta r)+(\varepsilon_1^k)^2(|y|^2\Theta_{2,j}+\mathcal{R}_{2,j})+O((\varepsilon_1^k)^3 |y|^3).\nonumber
\end{align}
Plugging these back into (\ref{ExpansionSecondOrderTotal}), we get
\begin{align}\label{Expansion:SecondOrderFinal}
\Delta(V_i^{M,k}&+c_{1,1,i,k})+\sum_j a_{ij}H_j^{M,k}( y)e^{V_j^{M,k}+c_{1,1,j,k}}\nonumber\\
&=\frac{1}{4}\sum_j a_{ij}H_j^{M,k}(0)e^{V_j^{M,k}}\Bigg((\varepsilon_1^k)^2 |y|^2\Delta_x\ln h_j^{M,k}(0)+\Big|(g_{1,1,i,k,1},g_{1,1,i,k,2})+\varepsilon_1^k |y|\nabla_x h_j^{M,k}(0)\Big|^2\Bigg)\nonumber\\
&\quad+\sum_j a_{ij}H_j^{M,k}(0)e^{V_j^{M,k}}\hat{\Theta}_{2,j}
\end{align}
with
\begin{align}\label{def:Theta2}
\hat{\Theta}_{2,j}&=\frac{1}{4}\Bigg[(\varepsilon_1^k)^2\partial_{x_1}^2\ln h_j^{M,k}(0)-(\varepsilon_1^k)^2\partial_{x_2}^2\ln h_j^{M,k}(0)+g_{1,1,i,k,1}^2-g_{1,1,i,k,2}^2 \nonumber\\
&\quad\quad+\varepsilon_1^k\Big(g_{1,1,i,k,1}\partial_1 H_j^{M,k}(0)-g_{1,1,i,k,2}\partial_2 H_j^{M,k}(0)\Big)\Bigg]\cos2
\theta\nonumber\\
&+\frac{1}{2}\Bigg[(\varepsilon_1^k)^2\partial_{x_1}\partial_{x_2}\ln h_j^{M,k}(0)+g_{1,1,i,k,1}g_{1,1,i,k,2}+\varepsilon_1^k\Big(g_{1,1,i,k,1}\partial_{x_2} h_j^{M,k}(0)+g_{1,1,i,k,2}\partial_{x_1} h_j^{M,k}(0)\Big)|y|\nonumber\\
&\quad\quad+(\varepsilon_1^k)^2\partial_{x_1} h_j^{M,k}(0)\cdot \partial_{x_2} h_j^{M,k}(0)|y|^2\Bigg]\sin2\theta.
\end{align}
Here, the functions $g_{1,1,i,k,1}$ and $g_{1,1,i,k,2}$ are defined as in (\ref{e:1Feqg1}) and (\ref{e:1Feqg2}).
Therefore, it is reasonable to 
define the second order terms of the error:
\begin{equation}\label{e:0FeqMain}
\left\{
\begin{array}{lr}
    c_{2,0,i,k}''+\frac{1}{r} c_{2,0,i,k}'+\sum_j a_{ij}H_j^{M,k}(0)e^{V_j^{M,k}}c_{2,0,i,k}\\
    \quad\quad\quad=\frac{1}{4}\sum_j a_{ij}H_j^{M,k}(0)e^{V_j^{M,k}}\Bigg((\varepsilon_1^k)^2 r^2\Delta_x\ln h_j^{M,k}(0)+\Big|(g_{1,1,i,k,1},g_{1,1,i,k,2})+\varepsilon_1^k r\nabla_x h_j^{M,k}(0)\Big|^2\Bigg),\\
    c_{2,0,i,k}(0)=c_{2,0,i,k}'(0)=0\mbox{ for }i\in I
\end{array}
\right.
\end{equation}
and
\begin{equation}\label{e:2FeqMain}
\left\{
\begin{array}{lr}
    \Delta c_{2,2,i,k}+\sum_j a_{ij}H_j^{M,k}(0)e^{V_j^{M,k}}c_{2,0,i,k}=\sum_j a_{ij}H_j^{M,k}(0)e^{V_j^{M,k}} \hat{\Theta}_{2,j},\\
    c_{2,2,i,k}(0)=|\nabla c_{2,2,i,k}(0)|=0\mbox{ for }i\in I.
\end{array}
\right.
\end{equation}
Here, the terms $\hat{\Theta}_{2,j}$ is defined as in (\ref{def:Theta2}). An ODE method and a Green's function argument implies that
\begin{align}
|c_{2,0,i,k}(y)|=O((\varepsilon_1^k)^2(1+|y|)^{4-m^*{+\delta}})\,\mbox{ for }\,0<|y|<\tau_M^k(\varepsilon_1^k)^{-1}\nonumber
\end{align}
and
\begin{align}
|c_{2,2,i,k}(y)|=O((\varepsilon_1^k)^2(1+|y|)^{4-m^*{+\delta}})\,\mbox{ for }\,y\in\Omega_{M,k}\nonumber
\end{align}
with $c_{2,0,i,k}|_{\partial\Omega_{M,k}}=0$ and $c_{2,2,i,k}|_{\partial\Omega_{M,k}}=0$.

Before we go to the next stage, we have the following result concerning $c_{2,0,k,i}$.
As similar computation as in Proposition \ref{prop:ciError} gives the following result.
\begin{claim}
It holds that
\begin{align}
|w_{i,k}(y)-c_{1,1,i,k}(y)-c_{2,0,i,k}(y)-c_{2,2,i,k}(y)|=O((\varepsilon_1^k)^2(1+|y|)^{4-m^*{+\delta}})\nonumber
\end{align}
for $y\in\Omega_{M,k}$.
\end{claim}

Moreover, for $c_{2,0,i,k}$, we get
\begin{claim}
Assume that $m^*<4$. Then, $\int_{\partial\Omega_{M,k}}\frac{\partial c_{2,0,i,k}}{\partial{\bf n}}=O(\tau^{4-m^*}\varepsilon^{m^*-2})$.
\end{claim}
This follows a direct computation. On the other hand, we get
\begin{claim}
Assume that $m^*=4$. Then,
\begin{align}
\int_{\partial\Omega_{M,k}}\frac{\partial c_{2,0,i,k}}{\partial{\bf n}}=b_i^k(\varepsilon_1^k)+o_\tau(1){\varepsilon^2}\nonumber
\end{align}
with $$b_i^k(\varepsilon_1^k)=-\sum_j a_{ij}m_j^k H_j^{M,k}(0){\Big[\Delta\ln h_i^{M,k}(0)+|\nabla_x h_i^{M,k}(0)|^2\Big]}\cdot(\varepsilon_1^k)^2\int_{\Omega_{M,k}}e^{V_j^{M,k}}Z_{0,j}|y|^2dy=O(\varepsilon^2\ln\frac{1}{\varepsilon}).$$
\end{claim}
\noindent{\bf Proof.}
Recall that
\begin{align}
\Delta\Big(\sum_j a^{ij}Z_{0,j}\Big)+H_i^{M,k}(0)e^{V_i^{M,k}}Z_{0,i}=0\mbox{ for }i=1,\cdots,n.\nonumber
\end{align}
Combining (\ref{e:0FeqMain}) and divergence theorem, we get
\begin{align}
&\sum_j a^{ij}\int_{\partial\Omega_{M,k}}\bigg (\frac{\partial c_{0,j}}{\partial{\bf n}}Z_{i,0}-\frac{\partial Z_{0,j}}{\partial{\bf n}}c_{i,0}\bigg )=\sum_j a^{ij}\int_{\Omega_{M,k}}(\Delta c_{0,j} Z_{i,0}-\Delta Z_{0,j} c_{0,i})=\int_{\Omega_{M,k}}E_{0,i}Z_{0,i}\nonumber\\
&=\frac{1}{4} H_i^{M,k}(0)\int_{\Omega_{M,k}}e^{V_i^{M,k}}Z_{0,i}\Bigg[(\varepsilon_1^k)^2 |y|^2{\Delta_x\ln h_i^{M,k}(0)}+\Big|(g_{1,1,i,k,1},g_{1,1,i,k,2})+\varepsilon_1^k y{\nabla_x h_i^{M,k}(0)}\Big|^2\Bigg].\nonumber
\end{align}
Here, $E_{0,i}$ denotes the right hand side of (\ref{e:0FeqMain}).
By a direct computation, we get
\begin{align}
\int_{\partial\Omega_{M,k}}\frac{\partial Z_{0,j}}{\partial{\bf n}}c_{i,0}=O(\varepsilon^{2+\delta})\nonumber
\end{align}
and
\begin{align}
\int_{\partial\Omega_{M,k}}\frac{\partial c_{0,j}}{\partial{\bf n}}Z_{i,0}=-m_i^k\int_{\partial\Omega_{M,k}} \frac{\partial c_{0,j}}{\partial{\bf n}}+O(\varepsilon^{2+\delta}).\nonumber
\end{align}
On the other hand,
\begin{align}
&\int_{\Omega_{M,k}}e^{V_i^{M,k}}Z_{0,i}\Big|(g_{1,1,i,k,1},g_{1,1,i,k,2})+\varepsilon_1^k |y|{\nabla_x h_i^{M,k}(0)}\Big|^2dy\nonumber\\
&=o_\tau(1)\varepsilon^2+(\varepsilon_1^k)^2{\partial_{x_1}h_i^{M,k}(0)}\int_0^{\tau_k^M/\varepsilon_1^k}r^2 e^{V_i^{M,k}}Z_{0,i}dr.\nonumber
\end{align}
Notice that in the last step, we use integration by parts. Moreover, we get
\begin{align}
(\varepsilon_1^k)^2\int_{\Omega_{M,k}}e^{V_i^{M,k}}Z_{0,i}|y|^2dy\leq 4\pi(\varepsilon_1^k)^2\int_0^{\tau_k^M/\varepsilon_1^k}e^{V_j}r^3dr=O\Big(\varepsilon^2\ln\frac{1}{\varepsilon}\Big).
\end{align}
Therefore, we get
\begin{align}
\int_{\partial\Omega_{M,k}}\frac{\partial c_{0,i}}{\partial{\bf n}}=-\sum_j a_{ij}m_j^k H_j^{M,k}(0){\Big[\Delta_x\ln h_i^{M,k}(0)+|\nabla_x h_i^{M,k}(0)|^2\Big]}\cdot(\varepsilon_1^k)^2\int_{\Omega_{M,k}}e^{V_j^{M,k}}Z_{0,j}|y|^2dy+o_\tau(1)\varepsilon^2.\nonumber
\end{align}
$\Box$

~~

\noindent{\bf Step 3. Higher order expansion: a general result}

~~

Next we obtain a $M$-th order approximation recursively. Denoting this term by $c_{M,i}^k$. It is known that
\begin{align}\label{ineq:cMki}
|c_{M,i}^k(y)|=O((\varepsilon_1^k)^{M}(1+|y|)^{M+\delta})\,\mbox{ in }\,\Omega_{M,k}
\end{align}
for $i=1,\cdots,n$. By a direct computation, we get
\begin{align}
\Delta(V_i^{M,k}+c_{M,i}^k)+\sum_j a_{ij}H_j^{M,k}( y)e^{V_i^{M,k}+c_{M,i}^k}=:A_{M,i}+B_{M,i}\nonumber
\end{align}
for $i=1,\cdots,n$. Here,
\begin{align}
A_{M,i}(y)=O((\varepsilon_1^k)^{M+1}|y|^{M'+1}(1+|y|)^{-m^*{+\delta}})\,\,\mbox{and}\,\,
B_{M,i}(y)=O((\varepsilon_1^k)^{M}(1+|y|)^{3M-(M+1)m^*{+\delta}})\nonumber
\end{align}
for $y\in\Omega_{M,k}$.
We define the following two 
auxiliary functions:
\begin{equation}\label{e:AuxFun12}
\left\{
\begin{array}{lr}
    \Delta \phi_{1,i}^{**,k}=0\mbox{ in }\Omega_{M,k},\\
    \phi_{1,i}^{*,k}|_{\partial\Omega_{M,k}}=-w_{i,k}+\frac{1}{2\pi\tau_k^M/\varepsilon_1^k}\int_{\partial\Omega_{M,k}}w_{i,k}\mbox{ for }i=1,\cdots,n
\end{array}
\right.
\end{equation}
and
\begin{equation}\label{e:AuxFun22}
\left\{
\begin{array}{lr}
\Delta \phi_{2,i}^{**,k}+\sum_j a_{ij}H_j^{M,k}(0)e^{V_j^{M,k}}\phi_{2,j}^{**,k}=A_{M',i}\mbox{ in }\Omega_{M,k},\\
\phi_{2,i}^{**,k}(0)=0,\,\phi_{2,i}^{**,k}|_{\partial\Omega_{M,k}}=\mbox{constant}\mbox{ for }i=1,\cdots,n.
\end{array}
\right.
\end{equation}
It is worth to be point out that the functions $c_{M',i}^k$, $\phi_{1,i}^{**,k}$ and $\phi_{2,i}^{**,k}$ are $C^1$ depends on $p_{i}^{M,k}$, $h_i$ and its derivatives up to $M'+2$ order. Defining the error
\begin{align}\label{Expansion:viMk}
W_{i,k}^*=v_i^{M,k}-V_i^{M,k}-c_{M,i}-\phi_{1,i}^{**,k}-\phi_{2,i}^{**,k},
\end{align}
we get
\begin{equation}\label{equ:W*-Final}
\left\{
\begin{array}{lr}
\Delta W^*_{i,k}+\sum_j a_{ij} H_j^{M,k}( y)e^{\overline{V_i^{M,k}}}W^*_{j,k}=O((\varepsilon_1^k)^{M+1}(1+|y|)^{3M-(1+M)m^*{+\delta}})\mbox{ in }\Omega_{M,k},\\
W^*_{i,k}|_{\partial\Omega_{M,k}}=\mbox{constant}\mbox{ for }i=1,\cdots,n.
\end{array}
\right.
\end{equation}
Here, $\overline{V_i^{M,k}}$ is due to Mean Value Theorem.
Denoting $(c_{M,i}+\phi_{1,i}^{**}-\phi_{2,i}^{**})_0-c_{2,0,i,k}=c_{*,0,i,k}$. Here, $(u)_0$ represents the radial part of $u$.
A routine estimate implies that
\begin{align}
|W_{i,k}^*(y)|=O((\varepsilon_1^k)^{M+1}(1+|y|)^{2+3M-(1+M)m^*{+\delta}})\nonumber
\end{align}
for $y\in\Omega_{M,k}$.
By a direct computation, we get for the case $m^*<4$
\begin{align}
\int_{\Omega_{M,k}}H_j^{M,k}( y)e^{v_j^{M,k}}&=\int_{\Omega_{M,k}}H_j^{M,k}(0)e^{V_j^{M,k}}+\int_{\partial\Omega_{M,k}}\frac{\partial c_{2,0,i,k}}{\partial{\bf n}}+\int_{\partial\Omega_{M,k}}\frac{\partial c_{*,0,i,k}}{\partial{\bf n}}+O(\varepsilon_k^{(M+2)(m^*-2){-\delta}})\nonumber\\
&=\int_{\Omega_{M,k}}H_j^{M,k}(0)e^{V_j^{M,k}}+T_{M,k}^{in}+O(\varepsilon_k^{(M+2)(m^*-2){-\delta}}).\nonumber
\end{align}
Here, the term $\,T^{in}_{M,k}=o_\tau(1)\varepsilon^{m^*-2}$ and is $C^2$ depends on $p_{i}^{M,k}$, $h_i$ and its derivatives up to $M'+3$ order.
As similar argument can be proceed for the case $m^*=4$.
The proof is completed.

\noindent{$\Box$}

\subsection{A refined estimate outside the bubbling disks}

In order to proceed with into the next step, we need to refine the estimate outside of the bubbling disks. A rough estimate is 
\begin{align}\label{e:ExpansionOutDisksROUGH}
u_i^k(x)=\int_M u_i^k+2\pi m_i^k \sum_{t=1}^N G(x;q^{M,k}_t)+O(\varepsilon_{k}^{m^*-2-\delta})
\end{align}
for $x\in M\backslash \cup_{t=1}^N B_{\tau_k^M}(q_t^{M,k})$ and
\begin{align}\label{e:AverageuikROUGH}
\int_M u_i^k=(1-\frac{m_i^k}{2})M_{k,t}-\ln(\rho_i^k h_i^k(q_t^k))-2\pi m_i^k G^*(q_t^{M,k};q_t^{M,k})+I_t^k+O(\varepsilon^{m^*-2-\delta})
\end{align}
for any $t=1,\cdots,N$.
Here, {$q_t^{M,k}$ are the coordinate in $M$ of $p_t^{M,k}$, $q^{M,k}=(q_1^{M,k},\cdots,q_N^{M,k})$ and }
\begin{align}
m_i^k=\frac{1}{2\pi}\sum_{j=1}^n a_{ij}\int_{\Omega_{t,M,k}}h_i^k e^{u_i^k}d V_g.\nonumber
\end{align}
For the harmonic function defined as in (\ref{def:HarmonicFunction}), an immediate corollary is
\begin{align}\label{e:HarmonicLimitROUGH}
\phi_{t,i}^{0,k}(y)=2\pi m_i^*(G^*(y;q_t^*)-G^*(q_t^*;q_t^*))-f_i^k(y)+O(\varepsilon^{m^*-2-\delta})
\end{align}
for any $t=1,\cdots,N$, $i=1,\cdots,n$ and any $y\in B_{\tau_k^M}(q_t^{M,k})$.
This is an immediate corollary of (\ref{e:ExpansionOutDisks}), (\ref{e:Averageuik}), (\ref{e:HarmonicLimit}) and the shift due to  $p_t^{M,k}$.
Our refined estimates reads as follows.

\begin{proposition}\label{prop:OutAverageExpansion}
For any large integer $M$, it holds that
\begin{align}
 u_i^k(x)&=\Bigg(1-\frac{m_i^k}{2}\Bigg)M_{k,t}-\ln(\rho_i^k h_i^k(q_t^k))-2\pi m_i^k G^*(q_t^{M,k};q_t^{M,k})+I_t^k+T_{M,k}^{out}+O(\varepsilon^{M(m^*-2) {-\delta}})\nonumber
\end{align}
and
\begin{align}
\int_M u_i^k=-\frac{m_{itv}^k-2}{2}\ln\frac{1}{\varepsilon_t^k}-2\pi m_{itv}^k G^*(q_t^{M,k};q_t^{M,k})+I_t^k-\ln(\rho_i^kh_i^k(q_t^{M,k}))+T_{M,k}^{avg}+O(\varepsilon^{M(m^*-2){-\delta}}).\nonumber
\end{align}
Here, $T_{M,k}^{out}$ and $T_{M,k}^{avg}$ 
depend $C^2$-continuously in $p_{t}^{M,k}$, $\varepsilon_t$ and the derivatives of $h_i^k$ up to $M+3$-th order.
\end{proposition}
\noindent{\bf Proof.}
We start by a Green representation.
\begin{align}
    u_i^k(x)&=\int_M u_i^k+\int_M G(x,\eta)\sum_{j=1}^n a_{ij}\rho_j^k h_j^k e^{u_j^k}d V_g(\eta)\nonumber\\
    &=\int_M u_i^k+\sum_{t=1}^N\int_{B_{\tau_k^M}(q_t^{M,k})}+\int_{M\backslash(\cup_t B_{\tau_k^M}(q_t^{M,k}))} G(x,\eta)\sum_{j=1}^n a_{ij}\rho_j^k h_j^k e^{u_j^k}d V_g(\eta).\nonumber
\end{align}
Here,
for any $t=1,\cdots,N$, we get
\begin{align}\label{Expansion:GreenInside}
\int_{B_{\tau_k^M}(q_t^{M,k})}&G(x,\eta)\sum_{j=1}^n a_{ij}\rho_j^k h_j^k e^{u_j^k}dV_g(\eta)\nonumber\\
&=\sum_{j=1}^n a_{ij}\rho_j^k\int_{B_{\tau_k^M}(q_t^{M,k})}G(x,\eta) H_j^{M,k}(\eta)e^{V_j^{M,k}+w_j^{M,k}(\eta)}d V_g(\eta)\nonumber\\
&=\sum_{j=1}^n a_{ij} \rho_j^k\int_{B_{\tau_k^M}(q_t^{M,k})}G(x,\eta) H_j^{M,k}(\eta)e^{V_j^{M,k}+c_{t,k}}dV_g+O(\varepsilon^{(M'+2)(m^*-2){-\delta}})\nonumber\\
&=:m_{itv}^kG(x,q_t^{M,k})+T_{M,k}^{G}+O(\varepsilon^{(M'+2)(m^*-2){-\delta}}).
\end{align}
Here, $T_{M,k}^{G}$ is a term of order $o_\tau(\varepsilon^{m^*-2})$. Notice that we used the frequency analysis in Subsection \ref{subsection:BubblingInside}.
Now we derive the expansion outside the bubbling disks. By a classical estimate \cite[Proposition 3.1]{LinZhang2011}, we get
\begin{align}
u_j^k(x)=-(m_{it}^k-2)\ln\frac{1}{\varepsilon_{1}^k}+O(1).\nonumber
\end{align}
This implies that
\begin{align}
&\int_{M\backslash(\cup_t B_{\tau_k^M}(q_t^{M,k}))} G(x,\eta)\sum_{j=1}^n a_{ij}\rho_j^k h_j^k e^{u_j^k}d V_g(\eta)=O(\varepsilon^{m^*-2{-\delta}}).\nonumber
\end{align}
Therefore, we get
\begin{align}
u_i^k(x)=\int_M u_i^k+2\pi\sum_{t=1}^N m_{itv}^k G(x;q_t^{M,k})+O(\varepsilon^{m^*-2-\delta}).\nonumber
\end{align}
Around each $p_t^{M,k}$, we rewrite it as
\begin{align}\label{Expansion:Outside01}
u_i^k(x)=\int_M u_i^k-2\pi m_{itv}^k\ln|x|+2\pi m_{itv}^k G^*(q_t^{M,k};q_t^{M,k})+O(\varepsilon^{m^*-2{-\delta}}).
\end{align}
In our next step, we provide an estimate on $\int_{M}u_i^k$ as follows. As we can find in Subsection \ref{subsection:BubblingInside}, we get the following expansion on the boundary of the $t$-th bubbling disk:
\begin{align}\label{Expansion:Boundary}
u_{it}^k(x)=-m_{itv}^k\ln|x|-(m_{itv}^k-2)\ln\frac{1}{\varepsilon_t^k}+I_t^k+T^{u,in}_{M',k}-\ln(\rho_i^k h_i^k(q_t^{M,k}))+O(\varepsilon^{(M'+2)(m^*-2){-\delta}}).
\end{align}
Here, $T^{u,in}_{M',k}$ is a higher order expansion of $u_{it}^k$.
Matching (\ref{Expansion:Outside01}) with (\ref{Expansion:Boundary}), we get
\begin{align}\label{Expansion:AverageRough}
\int_M u_i^k=-(m_{itv}^k-2)\ln\frac{1}{\varepsilon_t^k}-2\pi m_{itv}^k G^*(p_t^{M,k};p_t^{M,k})+I_t^k-\ln(\rho_i^kh_i^k(p_t^{M,k}))+O(\varepsilon^{m^*-2{-\delta}}).
\end{align}
Plugging this back into the integral outside of bubbling disks, we get new expansions of $u_i^k(x)$ and $\int_M u_i^k$. Repeating this calculation and we complete the proof.

\noindent{$\Box$}

An immediate corollary is that

\begin{corollary}\label{coro:RefinedRate}
For any large integer $M$, it holds that 
\begin{align}
\Bigg(\frac{\varepsilon^k_t}{\varepsilon^s_k}\Bigg)^{m_i^k-2}=\frac{e^{H_{i,t}(q^{M,k})}}{e^{H_{i,s}(q^{M,k})}} +T_{M,k,i,t,s}^{ratio}+O(\varepsilon^{M(m^*-2){-\delta}}).\nonumber
\end{align}
Here, $T_{M,k,i,t,s}^{ratio}$ depends $C^2$-continuously in $p_{t}^{M,k}$, $\varepsilon_t$ and the derivatives of $h_i^k$ up to $M+3$-th order for $t=1,\cdots,N$.
\end{corollary}
This is a refined version of (\ref{e:VarepsilonTVarespilonS}). Another important corollary is as follows.
{
\begin{corollary}\label{coro:TauMk}
For any integer $M$, it holds that
\begin{align}
\tau_k^M=\tau_0+T^\tau_{M,k,t}+O(\varepsilon^{M(m^*-2)-\delta}).\nonumber
\end{align}
Here, $T_{M,k,t}^{\tau}=O(\varepsilon^2)$ and depends $C^2$-continuously in $p_{t}^{M,k}$, $\varepsilon_t$ and the derivatives of $h_i^k$ up to $M+3$-th order for $t=1,\cdots,N$.
\end{corollary}
}

\subsection{A refined expansion of $\Lambda_{I,N}(\rho^k)$}

In this subsection, we obtain a refined expansion of $\Lambda_{I,N}(\rho^k)$.
\begin{theorem}\label{t:ExpansionLambda}
Under the assumptions of Theorem \ref{t:MAIN} and if $m^*<4$, we get
\begin{align}
\Lambda_{I,N}(\rho^k)&=\sum_{i=1}^n \sum_{t=1}^N \frac{2-m^*_i}{ N} D_{it}^{k} \frac{e^{H_{i,t}(q^{M,k})}}{e^{H_{i,1}(q^{M,k})}}(\varepsilon_1^k)^{m_i^k-2}+T^\Lambda_{M,k}+O(\varepsilon^{M(m^*-2){-\delta}}).\nonumber
\end{align}
Here,
\begin{align}
D^{k}_{it}=\lim_{\tau_k^M\to0+}\Bigg[ e^{I_i}\int_{\Omega_t\backslash B^M_{\tau_k^M}(q_t^{M,k})} \frac{h_i^k(x)}{h_i^k(q_t^{M,k})}e^{2\pi m_i^k\sum_{l=1}^N[G(x,q_l^{M,k})-G^*(q_l^{M,k},q_t^{M,k})]}dV_g-\frac{2\pi e^{I_i}}{m_i^*-2}\frac{1}{(\tau_k^M)^{m_i^k-2}}\Bigg]\nonumber
\end{align}
and $T_{M,k}^\Lambda=o_\tau(1)\varepsilon^{m^*-2}$ and is $C^2$ depends on $p_{i}^{M,k}$, $h_i$ and its derivatives up to $M+3$ order.
If $m^*=4$, we get
\begin{align}
\Lambda_{I,N}(\rho^k)=-\frac{2}{N}\sum_{i=1}^n\sum_{t=1}^N\Bigg[b^k_{it}(\varepsilon_t^k)+D_{it}^k\frac{e^{H_{it}(q^{M,k})}}{e^{H_{i1}(q^{M,k})}}(\varepsilon_1^k)^2\Bigg]+T^{\Lambda,*}_{M,k}+O(\varepsilon^{2M-\delta}).
\end{align}
Here, $b_{it}^k$ is defined as in (\ref{def:bik}).
\end{theorem}

\noindent{\bf Proof.}
{We only study the case $m^*<4$ since the cases $m^*=4$ is similar.}
Using the notations in Subsection
\ref{subsection:ClassicalBubbling}.
We begin with the expansion of $\rho^k$. By a direct computation, we find
\begin{align}
\rho_i^k&=\rho_i^k \int_M h_i^k e^{u_i^k}dV_g=\rho_i^k \sum_{t=1}^N \int_{B^M_{\tau_k^M}(q_t^{M,k})}+ \rho_i^k \int_{M\backslash \cup_{t=1}^N B^M_{\tau_k^M}(q_t^{M,k})} h_i^k e^{u_i^k}dV_g =:\sum_{t=1}^N \rho_{ti}^k+\rho_{bi}^k.\nonumber
\end{align}

We begin with $\rho_{bi}^k$.
Recall that we divide the domain into $N$ parts $\Omega_1,\cdots,\Omega_N$ such that $p_t^*\in\Omega_t$, $\overline{\cup_{t=1}^N\Omega_t}=M$ and $\Omega_{t_1}\cap\Omega_{t_2}=\emptyset$ for $t_1,t_2=1,\cdots,N$ with $t_1\neq t_2$. By Proposition \ref{prop:OutAverageExpansion}, we get
\begin{align}\label{ineq:rhoib}
\rho_{bi}^k &=\rho_i^k \int_{M\backslash \cup_{t=1}^N B^M_{\tau_k^M}(q_t^{M,k})} h_i^k e^{u_i^k}dV_g =\sum_{t=1}^N\rho_i^k\int_{\Omega_t\backslash B^M_{\tau_k^M}(q_t^{M,k})} h_i^k e^{u_i^k}dV_g\nonumber\\
&=e^{I_i}\sum_{t=1}^N(\varepsilon_t^k)^{m_i^k-2}\int_{\Omega_t\backslash B^M_{\tau_k^M}(q_t^{M,k})} \frac{h_i^k(x)}{h_i^k(q_t^{M,k})}e^{2\pi m_i^k\sum_{l=1}^N[G(x,q_l^{M,k})-G^*(q_l^{M,k},q_t^{M,k})]}dV_g+T_{M,k,i}^{\rho,out}+O(\varepsilon^{M(m^*-2){-\delta}}).
\end{align}
Here, $T_{M,k,i}^{\rho,out}=O(\varepsilon^{2m^*-4})$ and depends $C^2$-continuously on $p_{t}^{M,k}$ and the derivatives of $h_i^k$ up to $M+3$-th order.
As for the integrals $\rho_{ti}^k$'s, by Theorem \ref{t:IntegralExpansion}, we get
\begin{align}
\rho_{it}^k
&=\rho_i^k\int_{B_{\tau_k^M}^M(q_t^{M,k})}h_i^k e^{u_i^k}dV_g=\int_{B_{{\tau_k^M/\varepsilon_t^k}}(0)}H_{t,i}^{M,k}( y)e^{\phi_{i,t}^k(q_t^{M,k}+\varepsilon_t^k y)}e^{v_i^k(y)}dy\nonumber\\
&=H_{t,i}^{M,k}(0)\int_{\mathbb{R}^2}e^{V_{t,i}^k}-{\frac{2\pi H_{t,i}^{1,k}(0) e^{I_i^k}}{m_i^k -2}\Big(\frac{\varepsilon_t^k}{\tau_k^M}\Big)^{m_i^k-2}}+T_{M,k,i,t}^{in}+O(\varepsilon^{M(m^*-2){-\delta}}).\nonumber
\end{align}
Here, $H_{t,i}^k(\cdot)=\frac{h^k_i(q_t^{M,k}+\cdot)}{h_i^k(q_t^{M,K})}e^{\psi+f_i^k}$ and
$T_{M,k,i,t}^{in}$ is the term defined similarly as in Theorem \ref{t:IntegralExpansion}.
Then, we get
\begin{align}
\rho_i^k=&\sum_{t=1}^N H_{t,i}^{M,k}(0)\int_{\mathbb{R}^2}e^{V_{t,i}^k}+e^{I_i}\sum_{t=1}^N(\varepsilon_t^k)^{m_i^k-2}\int_{\Omega_t\backslash B^M_{\tau_k^M}(q_t^k)} \frac{h_i^k(x)}{h_i^k(q_t^{M,k})}e^{2\pi m_i^k\sum_{l=1}^N[G(x,q_l^{M,k})-G^*(q_l^{M,k},q_t^{M,k})]}dV_g\nonumber\\
&\quad-\sum_{t=1}^N\frac{2\pi H_i^{M,k}(0) e^{I_i^k}}{m_i^k -2}\Big(\frac{\varepsilon_t^k}{\tau_k^M}\Big)^{m_i^k-2}+\sum_{t=1}^N T_{M,k,i,t}^{in}+T_{M,k,i}^{\rho,out}+O(\varepsilon^{M(m^*-2){-\delta}}).\nonumber
\end{align}
By (\ref{def:psi}), (\ref{def:def:fik}) and Lemma \ref{l:pik}, we get
\begin{align}\label{ineq:HinteV}
\int_{\mathbb{R}^2}e^{V_{t,i}^k}=\int_{\mathbb{R}^2}e^{V_{1,i}^k}+O(\varepsilon_k^{M(m^*-2)}).
\end{align}
Therefore,
\begin{align}
&4\sum_{i=1}^n\frac{\rho_i^k}{2\pi N}=\nonumber\\
&\quad4\sum_{i=1}^n\frac{\int_{\mathbb{R}^2}e^{V_{1,i}^k}}{2\pi}+\frac{4}{2\pi N}\sum_{i=1}^n e^{I_i}\sum_{t=1}^N(\varepsilon_t^k)^{m_i^k-2}\int_{\Omega_t\backslash B^M_{\tau_k^M}(q_t^{M,k})} \frac{h_i^k(x)}{h_i^k(q_t^{M,k})}e^{2\pi m_i^k\sum_{l=1}^N[G(x,q_l^{M,k})-G^*(q_l^{M,k},q_t^{M,k})]}dV_g\nonumber\\
&\quad-\frac{4}{N}\sum_{i=1}^n\sum_{t=1}^N\frac{e^{I_i^k}}{m_i^k -2}\Big(\frac{\varepsilon_t^k}{\tau_k^M}\Big)^{m_i^k-2}+\frac{4}{2\pi N}\sum_{t=1}^N\sum_{i=1}^n \big(T_{M,k,i,t}^{in}+T_{M,k,i}^{\rho,out}\big)+O(\varepsilon^{M(m^*-2){-\delta}}).\nonumber
\end{align}
On the other hand, by a direct computation, we get
\begin{align}
\sum_{i,j=1}^n a_{ij}&\frac{\rho_i^k}{2\pi N} \frac{\rho_j^k}{2\pi N}\nonumber\\
&=4\sum_{i=1}^n\frac{\int_{\mathbb{R}^2}e^{V_{1,i}^k}}{2\pi}+\sum_{i=1}^n \sum_{t=1}^N \frac{2m_{t,i}^i}{2\pi N} e^{I_i}(\varepsilon_t^k)^{m_i^k-2}\int_{\Omega_t\backslash B^M_{\tau_k^M}(q_t^{M,k})} \frac{h_i^k(x)}{h_i^k(q_t^{M,k})}e^{2\pi m_i^k\sum_{l=1}^N[G(x,q_l^{M,k})-G^*(q_l^{M,k},q_t^{M,k})]}dV_g\nonumber\\
&\quad-\frac{1}{N}\sum_{i=1}^n\sum_{t=1}^N\frac{2m_{t,i}^k}{m_{t,i}^k-2}\Big(\frac{\varepsilon_t^k}{\tau_k^M}\Big)^{m_i^k-2}+\frac{2}{2\pi N} \sum_{t=1}^N\sum_{i=1}^n m_{t,i}^k \big(T_{M,k,i,t}^{in}+T_{M,k,i}^{\rho,out}\big)+O(\varepsilon^{M(m^*-2){-\delta}}).\nonumber
\end{align}
Therefore, passing $\tau\to0+$, we get
\begin{align}
\Lambda_{I,N}(\rho^k)=4\sum_{i=1}^n\frac{\rho_i^k}{2\pi N}-\sum_{i,j=1}^n a_{ij}\frac{\rho_i^k}{2\pi N} \frac{\rho_j^k}{2\pi N}=:\sum_{i=1}^n \sum_{t=1}^N D_{it}^k\frac{e^{{H_{i,t}(q^{M,k})}}}{e^{H_{i,1}(q^{M,k})}} (\varepsilon_1^k)^{m_i^k-2}+T^\Lambda_{M,k}+O(\varepsilon^{M(m^*-2){-\delta}})\nonumber
\end{align}
with $T_{M,,k}^\Lambda=o_\tau(1)\varepsilon^{m^*-2}$ and is $C^2$ depends on $p_{i}^{M,k}$, $h_i$ and its derivatives up to $M+3$ order.
Notice that in the last step we use Corollary \ref{coro:RefinedRate}.
This proves the result.

\noindent{$\Box$}

{
\begin{remark}
This is a refined version of \cite[Theorem 1.2]{LinZhang2013} and \cite[Theorem 1.3]{HuangZhang2022}. Notice that the expansion of $\Lambda_{I,N}$ for $m^*=4$ fits the results in \cite[Theorem 1.2]{LinZhang2013} and \cite[Theorem 1.3]{HuangZhang2022} due to Remark \ref{r:LeadingTermL}.
\end{remark}
}

\subsection{A refined Pohozaev identity}
Let us recall the classical Pohozaev identity for Liouville systems (see \cite[pp. 2620]{LinZhang2013}):
\begin{align}\label{e:PohozaecClassical}
\sum_i\int_{\Omega}&\partial_{x_s}H_i^{M,k}(y)e^{v_i^{M,k}(y)}\nonumber\\
&=\int_{\partial\Omega}\nu_s\sum_j e^{v_i^{M,k}}H_i^{M,k}+\sum_{ij}a^{ij}\Bigg(\partial_{x_s} v_i^{M,k}\partial_{\nu}v_i^{M,k}-\frac{1}{2}\nabla v_i^{M,k}\cdot\nabla v_j^{M,k}\nu_s\Bigg).
\end{align}
Here, $s=1,2$.
Now we derive a 
\begin{proposition}\label{prop:Pohozaev}
It holds that
\begin{align}\label{e:Pohozaev}
\int_{B_{\tau_M^k/\varepsilon_1^k}(p_t^{M,k})}\partial_1 H_i^{M,k}(y)e^{V_{t,i}^{M,k}}+T^{Po}_{M,k}+O(\varepsilon^{(M+1)(m^*-2)-\delta})=0
\end{align}
for certain term $T^{Po}_{M,k}$ { of order $\varepsilon^{m^*-2}$ for the case $m^*<4$ and $\varepsilon^2\ln\frac{1}{\varepsilon}$ for $m^*=4$} and depends $C^1$-continuously in $p_t^{M,k}$, $H^{M,k}_j$ and its derivatives up to $M+2$ order.
\end{proposition}
\noindent{\bf Proof.}
For the left hand side of (\ref{e:PohozaecClassical}), we get
\begin{align}
\partial_1 H_i^{M,k}(y)e^{v_i^{M,k}(y)}
&=\partial_1 H_i^{M,k}(y) e^{V_{t,i}^{M,k}}\Bigg(1+A^M_k+O((\varepsilon_1^k)^{M+1}(1+|y|)^{(M+1)(3-m^*){+\delta}})\Bigg)\nonumber\\
&=\partial_1 H_i^{M,k}(y)e^{V_{t,i}^{M,k}}+\overline{A}_k^M+O((\varepsilon_1^k)^{M+1}(1+|y|)^{(M+1)(3-m^*){+\delta}}).\nonumber
\end{align}
Here, $A_k^M$ and $\overline{A}_k^M$ are two terms can be defined inductively.
Integrating the above on $B_{\tau\varepsilon^{-1}}(p_t^{M,k})$, we get
\begin{align}
\int_{B_{\tau_M^k/\varepsilon_1^k}(p_t^{M,k})}\partial_1 H_i^{M,k}(y)e^{v_i^{M,k}(y)}=&\int_{B_{\tau_M^k/\varepsilon_1^k}(p_t^{M,k})}\partial_1 H_i^{M,k}(y)e^{V_{t,i}^{M,k}}\nonumber\\
&\quad+\int_{B_{\tau_M^k/\varepsilon_1^k}(p_t^{M,k})}\overline{A}_k^M+O(\varepsilon^{(M+1)(m^*-2)-\delta}).\nonumber
\end{align}
Here, the integral $\int_{B_{\tau_M^k/\varepsilon_1^k}(p_t^{M,k})}\overline{A}_k^M$ is of $\varepsilon^{m^*-2}$ order and depends $C^1$-continuously on $p_t^{M,k}$, $H_j^{M,k}$ and its derivatives up to $M+3$ order.

On the other hand, by the decomposition (\ref{Expansion:viMk}) of $v_i^{M,k}$, we get
\begin{align}
v_{i}^{M,k}=V_{i,t}^{M,k}+c_{M,i}+\phi^{**,k}_{1,i}+\phi^{**,k}_{2,i}+W^*_{i,k}\nonumber
\end{align}
and
\begin{align}
\nabla v_{i}^{M,k}=\nabla V_{i,t}^{M,k}+\nabla c_{M,i}+\nabla \phi^{**,k}_{1,i}+\nabla \phi^{**,k}_{2,i}+\nabla W^*_{i,k}.\nonumber
\end{align}
By these, we get
\begin{align}
&\int_{\partial B_{\tau_M^k/\varepsilon_1^k}(p_t^{M,k})}\sum_{ij}a^{ij}\Bigg(\partial_1 v_i^{M,k}\partial_{\upsilon}v_j^{M,k}-\frac{1}{2}\nabla v_i^{M,k}\cdot\nabla v_j^{M,k}\upsilon_s\Bigg)\nonumber\\
&=\int_{\partial B_{\tau_M^k/\varepsilon_1^k}(p_t^{M,k})}\sum_{ij}a^{ij}\Bigg[\partial_1\Big(V_{1,i}^{M,t}+c_{M,i}\Big)\partial_\upsilon\Big(V_{1,j}^{M,t}+c_{M,j}\Big)-\frac{1}{2}\nabla\Big(V_{1,i}^{M,t}+c_{M,i}\Big)\cdot\nabla\Big(V_{1,j}^{M,t}+c_{M,j}\Big)\Bigg]\nonumber\\
&\quad+O(\varepsilon^{(M+1)(m^*-2)-\delta})\nonumber
\end{align}
and
\begin{align}
\int_{\partial B_{\tau_M^k/\varepsilon_1^k}(p_t^{M,k})}\upsilon_s\sum_i e^{v_i^{M,k}}H_i^{M,k}=\int_{\partial B_{\tau_M^k/\varepsilon_1^k}(p_t^{M,k})}\upsilon_s\sum_{i}H_i^{M,k}e^{V_{i,t}^{M,t}+c_{M,i}}+O(\varepsilon^{(M+1)(m^*-2)-\delta}).\nonumber
\end{align}
It should be noticed that, by the symmetry and the Fourier decomposition, we get
\begin{align}
\int_{\partial B_{\tau_M^k/\varepsilon_1^k}(p_t^{M,k})}\sum_{ij}a^{ij}\Bigg(\partial_1 v_i^{M,k}\partial_{\upsilon}v_j^{M,k}-\frac{1}{2}\nabla v_i^{M,k}\cdot\nabla v_j^{M,k}\upsilon_s\Bigg)=\left\{
\begin{aligned}
&O(\varepsilon_1^{m^*}) &\mbox{ if } m^*<4;\nonumber\\
&O(\varepsilon_1^2\ln\frac{1}{\varepsilon_1})&\mbox{ if } m^*= 4\nonumber
\end{aligned}
\right.
\end{align}
and
\begin{align}
\int_{\partial B_{\tau_M^k/\varepsilon_1^k}(p_t^{M,k})}\upsilon_s\sum_i e^{v_i^{M,k}}H_i^{M,k}=\left\{
\begin{aligned}
&O(\varepsilon_1^{m^*}) &\mbox{ if } m^*<4;\nonumber\\
&O(\varepsilon_1^2\ln\frac{1}{\varepsilon_1})&\mbox{ if } m^*= 4.\nonumber
\end{aligned}
\right.
\end{align}
Combining the above, we get
\begin{align}
\int_{B_{\tau_M^k/\varepsilon_1^k}(p_t^{M,k})}\partial_1 H_i^{M,k}(y)e^{V_{t,i}^{M,k}}+T^{Po}_{M,k}+O(\varepsilon^{(M+1)(m^*-2)-\delta})=0\nonumber
\end{align}
for {certain term $T_{M,k}^{Po}$ of order $\varepsilon^{m^*-2}$ for the case $m^*<4$ and $\varepsilon^2\ln\frac{1}{\varepsilon}$ for $m^*=4$} and depends $C^1$-continuously in $p_t^{M,k}$, $H^{M,k}_j$ and its derivatives up to $M+2$ order.

\noindent{$\Box$}

\section{Preliminary estimates on the difference between solutions}

Since we prove Theorem \ref{t:MAIN} by contradiction, we study two sequences of blowup solutions to Problem (\ref{e:00k}) from now on. To be precise, we assume that there are two sequences of solutions
\begin{align}
u^{(1),k}=(u^{(1),k}_1,\cdots,u^{(1),k}_n)\mbox{ and }u^{(2),k}=(u^{(2),k}_1,\cdots,u^{(2),k}_n)\nonumber
\end{align}
solving Problem (\ref{e:00k}), i.e.
\begin{equation}
    \left\{
   \begin{array}{lr}
     -\Delta u^{(l),k}_i=\sum_{j=1}^n a_{ij}\rho^k_j\Big(\frac{h^k_j e^{u^{(l),k}_j}}{\int_M h^k_j e^{u^{(l),k}_j}}-1\Big),\nonumber\\
     u^{(l),k}_i\in \mathring{H}^1(M)\mbox{ for }i=1,\cdots,n,\,\,l=1,2.\nonumber
   \end{array}
   \right.
\end{equation}
Before going further, we fix the following notations:
\begin{itemize}
  \item $q^*:=(q_1^*,\cdots,q_N^*)\in M^N$ denotes the blowup points of $\{u^{(1),k}\}_{k=1}^\infty$ and $\{u^{(2),k}\}_{k=1}^\infty$;
  \item For any $l=1,2$, $k=1,2,\cdots$ and $t=1,\cdots,N$, we denote $M_t^{(l),k}=\max_{x\in B_{\tau_0}(p_t^*)}u_{t}^{(l),k}(x)$ and $p_t^{(l),k}$ is the point where the maximum is achieved. Denote $M^{(l),k}=\max_{t=1,\cdots,N} M_t^{(l),k}$. Denote $q_t^{(l),k}:=\varepsilon_t^{(l),k}p_t^{(l),k}$;
  \item the corresponding radii $\tau_k^M$ are denoted as $\tau_k^{(1)}$ and $\tau_k^{(2)}$, respectively;
  \item $\varepsilon^{(l),k}:=e^{-\frac{M^{(l),k}}{2}}$ for any $l=1,2$ and $k=1,2,\cdots$;
  \item Without loss of generality, we assume that $\int_M h_i e^{u_i^{(l),k}}=1$ for $l=1,2$ and $i=1,\cdots,n$.
\end{itemize}


\subsection{Comparison between $\varepsilon_t^{(1),k}$ and $\varepsilon_t^{(2),k}$}

\begin{proposition}\label{p:varepsilon12}
Under the assumptions in Theorem \ref{t:MAIN},
we get
\begin{align}
\Bigg|1-\frac{\varepsilon_t^{(2),k}}{\varepsilon_t^{(1),k}}\Bigg|\leq C\varepsilon\Bigg| p_t^{(1),k}-\frac{\varepsilon_t^{(2),k}}{\varepsilon_t^{(1),k}}p_t^{(2),k}\Bigg|+O(\varepsilon^{M(m^*-2){-\delta}}).\nonumber
\end{align}
\end{proposition}
\noindent{\bf Proof.}
{We only study the case $m^*<4$ since the cases $m^*=4$ is similar.}
Due to the assumption, we get
\begin{align}\label{e:D=D01}
\Lambda_{I,N}(\rho^k)&=\sum_{i=1}^n \sum_{t=1}^N D_{it}^{k,(1)} \frac{e^{H_{it}(p^{(1),k})}}{e^{H_{i1}(p^{(1),k})}}(\varepsilon_1^{(1),k})^{m_i^k-2}+T_{M,k}^{(1),\Lambda}+O((\varepsilon_t^{(2),k})^{M(m^*-2){-\delta}})\nonumber\\
&=\sum_{i=1}^n \sum_{t=1}^N D_{it}^{k,(2)} \frac{e^{H_{it}(p^{(2),k})}}{e^{H_{i1}(p^{(2),k})}}(\varepsilon_1^{(2),k})^{m_i^k-2}+T_{M,k}^{(2),\Lambda}+O((\varepsilon_t^{(2),k})^{M(m^*-2){-\delta}}).
\end{align}
Here, the terms $D_{it}^{k,(1)}$, $D_{it}^{k,(2)}$, $T_{M,k}^{(1),\Lambda}$, $T_{M,k}^{(2),\Lambda}$ and are defined as in Theorem \ref{t:ExpansionLambda}.

An elementary observation implies that $\varepsilon_t^{(1),k}/\varepsilon_t^{(2),k}\to 1$ as $k\to\infty$. From now on, for the sake of simplicity, we refer to these two quantities collectively as $\varepsilon$ in the error parts. 

\begin{center}
\begin{figure}[htbp]
\centering
\includegraphics[width=0.8\textwidth]{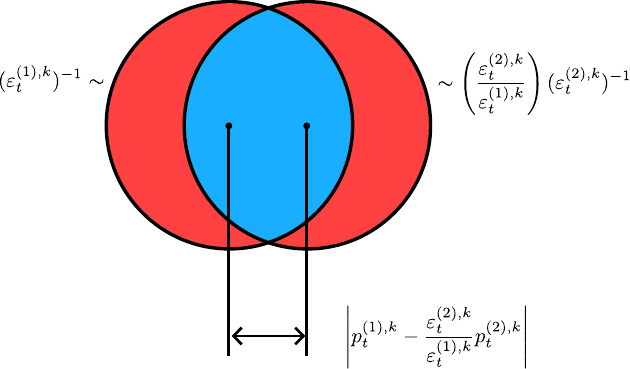}
\caption{The effect of  the symmetric difference of the domain}
\label{fig:DifferentBalls}
\end{figure}
\end{center}

Recall that due to the notions in 
Theorem \ref{t:ExpansionLambda}, we get
\begin{align}
T_{M,k}^{(l),\Lambda}=\frac{1}{2\pi N} \sum_{t=1}^N\sum_{i=1}^n (2-m_{t,i}^k )\big(T_{M,k,i,t}^{(l),in}+T_{M,k,i}^{(l),\rho,out}\big)\nonumber
\end{align}
for $l=1,2$.
Notice that the quantities $T_{M,k,i,t}^{in}$ and $T_{M,k,i}^{\rho,out}$ depend on $p^{(l),k}_t$ and $\varepsilon_t^{(l),k}$ for $l=1,2$. Then, a direct computation implies that
\begin{align}
|T_{M,k,i,t}^{(1),in}-T_{M,k,i,t}^{(2),in}|\leq A+B.\nonumber
\end{align}
Here, $A$ denotes the difference due to the symmetric difference of the domain and $B$ denotes the difference due to the difference of $p_t^{(1),k}$ and $p_t^{(2),k}$.  We use Figure \ref{fig:DifferentBalls} to show the effect of  the symmetric difference of the domain. 
As it is shown in Figure \ref{fig:DifferentBalls}, $A$ denotes the difference in the red domains and $B$ denotes the difference in the blue domains.
Then, we get
\begin{align}
|A|\leq C\varepsilon^2\int_{\frac{\min(\tau_k^{(1)},\tau_k^{(2)})}{\varepsilon_t^{(1),k}}}^{\Big| p_t^{(1),k}-\frac{\varepsilon_t^{(2),k}}{\varepsilon_t^{(1),k}}p_t^{(2),k}\Big|+\frac{\tau_k^{(1)}+\tau_k^{(2)}}{\varepsilon_t^{(1),k}}}r^3(1+r)^{-m^*}dr=O\Bigg(\varepsilon^{m^*-1}\Bigg| p_t^{(1),k}-\frac{\varepsilon_t^{(2),k}}{\varepsilon_t^{(1),k}}p_t^{(2),k}\Bigg|\Bigg)\nonumber
\end{align}
and
\begin{align}
|B|\leq C\varepsilon^3\int_0^{\tau_k^{(1)}+\tau_k^{(2)}}r^3(1+r)^{-m^*}dr\times\Bigg| p_t^{(1),k}-\frac{\varepsilon_t^{(2),k}}{\varepsilon_t^{(1),k}}p_t^{(2),k}\Bigg|=O\Bigg(\varepsilon^{m^*-1}\Bigg| p_t^{(1),k}-\frac{\varepsilon_t^{(2),k}}{\varepsilon_t^{(1),k}}p_t^{(2),k}\Bigg|\Bigg).\nonumber
\end{align}
The estimate on $B$ is due to the $C^1$-dependence on the parameters $p_t^{(1),k}$ and $p_t^{(2),k}$.
These implies that
\begin{align}
|T_{M,k,i,t}^{(1),in}-T_{M,k,i,t}^{(2),in}|=O\Bigg(\varepsilon^{m^*-1}\Bigg| p_t^{(1),k}-\frac{\varepsilon_t^{(2),k}}{\varepsilon_t^{(1),k}}p_t^{(2),k}\Bigg|\Bigg).\nonumber
\end{align}
For terms $T_{M,k,i}^{(1),\rho,out}$ and $T_{M,k,i}^{(2),\rho,out}$, we get a similar result.
\begin{align}
|T_{M,k,i}^{(1),\rho,out}-T_{M,k,i}^{(2),\rho,out}|=O\Bigg(\varepsilon^{m^*-1}\Bigg| p_t^{(1),k}-\frac{\varepsilon_t^{(2),k}}{\varepsilon_t^{(1),k}}p_t^{(2),k}\Bigg|\Bigg).\nonumber
\end{align}

Plugging these into (\ref{e:D=D01}), we get
\begin{align}
\Bigg|1-\frac{\varepsilon_t^{(2),k}}{\varepsilon_t^{(1),k}}\Bigg|\leq C\varepsilon\Bigg| p_t^{(1),k}-\frac{\varepsilon_t^{(2),k}}{\varepsilon_t^{(1),k}}p_t^{(2),k}\Bigg|+O(\varepsilon_k^{M(m^*-2)}).\nonumber
\end{align}
This completes the proof.

\noindent{$\Box$}

\begin{remark}
From now on, for the sake of simplicity, we refer to these two quantities collectively as $\varepsilon$ in the error parts. 
\end{remark}

\subsection{Comparison between $p_t^{(1),k}$ and $p_t^{(2),k}$}
The difference between $p_t^{(1),k}$ and $p_t^{(2),k}$ is due to the comparison of Pohozaev identity (\ref{e:Pohozaev}) and the non-degeneracy assumption.
\begin{proposition}\label{prop:p1-p2}
Under the assumptions of Theorem \ref{t:MAIN}, we get
\begin{align}
|\varepsilon_t^{(1),k}p_t^{(1),k}-\varepsilon_t^{(2),k}p_t^{(2),k}|\leq o(1)\Bigg|1-\frac{\varepsilon_t^{(1),k}}{\varepsilon_t^{(2),k}}\Bigg|+O(\varepsilon^{(M+1)(m^*-2)-\delta}).\nonumber
\end{align}
\end{proposition}
\noindent{\bf Proof.}
For the solutions $u^{(1),k}=(u^{(1),k}_1,\cdots,u^{(1),k}_n)$ and $u^{(2),k}=(u^{(2),k}_1,\cdots,u^{(2),k}_n)$, by Proposition \ref{prop:Pohozaev}, we get the following two Pohozaev identity:
\begin{align}
\int_{B_{\tau_k^M/\varepsilon_t^{(1),k}}(p_t^{(1),k})}\partial_1 H_i^{M,k}(y)e^{V_{t,i}^{(1),k}}=-T^{(1),Po}_{M,k}+O(\varepsilon^{(M+1)(m^*-2)-\delta})\nonumber
\end{align}
and
\begin{align}
\int_{B_{\tau_k^M/\varepsilon_t^{(2),k})}(p_t^{(2),k})}\partial_1 H_i^{M,k}(y)e^{V_{t,i}^{(2),k}}=-T^{(2),Po}_{M,k}+O(\varepsilon^{(M+1)(m^*-2)-\delta}).\nonumber
\end{align}
We will prove the proposition by differencing these two identities. Notice that
\begin{align}
\int_{B_{\tau_k^M/\varepsilon_t^{(1),k})}(p_t^{(1),k})}&\partial_1 H_i^{M,k}(y)e^{V_{t,i}^{(1),k}}\nonumber\\
&=\partial_1 H_i^{M,k}(p_t^{(1),k})\int_{\mathbb{R}^2}e^{V_{t,i}^{(1),k}}+\int_{B_{\tau_k^M/\varepsilon_t^{(1),k}}(p_t^{(1),k})}(\partial_1 H_i^{M,k}(y)-\partial_1 H_i^{M,k}(0))e^{V_{t,i}^{(1),k}}\nonumber\\
&\quad+\partial_1 H_i^{M,k}(p_t^{(1),k})\int_{\mathbb{R}^2\backslash B_{\tau_k^M/\varepsilon_t^{(1),k}}(p_t^{(1),k})}e^{V_{t,i}^{(1),k}}.\nonumber
\end{align}
A similar identity hold for solution $u^{(2),k}$. Therefore, we get
\begin{align}\label{ineq:p1-p2-01}
\int_{B_{\tau_k^M/\varepsilon_t^{(1),k}}(p_t^{(1),k})}&\partial_1 H_i^{M,k}(y)e^{V_{t,i}^{(1),k}}-\int_{B_{\tau_k^M/\varepsilon_t^{(2),k}}(p_t^{(2),k})}\partial_1 H_i^{M,k}(y)e^{V_{t,i}^{(2),k}}\nonumber\\
&\geq C|\varepsilon_t^{(1),k}p_t^{(1),k}-\varepsilon_t^{(2),k}p_t^{(2),k}|-O\Bigg(\varepsilon^{m^*-2}\Bigg|1-\frac{\varepsilon_t^{(1),k}}{\varepsilon_t^{(2),k}}\Bigg|\Bigg)-O(\varepsilon^{(M+1)(m^*-2)-\delta}).
\end{align}
On the other hand, using the notions in Proposition \ref{prop:Pohozaev}, we get
\begin{align}\label{ineq:p1-p2-02}
-T^{(1),Po}_{M,k}+T^{(2),Po}_{M,k}=O(\varepsilon^{m^*-2})\Bigg(\Bigg|1-\frac{\varepsilon_t^{(1),k}}{\varepsilon_t^{(2),k}}\Bigg|+|\varepsilon_t^{(1),k}p_t^{(1),k}-\varepsilon_t^{(2),k}p_t^{(2),k}|\Bigg).
\end{align}
Combining (\ref{ineq:p1-p2-01}) and (\ref{ineq:p1-p2-02}), we finish the proof.

\noindent{$\Box$}

An important corollary of Proposition \ref{p:varepsilon12} and Proposition \ref{prop:p1-p2} is that
\begin{corollary}
It holds that $\Big|1-\frac{\varepsilon_t^{(1),k}}{\varepsilon_t^{(2),k}}\Big|=O(\varepsilon^{M(m^*-2)})$ and $\Big|p_t^{(1),k}-\frac{\varepsilon_t^{(2),k}}{\varepsilon_t^{(1),k}}p_t^{(2),k}\Big|=O(\varepsilon^{M(m^*-2)})$ for $t=1,\cdots,N$ and a sufficiently large positive constant $M$. 
\end{corollary}

\subsection{Comparison between $u^{(1),k}$ and $u^{(2),k}$}

Now we estimate $\|u_i^{(1),k}-u_i^{(2),k}\|_{L^\infty(M)}$.

\begin{proposition}\label{prop:u1-u2}
Under the assumptions of Theorem \ref{t:MAIN}, we get
\begin{align}
\|u_i^{(1),k}-u_i^{(2),k}\|_{L^\infty(M)}=O(\varepsilon^{M(m^*-2)-2-\delta})
\end{align}
for sufficiently large $M$.
\end{proposition}
\noindent{\bf Proof.}
Our proof is divided into two parts.

~

\noindent{\bf Step 1. Estimate in $\cup_{t=1}^N B_\tau(p_t^*)$.}

~

By a direct computation, we get
\begin{align}\label{e:Theta-ThetaInner}
&u_i^{(1),k}(x)-u_i^{(2),k}(x)\nonumber\\
&=V_{it}^{(1),k}\Big(\frac{x}{\varepsilon_t^{(1),k}}-p_t^{(1),k}\Big)-V_{it}^{(2),k}\Big(\frac{x}{\varepsilon_t^{(2),k}}-p_t^{(2),k}\Big)+2\ln\frac{\varepsilon_t^{(2),k}}{\varepsilon_t^{(1),k}}+\phi_i^{(1),k}(x- \varepsilon_t^{(1),k} p_t^{(1),k})-\phi_i^{(2),k}(x-\varepsilon_t^{(2),k} p_t^{(2),k}))\nonumber\\
&\quad+c_{M,i}^{(1)}\Big(\frac{x}{\varepsilon_t^{(1),k}}-p_t^{(1),k}\Big)-c_{M,i}^{(2)}\Big(\frac{x}{\varepsilon_t^{(2),k}}-p_t^{(2),k}\Big)+\phi_{1,i}^{(1),**,k}\Big(\frac{x}{\varepsilon_t^{(1),k}}-p_t^{(1),k}\Big)-\phi_{1,i}^{(2),**,k}\Big(\frac{x}{\varepsilon_t^{(2),k}}-p_t^{(2),k}\Big)\nonumber\\
&\quad+\phi_{2,i}^{(1),**,k}\Big(\frac{x}{\varepsilon_t^{(1),k}}-p_t^{(1),k}\Big)-\phi_{2,i}^{(2),**,k}\Big(\frac{x}{\varepsilon_t^{(2),k}}-p_t^{(2),k}\Big)+W_{i,k}^{(1),*}\Big(\frac{x}{\varepsilon_t^{(1),k}}-p_t^{(1),k}\Big)-W_{i,k}^{(2),*}\Big(\frac{x}{\varepsilon_t^{(2),k}}-p_t^{(2),k}\Big)\nonumber\\
&=O(\varepsilon^{M(m^*-2)-\delta})
\end{align}
for any $x\in B_\tau(p_t^*)$ for any $t=1,\cdots,N$.
Here, we use (\ref{Expansion:viMk}).

~

\noindent{\bf Step 2. Estimate in $M\backslash \cup_{t=1}^N B_\tau(p_t^*)$.}

~

For suitable small $\tau_*>0$, we get $M\backslash \cup_{t=1}^N B_{\tau}(p_t^*)\subset M\backslash \cup_{t=1}^N B_{2\tau_*}(p_t^{(1),k})$. Therefore, it is only need to estimate $\|u_{i}^{(1),k}-u_i^{(2),k}\|_{L^\infty(M\backslash \cup_{t=1}^N B_{2\tau_0}(p_t^{(1),k}))}$. By Green's representation, we get for $x\in M\backslash \cup_{t=1}^N B_{\tau_0}(p_t^{(1),k})$
\begin{align}
(u_{i}^{(1),k}&(x)-\overline{u_{i}^{(1),k}})-(u_i^{(2),k}(x)-\overline{u_{i}^{(2),k}})\nonumber\\
&=\int_M G(x,\eta)\sum_{j=1}^n a_{ij}\rho_j^k h_j\big(e^{u_{j}^{(1),k}}-e^{u_{j}^{(1),k}}\big)d V_g(\eta)\nonumber\\
&=\sum_{t=1}^N\int_{B_{\tau_*}(p_t^{(1),k})}(G(x,\eta)-G(x,p_t^{(1),k}))\sum_{j=1}^n a_{ij}\rho_j^k h_j\big(e^{u_{j}^{(1),k}}-e^{u_{j}^{(2),k}}\big)d V_g(\eta)\nonumber\\
&\quad+\sum_{t=1}^N\int_{B_{\tau_*}(p_t^{(1),k})}G(x,p_t^{(1),k})\sum_{j=1}^n a_{ij}\rho_j^k h_j\big(e^{u_{j}^{(1),k}}-e^{u_{j}^{(2),k}}\big)d V_g(\eta)\nonumber\\
&\quad+\int_{M\backslash \cup_{t=1}^N B_{\tau_*}(p_t^{(1),k})}G(x,\eta)\sum_{j=1}^n a_{ij}\rho_j^k h_j\big(e^{u_{j}^{(1),k}}-e^{u_{j}^{(2),k}}\big)d V_g(\eta).\nonumber
\end{align}
Here,
for $x\in M\backslash \cup_{t=1}^N B_{2\tau_*}(p_t^{(1),k})$ and $\eta\in B_{\tau_*}(p_t^{(1),k})$, we get
\begin{align}
G(x,\eta)-G(x,p_t^{(1),k})=\nabla_\eta G(x,p_t^{(1),k})\cdot (\eta -p_t^{(1),k})+O(|\eta -p_t^{(1),k}|^2).\nonumber
\end{align}
Then, we get
\begin{align}\label{ineq:Theta-Theta-OUT-1}
\sum_{t=1}^N&\int_{B_{\tau_*}(p_t^{(1),k})}(G(x,\eta)-G(x,p_t^{(1),k}))\sum_{j=1}^n a_{ij}\rho_j^k h_j\big(e^{u_{j}^{(1),k}}-e^{u_{j}^{(2),k}}\big)d V_g(\eta)\nonumber\\
&=\sum_{t=1}^N\int_{B_{\tau_*}(p_t^{(1),k})}[\nabla_\eta G(x,p_t^{(1),k})\cdot (\eta -p_t^{(1),k})+O(|\eta -p_t^{(1),k}|^2)] \sum_{j=1}^n a_{ij}\rho_j^k h_je^{u_{j}^{(1),k}}\big(1-e^{u_{j}^{(2),k}-u_{j}^{(1),k}}\big)d V_g(\eta)\nonumber\\
&=O(\varepsilon^{M(m^*-2)-2-\delta}).
\end{align}
Similarly, we get
\begin{align}\label{ineq:Theta-Theta-OUT-2}
\sum_{t=1}^N
\int_{B_{\tau_*}(p_t^{(1),k})}G(x,p_t^{(1),k})\sum_{j=1}^n a_{ij}\rho_j^k h_j\big(e^{u_{j}^{(1),k}}-e^{u_{j}^{(2),k}}\big)d V_g(\eta)=O(\varepsilon^{M(m^*-2)-\delta}).
\end{align}
On the other hand, by Proposition \ref{prop:OutAverageExpansion}, we get
\begin{align}\label{ineq:Theta-Theta-OUT-3}
\int_{M\backslash \cup_{t=1}^N B_{\tau_*}(p_t^{(1),k})}G(x,\eta)\sum_{j=1}^n a_{ij}\rho_j^k h_j\big(e^{u_{j}^{(1),k}}-e^{u_{j}^{(2),k}}\big)d V_g(\eta)=O(\varepsilon^{M(m^*-2)-\delta}).
\end{align}
Combining (\ref{ineq:Theta-Theta-OUT-1}), (\ref{ineq:Theta-Theta-OUT-2}) and (\ref{ineq:Theta-Theta-OUT-3}), we get
\begin{align}\label{ineq:Theta-Theta-OUT-TOTAL1}
(u_{i}^{(1),k}&(x)-\overline{u_{i}^{(1),k}})-(u_i^{(2),k}(x)-\overline{u_{i}^{(2),k}})=O(\varepsilon^{M(m^*-2)-2-\delta}).
\end{align}
Now, by Proposition \ref{prop:OutAverageExpansion}, we get
\begin{align}\label{ineq:Theta-Theta-OUT-TOTAL2}
\overline{u_{i}^{(1),k}}-\overline{u_{i}^{(2),k}}=O(\varepsilon^{M(m^*-2)-\delta}).
\end{align}
(\ref{ineq:Theta-Theta-OUT-TOTAL1}) and (\ref{ineq:Theta-Theta-OUT-TOTAL2}) imply that
\begin{align}\label{ineq:Theta-Theta-OUT-TOTAL}
u_i^{(1),k}(x)-u_{i}^{(2),k}(x)=O(\varepsilon^{M(m^*-2)-2-\delta}).
\end{align}
(\ref{e:Theta-ThetaInner}) and (\ref{ineq:Theta-Theta-OUT-TOTAL}) complete the proof.

\noindent{$\Box$}


Based on the above consideration, in this subsection we introduce a limit problem concerning $u_i^{(1),k}-u_i^{(2),k}$. To be precise, let us define the function
\begin{align}\label{def:xiik}
\xi_i^k(x)=\frac{u_i^{(1),k}(x)-u_i^{(2),k}(x)}{\sup_{i\in I}\|u_i^{(1),k}-u_i^{(2),k}\|_{L^\infty(M)}}
\end{align}
for $i=1,\cdots,n$ and $k=1,2,\cdots$, which solves the equation
\begin{align}\label{e:xi}
-\Delta \xi_i^k=\sum_{j=1}^n a_{ij}\rho_j^k h_j c_j^k\cdot\xi_j^k=:f_i^k\mbox{ in }M
\end{align}
with $c_j^k(x)=\frac{e^{u_i^{(1),k}(x)}-e^{u_i^{(2),k}(x)}}{u_i^{(1),k}(x)-u_i^{(2),k}(x)}$. Notice that $\|\xi_i^k\|_{L^\infty(M)}\leq n$ and
\begin{align}\label{e:cExpansion}
c_i^k(x)=e^{u_i^{(2),k}(x)}(1+\sup_{i\in I}\|u_i^{(1),k}-u_i^{(2),k}\|_{L^\infty(M)})=e^{u_i^{(2),k}(x)}(1+O(\varepsilon^{M(m^*-2)-2-\delta}))
\end{align}
Now we rescale them as follows. For any $t=1,\cdots,N$, $i=1,\cdots,n$ and $k=1,\cdots$, we define
\begin{align}
&\xi_{it}^k(y)=\xi_{i}^k(\varepsilon_t^{(2),k}y+q_t^{(2),k}),\nonumber\\
&c_{it}^k(y)=c_{i}^k(\varepsilon_t^{(2),k}y+q_t^{(2),k}),\nonumber\\
&f_{it}^k(y)=f_{i}^k(\varepsilon_t^{(2),k}y+q_t^{(2),k}).\nonumber
\end{align}
Notice that $\xi_{it}^k(y)$ satisfies the equation
\begin{align}\label{e:xiScaled}
-\Delta_y\xi_{it}^k(y)=(\varepsilon_t^{(2),k})^2\sum_{j=1}^n a_{ij}\rho_j^k h_j^k(\varepsilon_t^{(2),k}y+q_t^{(2),k})c_j^l(\varepsilon_t^{(2),k}y+q_t^{(2),k})\xi_j^k(\varepsilon_t^{(2),k}y+q_t^{(2),k})\mbox{ for }|y|\leq\tau(\varepsilon_t^{(2),k})^{-1}.
\end{align}

Our first observation is that

\begin{lemma}\label{l:GlobalCancelation}
It holds that $\sum_{j=1}^na_{ij}\rho_j^k\int_M c_j^k\cdot\xi_j^k=0$.
\end{lemma}

Then, we get

\begin{lemma}\label{l:Limitbphi}
Under the above assumptions, there exist constants $b_{0,t},b_{1,t},b_{2,t}$ such that
\begin{align}
\xi_{it}^k(y)\to b_{0,t}Z_{0,i}(y)+b_{1,t}Z_{1,i}(y) +b_{2,t}Z_{2,i}(y)\mbox{ in }C_{loc}(\mathbb{R}^2)\nonumber
\end{align}
for $t=1,\cdots,N$ and $i=1,\cdots,n$.
Here, the functions $Z_{0,i}$, $Z_{1,i}$ and $Z_{2,i}$ for $i=1,\cdots,n$ are defined as in Proposition \ref{p:Kernel}.
\end{lemma}
\noindent{\bf Proof.}
By (\ref{e:cExpansion}), we get
\begin{align}
(\varepsilon_t^{(2),k})^2& c_i^k(q_t^{(2),k}+\varepsilon_t^{(2),k}y)\nonumber\\
&=e^{V_{it}^{(2),k}(y)+\phi_{it}^{(2),k}((\varepsilon_{it}^{(2),k})^{-1}y)+w_{it}^{(2),k}(y)+f_{it}^{(2),k}(\varepsilon_{it}^{(2),k}y)}(1+O(\varepsilon^{M(m^*-2)-2-\delta}+\varepsilon^{M-3-\delta}))\nonumber\\
&\to e^{V_i(y)}\mbox{ in }C_{loc}(\mathbb{R}^2)\nonumber
\end{align}
Plugging this into (\ref{e:xiScaled}), we know that $\xi_{it}^k$ converge to the bounded solution $(\xi^*_{1t},\cdots,\xi^*_{nt})$ of
\begin{align}
-\Delta_y \xi_{it}^*(y)=\sum_{j=1}^n a_{ij}\rho_j^* h_j(p_t^*)e^{V_j(y)}\xi_{it}^*(y)\mbox{ in }\mathbb{R}^2\nonumber
\end{align}
in $C_{loc}(\mathbb{R}^2)$-sense. Notice that
\begin{align}
-\Delta_y \Big(V_i(y)+\ln(\rho_i^* h_i(p_t^*))\Big)=\sum_{j=1}^n a_{ij} e^{V_j(y)+\ln(\rho_j^* h_j(p_t^*))}\mbox{ in }\mathbb{R}^2.\nonumber
\end{align}
The result follows Proposition \ref{p:Kernel} immediately.

\noindent{$\Box$}

\begin{lemma}\label{l:Limitb0}
Under the above assumptions, there exists a constant $B_0$ such that
\begin{align}
\xi_{i}^k\to B^0_i \mbox{ in }C_{loc}(M\backslash\{q^*,\cdots,q_N^*\})\nonumber
\end{align}
Moreover,
for any $i=1,\cdots,n$ and any $t=1,\cdots,N$. Moreover,
$B_i^0=(2-m_i^*)b_{0t}$
for $i=1,\cdots,n$ and $t=1,\cdots,N$.
\end{lemma}
\noindent{\bf Proof.}
By the definition (\ref{e:cExpansion}) of $c_i^k$, we know that
$c_i^k\to 0\mbox{ in }C_{loc}(M\backslash\{q_1^*,\cdots,q_N^*\})$.
On the other hand, $|\xi_i^k|\leq 1$ for any $i=1,\cdots,n$ and any $k=1,2,\cdots$, we get that
$\xi_i^k$ converges to a function $\xi_i^0$ satisfying
$\Delta \xi_i^0=0$ in $M\backslash\{q_1^*,\cdots,q_N^*\}$.
Since $|\xi_i^k(x)|\leq 1$ for any $x\in M$, any $i=1,\cdots,n$ and any $k=1,2,\cdots$, we know that $\xi_i^0\equiv B_i^0$ in $M$ and hence
\begin{align}\label{limit:xitoB0}
    \xi_i^k\to B_i^0\mbox{ in }C_{loc}(M\backslash\{q_1^*,\cdots,q_N^*\}).
\end{align}

Now we prove the relation between $b_{0,t}$ and $B_i^0$. To do so, let us introduce the functions $\zeta_{it}^k$ for any $i=1,\cdots,n$, $t=1,\cdots,N$ and $k=1,2,\cdots$ satisfying
\begin{align}
-\Delta \zeta_{it}^k=\sum_{j=1}^n a_{ij}\rho_j^k h_j(\overline{p}_t^{(2),k})e^{V_{jt}^{(2),k}(\frac{\cdot}{\varepsilon_{t}^{(2),k}})+2\ln\frac{1}{\varepsilon_t^{(2),k}}}\zeta_{it}^k\mbox{ in }B_\tau(q_t^{(2),k}).\nonumber
\end{align}
Rewriting the above and (\ref{e:xi}) into
\begin{align}
-\Delta \zeta_{t}^{ik}= \rho_i^k h_i(q_t^{(2),k})e^{V_{it}^{(2),k}(\frac{\cdot}{\varepsilon_{t}^{(2),k}})+2\ln\frac{1}{\varepsilon_t^{(2),k}}}\zeta_{it}^k\mbox{ in }B_\tau(q_t^{(2),k})\nonumber
\end{align}
and
\begin{align}
-\Delta \xi^{ik}=\rho_i^k h_i c_i^k\cdot\xi_i^k.\nonumber
\end{align}
Here, $\zeta_t^{ik}=\sum_{j=1}^n a^{ij}\zeta_{jt}^k$ and $\xi_t^{ik}=\sum_{j=1}^n a^{ij}\xi_{jt}^k$.
Combining these and divergence theorem, we get
\begin{align}
\int_{\partial B_\tau(q_t^{(2),k})}&\Big(\zeta_t^{ik}\frac{\partial\xi^{ik}}{\partial n}  -\xi^{ik}\frac{\partial\zeta_t^{ik}}{\partial n}\Big)d\mathcal{H}^1=\int_{B_\tau(q_t^{(2),k})}\Big(\zeta_t^{ik}\Delta\xi^{ik}-\xi^{ik}\Delta\zeta_t^{ik}\Big)dx\nonumber\\
&=\int_{B_r(p_t^{(2),k})}\rho_i^k\zeta_{it}^k\xi^k\Big[h_i(p_t^{(1),k})e^{V_{it}^{(2),k}(\frac{\cdot}{\varepsilon_{t}^{(2),k}})+2\ln\frac{1}{\varepsilon_t^{(1),k}}}-h_j e^{u_i^{(2),k}}(1+O(\|u^{(1),k}-u^{(2),k}\|_{L^\infty(M)}))\Big]dx\nonumber\\
&=O(\varepsilon^{m^*-2}+\varepsilon)\nonumber
\end{align}
Here, we use (\ref{e:cExpansion}).
We define
\begin{align}
(\xi^{ik})^*(r)=\frac{1}{|\partial B_R(p_t^{(1),k})|_1}\int_{\partial B_R(p_t^{(1),k})}\xi^{ik}d\mathcal{H}^1,\nonumber
\end{align}
which satisfies
\begin{align}\label{e:zizetaRadial}
\frac{d}{dr}(\xi^{ik})^*\zeta_{t}^{ik}-(\xi^{ik})^*\frac{d}{dr}\zeta_{t}^{ik}=\frac{1}{r}O(\varepsilon^{m^*-2}+\varepsilon)
\end{align}
for $r\in(R\varepsilon,\tau)$.
By Proposition \ref{p:Kernel}, we get
\begin{align}
\zeta_t^{ik}(r)&=\sum_{j=1}^n a^{ij}(2-m_j^*)+\frac{e^{I_i}(\varepsilon_t^{(1),k})^{m^*_i-2}}{m_i^* -2}r^{2-m_i^*}+O(\varepsilon^{2m^*-3}r^{3-2m^*}),\nonumber\\
\frac{d}{dr}\zeta_t^{ik}(r)&=\frac{O(\varepsilon^{m_i^*-2})}{r^{m_i^*-1}}\nonumber
\end{align}
for $r\in(R\varepsilon,\tau)$.
Then, we get
\begin{align}
\frac{d}{dr}(\xi_t^{ik})^*(r)\leq C\Bigg(\frac{\varepsilon^{m^*-2}+\varepsilon}{r}+\frac{\varepsilon^{m^*-2}}{r^{m^*-1}} \Bigg).\nonumber
\end{align}
Integrating over $(R\varepsilon_t^{(1),k},\tau)$, we find
\begin{align}
(\xi_t^{ik})^*(r)=(\xi_t^{ik})^*(R\varepsilon_t^{(1),k})+o_\varepsilon(1)+o_R(1)\nonumber
\end{align}
with $\lim_{R\to+\infty}o_R(1)=\lim_{\varepsilon\to0}o_\varepsilon(1)=0$. Notice that $\varepsilon_k<< R$. By Lemma \ref{l:Limitbphi} and (\ref{limit:xitoB0}), we get
$\sum_{j=1}^n a^{ij} B_j^0 =\sum_{j=1}^n a^{ij}(2-m_j^*)b_{0t}+o(1$),
which implies that
$B_i^0=(2-m_i^*)b_{0t}$
for $i=1,\cdots,n$ and $t=1,\cdots,N$.

\noindent{$\Box$}


\section{Proof of Theorem \ref{t:MAIN} and Theorem \ref{non-deg}}

By Lemma \ref{l:Limitbphi}, to prove the uniqueness, we only need to check that $b_{0t}=b_{1t}=b_{2t}=0$. In this part, we prove that $b_{0t}=0$.

\subsection{A preliminary expansion of the difference between $u^{(1),k}$ and $u^{(2),k}$}\label{Subsection:b01}

The aim of this section is to prove that
$b_{0t}\equiv 0$
for any $t=1,\cdots,N$. Using the symbols from the previous section, we proceed a analysis for the function $\xi^k=(\xi_1^k,\cdots,\xi_n^k)$. Recall that 
\begin{align}\label{e:xiik}
-\Delta_y\xi_{it}^k(y)=(\varepsilon_t^{(2),k})^2\sum_{j=1}^n a_{ij}\rho_j^k h_j^k(\varepsilon_t^{(2),k}y+q_t^{(2),k})c_j^k(\varepsilon_t^{(2),k}y+q_t^{(2),k})\xi_j^k(\varepsilon_t^{(2),k}y+q_t^{(2),k})\mbox{ for }|y|\leq\tau(\varepsilon_t^{(2),k})^{-1}.
\end{align}
In the following, we choose a blowup point to be the origin of the coordinate system. {From now on, for simplicity, we denote the bubbling disk around $p_t^*$ by $B_t$ for $t=1,\cdots,N$. $\lambda B_t$ denotes the ball with the same center with $B_t$ and $\lambda$ times radii. Denote $\sigma_k=\sup_i\|u_i^{(1),k}-u_i^{(2),k}\|$. We only study the case $m^*<4$ since $m^*=4$ is similar.}

~~

\noindent{\bf Step 1. The leading term of $\xi^k$}

~~

By a similar observation as in \cite{BartJevLeeYang2019}, for a bubbling disk $B_t$ for any $t=1,\cdots,N$, we choose 
\begin{align}\label{def:Xik}
\Xi^k(y)=\sum_{l=0}^2 b_{lt}^k Z^k_l
\end{align}
with suitable $b_{lt}^k$ such that
\begin{align}\label{ineq:b0b1b2Projection}
\sum_{i=1}^n\int_{(\varepsilon_k^{(2),k)})^{-1}\cdot B_t}e^{V_{it}^{(2),k}}Z_{i,l}^k(\xi_i^k-b_{lt}^k Z_{i,l}^k)dy=0
\end{align}
for $l=0,1,2$ and $t=1,\cdots,N$. It is clear that
\begin{align}\label{e:Xiik}
\Delta\Xi_i^k+\sum_j a_{ij}h_j^k(q_t^{(2),k})e^{V_{t,j}^{(2),k}}\Xi_j^k=0.
\end{align}
Evidently, $b_0^k$, $b_1^k$ and $b_2^k$ are bounded. Moreover, by Lemma \ref{l:Limitbphi}, we get
\begin{align}\label{e:bokLimit}
b_{0t}^k\to B_i^0(2-m_i^*)^{-1}
\end{align}
due to the local convergence of $\xi_i^k$ and Lemma \ref{l:Limitb0}.

Next, we define harmonic functions $\psi_{\xi,i}^k$ encoding the boundary oscillation of $\xi_i^k-\Xi_i^k$:
\begin{align}\label{def:psixi}
\Delta\psi_{\xi,i}^k=0\mbox{ in }B_t\mbox{, }\psi_{\xi,i}^k|_{\partial B_t}=\xi_i^k-\Xi_i^k-\frac{1}{2\pi\tau}\int_{\partial B_t}(\xi_i^k-\Xi_i^k) d\mathcal{H}^1
\end{align}
for any $i=1,\cdots,n$.
Due to Lemma \ref{l:Limitb0}, we get a rough estimate
\begin{align}
|\psi_{\xi,i}^k(\varepsilon_k y)|\leq o(\varepsilon_t^{(2),k})|y|.\nonumber
\end{align}
Denoting
\begin{align}
w_{\xi,i}^k(y):=\xi_i^k(\varepsilon_t^{(2),k}y+q_t^{(2),k})-\Xi_i^k(y)-\psi_{\xi,i}^k(\varepsilon_k^{(2),k}y),\nonumber
\end{align}
we will later derive a frequency analysis on $w_{\xi,i}^k$.
It is evident that
\begin{align}
w_{\xi,i}^k|_{\partial ((\varepsilon_t^{(2),k})^{-1}\cdot B_t)}=constant,\nonumber
\end{align}
\begin{align}\label{e:EZWERROR}
\sum_{i=1}^n\int_{(\varepsilon_t^{(2),k})^{-1}\cdot B_t}e^{V_{it}^{(2),k}}Z_{i,l}^k w_{\xi,i}^k=O(\varepsilon)
\end{align}
for $l=0,1,2$ and
\begin{align}
&\Delta \sum_{j=1}^n \Big(a^{ij}w_{\xi,j}^k\Big)+(\varepsilon_t^{(2),k})^2h_i c_i^kw^k_{\xi,i}\nonumber\\
&=h_i^k(q_t^{(2),k})\Bigg(e^{V_i^{(2),k}}-\frac{h_i^k(\varepsilon_t^{(2),k}y+q_t^{(2),k})}{h_i^k(p_t^{(2),k})}(\varepsilon_t^{(2),k})^2 c_i(\varepsilon_t^{(2),k}y)\Bigg)\Xi_i^k-h_i^k(\varepsilon_t^{(2),k}y+q_t^{(2),k})(\varepsilon_t^{(2),k})^2 c_i\psi_{\xi,i}^k\nonumber\\
&=O(\varepsilon_t^{(2),k}(1+|y|)^{1-m^*{+\delta}}).
\end{align}
Notice that (\ref{e:EZWERROR}) is due to
\begin{align}
\sum_{i=1}^n\int_{(\varepsilon_t^{(2),k})^{-1}\cdot B_t}e^{V_{it}^{(2),k}}Z_{i,0}^k w_{\xi,i}^k=\sum_{i=1}^n\int_{(\varepsilon_t^{(2),k})^{-1}\cdot B_t}e^{V_{it}^{(2),k}}Z_{i,0}^k \psi_{\xi,i}^k=0\nonumber
\end{align}
by (\ref{e:Xiik}) and its symmetry and
\begin{align}
\sum_{i=1}^n\int_{(\varepsilon_t^{(2),k})^{-1}\cdot B_t}e^{V_{it}^{(2),k}}Z_{i,1}^k w_{\xi,i}^k&=\sum_{i=1}^n\int_{(\varepsilon_t^{(2),k})^{-1}\cdot B_t}e^{V_{it}^{(2),k}}Z_{i,1}^k \psi_{\xi,i}^k\nonumber\\
&=o\Bigg(\varepsilon\sum_{i=1}^n\int_{(\varepsilon_t^{(2),k})^{-1}\cdot B_t}|\eta|e^{V_{it}^{(2),k}}Z_{i,1}^kd\eta\Bigg)\nonumber\\
&=o(\varepsilon).\nonumber
\end{align}
By a potential analysis argument, we get the following rough estimate for $w_{\xi,i}^k$:
\begin{equation}\label{ineq:wxiik}
|w_{\xi,i}^k(y)|=\left\{
\begin{array}{lr}
O(\varepsilon^{(2),k}_k(1+|y|)^\delta)\mbox{ if }m^*\geq3;\nonumber\\
O(\varepsilon^{(2),k}_k(1+|y|)^{3-m^*+\delta})\mbox{ if }m^*< 3.\nonumber
\end{array}
\right.
\end{equation}
for small $\delta>0$.

~~

\noindent{\bf Step 2. Oscillation of $\xi_i^k$ outside the bubbling disks}

Now we assume that
\begin{align}
b_{0t}^k\to b_{0t}\neq 0.\nonumber
\end{align}
An elementary observation is that $b_{0t}-b_{0s}=o(1)$ for $t,s=1,\cdots,N$ and $t\neq s$.
In follows, we give a finer estimate of the oscillation of $\xi_i^k$ outside of the bubbling disks, which lead us to a better bound for $\psi_{\xi,i}^k$.

\begin{lemma}\label{l:XiOsc}
For any $x_1,x_2\in M\backslash(\cup_{t=1}^N B_t)$, we get
\begin{align}
\xi_i^k(x_1)-\xi_i^k(x_2)=O\big(\varepsilon^{m^*-2-\delta}+\varepsilon^{1-\delta}\big).\nonumber
\end{align}
\end{lemma}
\noindent{\bf Proof.}
Denote $\overline{\xi_i^k}=\int_M\xi_i^kdx$. 
By Green's representation, we get
\begin{align}
\sum_{j=1}^n a^{ij}\xi_j^k(x)&=\sum_{j=1}^n a^{ij}\overline{\xi_j^k}+\int_M G(x,\eta)\rho_i^k h_ic_i^k(\eta)\xi_i^k(\eta)d\eta\nonumber\\
&=\sum_{j=1}^n a^{ij}\overline{\xi_j^k}+\sum_{t=1}^N\int_{B^M_{\tau/2}(p_t^{(2),k})} G(x,\eta)\rho_i^k h_ic_i^k(\eta)\xi_i^k(\eta)d\eta\nonumber\\
&\quad\quad+ \int_{M\backslash\cup_{t=1}B^M_{\tau/2}(p_t^{(2),k})} G(x,\eta)\rho_i^k h_ic_i^k(\eta)\xi_i^k(\eta)d\eta.\nonumber
\end{align}
It is clear that
\begin{align}\label{ineq:xiOut}
\int_{M\backslash\cup_{t=1}B^M_{\tau/2}(p_t^{(2),k})} G(x,\eta)\rho_i^k h_ic_i^k(\eta)\xi_i^k(\eta)d\eta=O(\varepsilon^{m^*-2{-\delta}}).
\end{align}
For any $t=1,\cdots,N$, we get
\begin{align}
&\int_{B^M_{\tau/2}(p_t^{(2),k})} G(x,\eta)\rho_i^k h_ic_i^k(\eta)\xi_i^k(\eta)d\eta\nonumber\\
&=\int_{B^M_{\tau/2}(p_t^{(2),k})} (G(x,\eta)-G(x,p_t^{(2),k})\rho_i^k h_ic_i^k(\eta)\xi_i^k(\eta)d\eta+G(x,p_t^{(2),k})\int_{B^M_{\tau/2}(p_t^{(2),k})}\rho_i^k h_ic_i^k(\eta)\xi_i^k(\eta)d\eta.\nonumber
\end{align}
Here,
\begin{align}
\int_{B(p_t^{(2),k},\tau/2)} (G(x,\eta)-G(x,p_t^{(2),k})\rho_i^k h_ic_i^k(\eta)\xi_i^k(\eta)d\eta\leq C\varepsilon\int_{B_{C(\varepsilon_t^{(2),k})^{-1}}(0)}|\eta|^{1-m^*}d\eta=O(\varepsilon_k^{m^*-2}+\varepsilon_k)\nonumber
\end{align}
and for $m^*<3$
\begin{align}
&\int_{B^M_{\tau/2}(p_t^{(2),k})}\rho_i^k h_ic_i^k(\eta)\xi_i^k(\eta)d\eta\nonumber\\
&=\int_{B^M_{\tau/(2\varepsilon_t^{(2),k})}(p_t^{(2),k})}\rho_i^k(h^k_i(0)+O(\varepsilon|\eta|))(e^{V_{t,j}^{(2),k}}(1+O(\varepsilon(1+|\eta|)^{3-m^*})))(\Xi_i^k+\psi_{\xi,i}^k+w_{\xi,i}^k)+O(\sigma_k)\nonumber\\
&=O(\varepsilon^{m^*-2})+\int_{B_{\tau/(2\varepsilon_t^{(2),k})}^M(p_t^{(2),k})}\rho_i^k h_i^k(0)e^{V_{t,i}^{(2),k}}(\Xi_i^k+\psi_{\xi,i}^k+w_{\xi,i}^k)+O(\sigma_k)\nonumber\\
&=O(\varepsilon^{m^*-2-\delta})\nonumber
\end{align}
while for $m^*\geq 3$,
\begin{align}
&\int_{B^M_{\tau/2}(p_t^{(2),k})}\rho_i^k h_ic_i^k(\eta)\xi_i^k(\eta)d\eta\nonumber\\
&=\int_{B^M_{\tau/(2\varepsilon_t^{(2),k})}(p_t^{(2),k})}\rho_i^k(h_i^k(0)+O(\varepsilon|\eta|))(e^{V_{t,j}^{(2),k}}(1+O(\varepsilon(1+|\eta|)^{\delta})))(\Xi_i^k+\psi_{\xi,i}^k+w_{\xi,i}^k)+O(\sigma_k)\nonumber\\
&=O(\varepsilon^{{-\delta}}).\nonumber
\end{align}

This proves the lemma.

\noindent{$\Box$}

Combining the definition of $\psi_{\xi,i}^k$ (see (\ref{def:psixi})), we get
\begin{align}\label{ineq:psiNewBdd}
|\psi_{\xi,i}^k(\varepsilon_i^{(2),k}y)|\leq C(\varepsilon^{m^*-1-\delta}+\varepsilon^{2-\delta})|y|.
\end{align}

~~

\noindent{\bf Step 3. A decomposition of $w_{\xi,i}^k$ by its frequency}

~~

Before we go to the next stage, let us recall that $w_{\xi,i}^k(y):=\xi_i^k(\varepsilon_t^{(2),k}y+q_t^{(2),k})-\Xi_i^k(y)-\psi_{\xi,i}^k(\varepsilon_k^{(2),k}y)$. Therefore, $w_{\xi,i}^k$ is constant on $\partial B_{\tau(\varepsilon_t^{(2),k})^{-1}}(0)$. By (\ref{e:Xiik}) and (\ref{e:xiik}), the corresponding equation is that
\begin{align}\label{e:wxiik}
&\Delta w_{\xi,i}^k+\sum_{j=1}^n a_{ij} h_j^k(q_t^{(2),k} ) e^{V_{t,j}^{(2),k}}w_{\xi,j}^k\nonumber\\
&=-\sum_j a_{ij}h_j^k(q_t^{(2),k} )e^{V_{t,j}^{(2),k}}\Bigg(\frac{h_j^k(\varepsilon_t^{(2),k}y+q_t^{(2),k})}{h_j^k(q_t^{(2),k})}e^{u_j^{(2),k}-V_{t,j}^{(2),k}}-1\Bigg)\xi_j^k-\sum_j a_{ij}h_j^k(q_t^{(2),k} )e^{V_{t,j}^{(2),k}}\psi_{\xi,i}^k\nonumber\\
&\quad+O(\sigma_k(1+|y|)^{-m^*{+\delta}})=:E_{\xi,i}^k.
\end{align}
Now we decompose $w_{\xi,i}^k$ by its frequency. To be precise, we define
\begin{align}\label{def:wxili}
\left\{
\begin{aligned}
&w_{\xi,0,i}^k(r):=\frac{1}{2\pi}\int_0^{2\pi}w_{\xi,l,i}^k(r\cos(\theta),r\sin(\theta))d\theta,\\
&w_{\xi,l,i,1}^k(r):=\frac{1}{2\pi}\int_0^{2\pi}w_{\xi,l,i}^k(r\cos(l\theta),r\sin(l\theta))\cdot \cos(l\theta)d\theta,\\
&w_{\xi,l,i,2}^k(r):=\frac{1}{2\pi}\int_0^{2\pi}w_{\xi,l,i}^k(r\cos(l\theta),r\sin(l\theta))\cdot \sin(l\theta)d\theta,\\
&w_{\xi,l,i}^k:=w_{\xi,l,i,1}^k(r)\cdot \cos(l\theta)+w_{\xi,l,i,2}^k(r)\cdot \sin
\end{aligned}
\right.
\end{align}
and
\begin{align}\label{def:Exili}
\left\{
\begin{aligned}
&E_{\xi,0,i}^k(r):=\frac{1}{2\pi}\int_0^{2\pi}E_{\xi,l,i}^k(r\cos(\theta),r\sin(\theta))d\theta,\\
&E_{\xi,l,i,1}^k(r):=\frac{1}{2\pi}\int_0^{2\pi}E_{\xi,l,i}^k(r\cos(l\theta),r\sin(l\theta))\cdot \cos(l\theta)d\theta,\\
&E_{\xi,l,i,2}^k(r):=\frac{1}{2\pi}\int_0^{2\pi}E_{\xi,l,i}^k(r\cos(l\theta),r\sin(l\theta))\cdot \sin(l\theta)d\theta,\\
&E_{\xi,l,i}^k:=E_{\xi,l,i,1}^k(r)\cdot \cos(l\theta)+E_{\xi,l,i,2}^k(r)\cdot \sin(l\theta)
\end{aligned}
\right.
\end{align}
for $l=1,2,\cdots$.
With these notions, we get
\begin{align}\label{e:wxilik}
\Delta w_{\xi,l,i}^k+\sum_{j=1}^n a_{ij}h_j^k(q_t^{(2),k} )e^{V_{t,j}^{(2),k}}w_{\xi,l,j}^k=E_{\xi,l,i}^k
\end{align}
for $l=0,1,2,\cdots$.
Moreover, $w_{\xi,0,i}^k|_{\partial B_{\tau(\varepsilon_k^{(2),k})^{-1}}(0)}=constant$ and $w_{\xi,l,i}^k|_{\partial B_{\tau(\varepsilon_k^{(2),k})^{-1}}(0)}=0$ for $l\geq 1$.

Among all of the $w_{\xi,l,i}^k$, the case of $l=0$ is of particular interest. Now we study it in detail.
To do this, we need to check $E_{\xi,0,i}^k$. By a direct computation, we get
\begin{align}\label{Expansion:Exiik}
E_{\xi,i}^k&=
-\sum_j a_{ij}h_j^k(q_t^{(2),k})e^{V_{t,j}^{(2),k}}\Bigg(\frac{h_j^k(\varepsilon_t^{(2),k}y+q_t^{(2),k})}{h_j^k(q_t^{(2),k})}e^{u_j^{(2),k}-V_{t,j}^{(2),k}}-1\Bigg)\xi_j^k-\sum_j a_{ij}h_j^k(q_t^{(2),k})e^{V_{t,j}^{(2),k}}\psi_{\xi,i}^k\nonumber\\
&\quad+O(\sigma_k(1+|y|)^{-m^*})\nonumber\\
&=-\sum_j a_{ij}h_j^k(q_t^{(2),k})e^{V_{t,j}^{(2),k}}\Bigg(\frac{h_j^k(\varepsilon_t^{(2),k}y+q_t^{(2),k})}{h_j^k(q_t^{(2),k})}e^{u_j^{(2),k}-V_{t,j}^{(2),k}}-e^{u_j^{(2),k}-V_{t,j}^{(2),k}}+e^{u_j^{(2),k}-V_{t,j}^{(2),k}}-1\Bigg)\times\nonumber\\
&\quad\times[b_{0t}^kZ_{0,i}^k+b_{1t}^kZ_{1,i}^k+b_{2t}^kZ_{2,i}^k+\psi_{\xi,i}^k+w_{\xi,i}^k]-\sum_j a_{ij}h_j^k(q_t^{(2),k})e^{V_{t,j}^{(2),k}}\psi_{\xi,i}^k+O(\sigma_k(1+|y|)^{-m^*{+\delta}}).
\end{align}
Taking the 0-frequency part of $E_{\xi,i}^k$, we get for $m^*<3$
\begin{align}
E_{\xi,0,i}^k=O(\varepsilon^2|b_{0t}^k|(1+|y|)^{2-m^*}+\varepsilon(1+|y|)^{-m^*{+\delta}}+\varepsilon^2(1+|y|)^{4-2m^*{+\delta}}).\nonumber
\end{align}
By a Green's function argument, we get for $m^*<3$
\begin{align}
|w_{\xi,0,i}^k(y)|=O(\varepsilon^2|b_{0t}^k|(1+|y|)^{4-m^*}+\varepsilon(1+|y|)^{2-m^*{+\delta}}+\varepsilon^2(1+|y|)^{6-2m^*{+\delta}}).\nonumber
\end{align}
Then, an ODE argument gives that for $m^*<3$
\begin{align}\label{ineq:wxi0ikm<3}
\Big|\frac{d}{dr}w_{\xi,0,i}^k(r)\Big|=O(\varepsilon^2|b_{0t}^k|(1+|y|)^{3-m^*}+\varepsilon(1+|y|)^{1-m^*{+\delta}}+\varepsilon^2(1+|y|)^{5-2m^*{+\delta}}).
\end{align}
For the case $m^*\geq3$, we get
\begin{align}\label{ineq:wxi0ikm>3}
\Big|\frac{d}{dr}w_{\xi,0,i}^k(r)\Big|=O(\varepsilon^2|b_{0t}^k|(1+|y|)^{3-m^*}+\varepsilon(1+|y|)^{1-m^*{+\delta}}).
\end{align}


~~

\noindent{\bf Step 4. The expansion of an integral}

~~

By a direct computation, we get
\begin{align}
\int_{B_{\tau/\varepsilon_t^{(2),k}}(0)}\sum_j a_{ij}\rho^k_jh_j^kc_j^k\xi_j^kd\eta
&=-\int_{B_{\tau/\varepsilon_t^{(2),k})}(0)}\Delta\xi_i^k d\eta=\int_{\partial B_{\tau/\varepsilon_t^{(2),k}}(0)}\frac{\partial \xi_i^k}{\partial{\bf n}}d\mathcal{H}_\eta^1\nonumber\\
&=b_0^k \int_{\partial B_{\tau/\varepsilon_t^{(2),k}}(0)}\frac{\partial Z_{0,i}^k}{\partial{\bf n}}d\mathcal{H}_\eta^1+ \int_{\partial B_{\tau/\varepsilon_t^{(2),k}}(0)} \frac{\partial w_{\xi,0,i}^k}{\partial{\bf n}}d\mathcal{H}_\eta^1.
\end{align}
Here, we use (\ref{e:xiik}), (\ref{def:psixi}), (\ref{def:Xik}) and (\ref{def:wxili}). By a direct computation, we get
\begin{align}
\frac{d}{dr}Z_{0,i}^k(r)=-\sum_{j=1}^n a_{ij}e^{I_j}r^{1-m^*_j}+O(r^{3-2m^*{+\delta}}).\nonumber
\end{align}
Then,
\begin{align}
b_{0t}^k \int_{\partial B_{\tau/\varepsilon_t^{(2),k}}(0)}\frac{\partial Z_{0,i}^k}{\partial{\bf n}}d\mathcal{H}_\eta^1=2\pi b_{0t}^k\sum_{j=1}^n a_{ij}e^{I_j}\Big(\frac{\varepsilon_t^{(2),k}}{\tau}\Big)^{m_j^*-2}+O(\varepsilon^{2m^*-4{-\delta}})\nonumber
\end{align}
and
\begin{align}
\int_{\partial B_{\tau/\varepsilon_t^{(2),k}}(0)} \frac{\partial w_{\xi,0,i}^k}{\partial{\bf n}}d\mathcal{H}_\eta^1=O(|b_{0t}^k|\tau^{4-m^*}\varepsilon^{m^*-2}+\varepsilon^{m^*-1{-\delta}}+\varepsilon^{2m^*-4{-\delta}})\nonumber
\end{align}
for some small $\tau$.
{
An expansion can be given for the $m^*=4$ case as well.}

The latter one is due to (\ref{ineq:wxi0ikm<3}) and (\ref{ineq:wxi0ikm>3}).
Therefore, we get
\begin{align}\label{ineq:XiInBall}
\int_{B_{\tau/\varepsilon_t^{(2),k}}(0)}\sum_j a_{ij}\rho^k_jh_j^kc_j^k\xi_j^kd\eta&=2\pi b_0^k\sum_{j=1}^n a_{ij}e^{I_j}\Big(\frac{\varepsilon_t^{(2),k}}{\tau}\Big)^{m_j^*-2}+O(|b_0^k|\tau^{4-m^*}\varepsilon^{m^*-2}+\varepsilon^{m^*-1{-\delta}}+\varepsilon^{2m^*-4{-\delta}}).
\end{align}

\subsection{Analysis of $b_{0t}$ for $t=1,\cdots,N$}\label{Subsection:b02}

\begin{claim}
Under the assumptions of Theorem \ref{t:MAIN}, we get $b_0^k\to0$.
\end{claim}

\noindent{\bf Proof.}
Let us assume $b_0\neq 0$ and consider the $m^*<4$ case.
To complete the proof, we estimate 
\begin{align}\label{ineq:XiOutBall}
&\int_{M\backslash \cup_{t=1}^N B_t}  \sum_j a_{ij}\rho^k_jh_j^kc_j^k\xi_j^kd\eta\nonumber\\
&=b_{0t}^k\sum_{j=1}^n a_{ij}e^{I_i}\sum_{t=1}^N(\varepsilon_t^k)^{m_i^k-2}\int_{\Omega_t\backslash B^M_\tau(q_t^{(2),k})} \frac{h_i^k(x)}{h_i^k(q_t^{(2),k})}e^{2\pi m_i^k\sum_{l=1}^N[G(x,p_l^k)-G^*(p_l^k,p_t^k)]}dV_g(1+o(1)).
\end{align}
Therefore, combining Lemma \ref{l:GlobalCancelation}, (\ref{ineq:XiInBall}) and (\ref{ineq:XiOutBall}), we get
\begin{align}\label{e:DXi}
0&=\int_{M}\sum_j a_{ij}\rho^k_jh_j^kc_j^k\xi_j^kd\eta\nonumber\\
&=\sum_{t=1}^N\int_{B_{\tau/\varepsilon_t^{(2),k}}(0)}\sum_j a_{ij}\rho^k_jh_j^kc_j^k\xi_j^kd\eta+\int_{M\backslash \cup_{t=1}^N B_t}\sum_j a_{ij}\rho^k_jh_j^kc_j^k\xi_j^kd\eta\nonumber\\
&=b_{0t}^k\Bigg(\sum_{j=1}^N a_{ij} e^{I_j} D_{jt}\frac{e^{H_{it}(p^{(2),k})}}{e^{H_{i1}(p^{(2),k})}}+o_\tau(1)\Bigg)(\varepsilon_1^{(2),k})^{m^*-2}(1+o(1)).
\end{align}
By the assumption, we see that $b_{0t}^k\to b_0=0$ for any $t=1,\cdots,N$. 

A similar argument holds for the $m^*=4$ case.

\noindent{$\Box$}

\subsection{Higher order expansions via induction}\label{Subsection:b121}

In this part, we prove that
$b_{1t}=b_{2t}=0$
for any $t=1,\cdots,N$. To do this, we assume that
\begin{align}
(b_1^k,b_2^k)\to(b_1,b_2)\neq(0,0).\nonumber
\end{align}

In this subsection, we always assume that $m^*<3$.

~~

\noindent{\bf Step 1. Initial estimate for the vanishing orders.}

~~

Matching the expansion of $\xi_i^k$ inside and outside of bubbling disks, we get
\begin{align}
\xi_i^k(x)=\overline{\xi}_i^k+O(\varepsilon^{m^*-2})=b_{0t}^kZ_{0,i}+b_{1t}^k Z_{1,i}+b_{2t}^{k}Z_{2,i}|_{\partial B_\frac{\tau}{\varepsilon}}+\psi_{\xi,i}^k+w_{\xi,i}^k
|_{B_\frac{\tau}{\varepsilon}}.\nonumber
\end{align}
Therefore, we get
\begin{align}\label{Expansion:XiOut02}
\xi_i^k(x)|_{M\backslash \cup_t B_t}=b_{t0}^k Z_{0,i}|_{\partial ((\varepsilon_t^{(2),k})^{-1}\cdot B_t)}+O(\varepsilon^{m^*-2})=(2-m^*)b_{0t}^k+O(\varepsilon^{m^*-2{-\delta}}).
\end{align}
Now we prove an vanishing rate of $b_{t0}^k$ as follows.
\begin{claim}\label{c:b0Vanishing02}
$|b_{t0}^k|=O(\varepsilon^{m^*-2{-\delta}})$.
\end{claim}
\noindent{\bf Proof.}
Plugging (\ref{Expansion:XiOut02}) into (\ref{ineq:XiOutBall}), we get
\begin{align}
&\int_{M\backslash \cup_{t=1}^N B_t} \sum_j a_{ij}\rho^k_jh_j^kc_j^k\xi_j^kd\eta\nonumber\\
&=b_{01}^k\sum_{j=1}^n a_{ij}e^{I_i}\sum_{t=1}^N(\varepsilon_t^k)^{m_i^k-2}\int_{\Omega_t\backslash B_t} \frac{h_i^k(x)}{h_i^k(p_t^k)}e^{2\pi m_i^k\sum_{l=1}^N[G(x,p_l^k)-G^*(p_l^k,p_t^k)]}dV_g(1+o(1))+O(\varepsilon^{2m^*-4{-\delta}}).\nonumber
\end{align}
Combining (\ref{ineq:XiInBall}), a similar computation as (\ref{e:DXi}), we get $|b_{t0}^k|=O(\varepsilon^{m^*-2})$ by { the assumption (\ref{nature-a})}.
$\Box$

In follows we prove a higher order expansion of $\xi_i^k$ outside the bubbling disks. 
\begin{claim}\label{c:ExpansionXiOut02}
It holds for any $x\in M\backslash\cup_t B_t$ that
\begin{align}
\xi_i^k(x)=\overline{\xi}_i^k+\overline{T^{\mbox{in}}_{i,k,2}}(x)+(2-m_i^k)b_{0t}^k\overline{T^{\mbox{out}}_{i,k,2}}\varepsilon^{m_i^k-2}+O(\varepsilon^{2m^*-4{-\delta}}).\nonumber
\end{align}
Moreover,
\begin{align}
\overline{T^{\mbox{in}}_{i,k,2}}(x)=O(\varepsilon^{m^*-2}|b_{0}^k|+\varepsilon(|b_1^k|+|b_2^k|))=O(\varepsilon^{2m^*-4{-\delta}}+\varepsilon)\nonumber
\end{align}
and
\begin{align}
(2-m_i^k)b_{0t}^k\overline{T^{\mbox{out}}_{i,k,2}}\varepsilon^{m_i^k-2}=O(\varepsilon^{m^*-2}|b_0^k|)=O(\varepsilon^{2m^*-4{-\delta}}).\nonumber
\end{align}
Here, we use $\overline{T^{\mbox{in}}_{i,k,2}}(x)$ and $\overline{T^{\mbox{out}}_{i,k,2}}$ to denote the higher order expansions.
\end{claim}
\noindent{\bf Proof.}
Using the Green's representation
\begin{align}
\sum_{j=1}^n a^{ij}\xi_j^k(x)&=\sum_{j=1}^n a^{ij}\overline{\xi_j^k}+\int_M G(x,\eta)\rho_i^k h_ic_i^k(\eta)\xi_i^k(\eta)d\eta\nonumber\\
&=\sum_{j=1}^n a^{ij}\overline{\xi_j^k}+\sum_{t=1}^N\int_{B_{\tau/2}(p_t^{(2),k})} G(x,\eta)\rho_i^k h_ic_i^k(\eta)\xi_i^k(\eta)d\eta\nonumber\\
&\quad+ \int_{M\backslash\cup_{t=1}B_{\tau/2}(p_t^{(2),k})} G(x,\eta)\rho_i^k h_ic_i^k(\eta)\xi_i^k(\eta)d\eta,\nonumber
\end{align}
we get by (\ref{Expansion:XiOut02}) that
\begin{align}
\int_{M\backslash\cup_t B_{\tau/2}(p_t^{(2),k})}&G(x,\eta)\rho_i^k h_i(\eta) c_i^k(\eta)\xi_i^k(\eta)d\eta\nonumber\\
&=\int_{M\backslash\cup_t B_{\tau/2}(p_t^{(2),k})}G(x,\eta)\rho_i^k h_i(\eta)e^{V_{it}^{(2),k}+c_{M,k}^{(2)}}b_{0t}^k (2-m^k_i)+O(\varepsilon^{2m^*-4{-\delta}})\nonumber\\
&=:(2-m_i^k)b_{0t}^k\overline{T^{\mbox{out}}_{i,k,2}}\varepsilon^{m_i^k-2}+O(\varepsilon^{2m^*-4{-\delta}}).\nonumber
\end{align}
On the other hand, we get for the integral inside the bubbling disks
\begin{align}
\int_{B_{\tau/2}(p_t^{(2),k})}&G(x,\eta)\rho_i^k h_i c_i^k(\eta)\xi_i^k(\eta)d\eta\nonumber\\
&=\int_{B_{\tau/2}(p_t^{(2),k})}G(x,\eta)\rho_i^k h_i e^{V_{it}^{(2),k}+c_{M,k}^{(2)}}(b_{0t}Z_{0,i}+b_{1t}Z_{1,i}+b_{2t}Z_{2,i})\nonumber\\
&\quad\quad+O\Bigg(\int_{B_\frac{\tau}{\varepsilon}}(1+|y|)^{-m^*}\psi_{\xi,i}^k w_{i}^k+\int_{B_\frac{\tau}{\varepsilon}}(1+|y|)^{-m^*}w_{\xi,i}^k w_i^k\Bigg)\nonumber\\
&=\int_{B_{\tau/2}(p_t^{(2),k})}G(x,\eta)\rho_i^k h_i e^{V_{it}^{(2),k}+c_{M,k}^{(2)}}(b_{0t}Z_{0,i}+b_{1t}Z_{1,i}+b_{2t}Z_{2,i})+O(\varepsilon^{2m^*-4{-\delta}})\nonumber\\
&=:\overline{T^{\mbox{in}}_{i,k,2}}(x)+O(\varepsilon^{2m^*-4{-\delta}}).\nonumber
\end{align}
Therefore, the expansion of $\xi_i^k$ outside of the bubbling disks reads
\begin{align}
\xi_i^k(x)=\overline{\xi}_i^k+\overline{T^{\mbox{in}}_{i,k,2}}(x)+(2-m_i^k)b_{0t}^k\overline{T^{\mbox{out}}_{i,k,2}}\varepsilon^{m_i^k-2}+O(\varepsilon^{2m^*-4{-\delta}}).\nonumber
\end{align}
By a direct computation, we get
\begin{align}
\overline{T^{\mbox{in}}_{i,k,2}}(x)=O(\varepsilon^{m^*-2}|b_{0}^k|+\varepsilon(|b_1^k|+|b_2^k|))=O(\varepsilon^{2m^*-4{-\delta}}+\varepsilon)\nonumber
\end{align}
and
\begin{align}
(2-m_i^k)b_{0t}^k\overline{T^{\mbox{out}}_{i,k,2}}\varepsilon^{m_i^k-2}=O(\varepsilon^{m^*-2}|b_0^k|)=O(\varepsilon^{2m^*-4{-\delta}}).\nonumber
\end{align}
\noindent{$\Box$}

An immediate corollary is about an improved estimate for $\psi_{\xi,i}^k$.
\begin{corollary}
It holds that $\sup_{x_1,x_2\in M\backslash \cup_t B_t}\xi_i^k(x_1)-\xi_i^k(x_2)=O(\varepsilon^{2m^*-4{-\delta}}+\varepsilon)$. Therefore,
\begin{align}\label{ineq:PhiXi02}
|\psi_{\xi,i}^k(y)|\leq C(\varepsilon^{2m^*-3{-\delta}}+\varepsilon^2)|y|
\end{align}
for $y\in B_{\tau/\varepsilon}$.
\end{corollary}

Now we have an improved decay of $b_0^k$:
\begin{claim}\label{c:b0Vanishing03}
It holds that $|b_{0t}^k|=O(\varepsilon+\varepsilon^{2m^*-4})$.
\end{claim}
\noindent{\bf Proof.}
With the help of (\ref{ineq:PhiXi02}), we get for the 0-frequency part of $E_{\xi,i}^k$
\begin{align}
E_{\xi,0,i}^k=O(\varepsilon^2(1+|y|)^{2-m^*}+\varepsilon(1+|y|)^{-m^*{+\delta}}+\varepsilon^{m^*-1}(1+|y|)^{3-2m^*{+\delta}}+(\varepsilon^2+\varepsilon^{2m^*-3})(1+|y|)^{1-m^*{+\delta}}).\nonumber
\end{align}
Here, $E_{\xi,i}^k$ is as (\ref{Expansion:Exiik}). It follows that
\begin{align}
|w_{\xi,0,i}^k(y)|=O(\varepsilon^2|b_{0t}^k|(1+|y|)^{4-m^*}+\varepsilon(1+|y|)^{2-m^*{+\delta}}+\varepsilon^{m^*-1}(1+|y|)^{5-2m^*{+\delta}}+(\varepsilon^2+\varepsilon^{2m^*-3})(1+|y|)^{3-m^*{+\delta}})\nonumber
\end{align}
and
\begin{align}
\int_{\partial B_{\tau/\varepsilon_t^{(2),k}}(0)} \frac{\partial w_{\xi,0,i}^k}{\partial{\bf n}}d\mathcal{H}_\eta^1=O(|b_{0t}^k|\varepsilon^{m^*-2})+O(\varepsilon^{m^*-1{-\delta}}+\varepsilon^{3m^*-6{-\delta}}).\nonumber
\end{align}
On the other hand, following (\ref{c:ExpansionXiOut02}) and similar computations as in Claim \ref{c:b0Vanishing02}, we get
$|b_{0t}^k|=O(\varepsilon+\varepsilon^{2m^*-4}$).

\noindent{$\Box$}

A routine computation gives

\begin{corollary}\label{coro:wXi02}
For $|y|\leq\tau/\varepsilon$,
\begin{align}
|w_{\xi,i}^k(y)|=O(\varepsilon(1+|y|)^{2-m^*}+\varepsilon^{2m^*-3}(1+|y|)^{3-m^*}),\nonumber
\end{align}
\begin{align}
|w_{\xi,i}^k|\Big|_{|y|\sim\frac{1}{\varepsilon}}=O(\varepsilon^{m^*-1}+\varepsilon^{3m^*-6})\mbox{ and }
|\nabla w_{\xi,i}^k|\Big|_{|y|\sim\frac{1}{\varepsilon}}=O(\varepsilon^{m^*}+\varepsilon^{3m^*-5}).\nonumber
\end{align}
\end{corollary}

~~

\noindent{\bf Step 2. Induction on the vanishing orders.}

~~

To begin with, we state several assumption for some $M\in\mathbb{Z}_+$ 
\begin{align}\tag{$\mbox{A}_1$}\label{ineq:b0VanishingM-1}
|b_{0t}^k|=O(\varepsilon^{M(m^*-2){-\delta}}+\varepsilon)
\end{align}
and
\begin{align}\tag{$\mbox{A}_2$}\label{Expansion:XiOutM-1}
\xi_i^k(x)=\overline{\xi}_i^k+\overline{T^{\mbox{in}}_{i,k,2}}(x)+(2-m_i^k)b_{0t}^k\overline{T^{\mbox{out}}_{i,k,M}}\varepsilon^{m_i^k-2}+O(\varepsilon^{M(m^*-2){-\delta}}+\varepsilon^2)
\end{align}
with
\begin{align}
\overline{T^{\mbox{in}}_{i,k,2}}(x)=O(\varepsilon^{m^*-2}|b_{0}^k|+\varepsilon(|b_1^k|+|b_2^k|))=O(\varepsilon^{M(m^*-2){-\delta}}+\varepsilon)\nonumber
\end{align}
and
\begin{align}
(2-m_i^k)b_{0t}^k\overline{T^{\mbox{out}}_{i,k,2}}\varepsilon^{m_i^k-2}=O(\varepsilon^{m^*-2}|b_0^k|)=O(\varepsilon^{m^*-1}+\varepsilon^{M(m^*-2){-\delta}}).\nonumber
\end{align}
By (\ref{Expansion:XiOutM-1}), we get
\begin{align}\label{ineq:PhiM-1}
|\psi_{\xi,i}^k(y)|\leq C(\varepsilon^{M(m^*-2)+1{-\delta}}+\varepsilon^2)|y|.
\end{align}
Assume that for $|y|\leq\frac{\tau}{\varepsilon}$,
\begin{align}\tag{$\mbox{A}_3$}\label{ineq:wXiM-1}
|w_{\xi,i}^k(y)|=O(\varepsilon(1+|y|)^{2-m^*}+\varepsilon^{M(m^*-2)+1}(1+|y|)^{3-m^*}).
\end{align}

Now we prove higher order versions of (\ref{ineq:b0VanishingM-1}), (\ref{Expansion:XiOutM-1}), (\ref{ineq:PhiM-1}) and (\ref{ineq:wXiM-1}).
By (\ref{ineq:PhiM-1}) and (\ref{ineq:wXiM-1}), we get for $x\in M\backslash\cup_t B_t$
\begin{align}
\overline{\xi_i^k}=b_0^k Z_{0i}+O(\varepsilon+\varepsilon^{M(m^*-2)}).\nonumber
\end{align}
On the other hand, a Green function argument \'a la Claim \ref{c:ExpansionXiOut02} gives
\begin{claim}\label{c:XiOutM}
For $x\in M\backslash\cup_t B_t$, we get
\end{claim}
\begin{align}
\xi_i^k(x)=\overline{\xi}_i^k+\overline{T^{\mbox{in}}_{i,k,2}}(x)+(2-m_i^k)b_{0t}^k\overline{T^{\mbox{out}}_{i,k,M}}\varepsilon^{m_i^k-2}+O(\varepsilon^{(M+1)(m^*-2)}+\varepsilon^2)\nonumber
\end{align}
with
\begin{align}
\overline{T^{\mbox{in}}_{i,k,2}}(x)=O(\varepsilon^{m^*-2}|b_{0}^k|+\varepsilon(|b_1^k|+|b_2^k|))=O(\varepsilon^{(M+1)(m^*-2)}+\varepsilon)\nonumber
\end{align}
and
\begin{align}
(2-m_i^k)b_{0t}^k\overline{T^{\mbox{out}}_{i,k,2}}\varepsilon^{m_i^k-2}=O(\varepsilon^{m^*-2}|b_0^k|)=O(\varepsilon^{m^*-1}+\varepsilon^{(M+1)(m^*-2)}).\nonumber
\end{align}
An immediate observation is that
\begin{align}\label{ineq:PhiM}
|\psi_{\xi,i}^k(y)|\leq C(\varepsilon^{(M+1)(m^*-2)+1}+\varepsilon^2)|y|.
\end{align}
Plugging these into $E_{\xi,i}^k$ and taking the 0-frequency part, we get
\begin{align}
\int_{\partial( (\varepsilon_t^{(2),k})^{-1}B_t)}\frac{\partial w_{\xi,0,i}^k}{\partial{\bf n}}=O(|b_0^k|\varepsilon^{m^*-2}+\varepsilon^{(M+2)(m^*-2)}).\nonumber
\end{align}
These implies that 
\begin{align}\label{ineq:b0VanishingM}
|b_0^k|=O(\varepsilon^{(M+1)(m^*-2){-\delta}}+\varepsilon).    
\end{align}
Now we are in a position to derive a new bound for $w_{\xi,i}^k$. By a routine computation, we get
\begin{align}\label{ineq:wXiM}
|w_{\xi,i}^k(y)|=O(\varepsilon(1+|y|)^{2-m^*}+\varepsilon^{(M+1)(m^*-2)+1}(1+|y|)^{3-m^*})
\end{align}
for $|y|\leq\frac{\tau}{\varepsilon}$.

Notice that (\ref{ineq:b0VanishingM-1}), (\ref{Expansion:XiOutM-1}), (\ref{ineq:PhiM-1}) and (\ref{ineq:wXiM-1}) are improved by (\ref{ineq:b0VanishingM}), Claim \ref{c:XiOutM}, (\ref{ineq:PhiM}) and (\ref{ineq:wXiM}), respectively. The induction is completed.

\begin{remark}
It should be pointed out that for $|y|\sim\frac{1}{\varepsilon}$
\begin{align}\label{Expansion:Phi}
\psi_{\xi,i}^k=b_1^kZ_{1,i}+b_2^k Z_{2,i}+O(\varepsilon^{m^*-1}+\varepsilon^{(M+3)(m^*-2)}).
\end{align}
\end{remark}

\subsection{Analysis of $b_{1t}$ and $b_{2t}$ for $t=1,\cdots,N$}\label{Subsection:b122}

\begin{claim}
Under the assumptions of Theorem \ref{t:MAIN}, we get $b_1^k,b_2^k\to0$. 
\end{claim}
\noindent{\bf Proof.}
By Pohozaev identity \cite{LinZhang2013}, we get
\begin{align}
\varepsilon_t^{(1),k}\sum_i\int_{\Omega}\partial_l &H^{M,k}_i(\varepsilon_t^{(1),k}y+q_t^{(1),k})e^{v_{i}^{(1),k}}\nonumber\\
&=\int_{\partial\Omega}{\bf n}_l\sum_i e^{v_i^{(1),k}}H_i^{M,k}+\sum_{ij}a^{ij}\Big(\partial_{\bf n}v_i^{(1),k}\partial_l v_i^{(1),k}-\frac{1}{2}\nabla v_i^{(1),k}\cdot\nabla v_j^{(1),k}{\bf n}_l\Big)\nonumber
\end{align}
and
\begin{align}
\varepsilon_t^{(2),k}\sum_i\int_{\Omega}\partial_l &H^{M,k}_i(\varepsilon_t^{(2),k}y+q_t^{(2),k})e^{v_{i}^{(2),k}}\nonumber\\
&=\int_{\partial\Omega}{\bf n}_l\sum_i e^{v_i^{(2),k}}H_i^{M,k}+\sum_{ij}a^{ij}\Big(\partial_{\bf n}v_i^{(2),k}\partial_l v_i^{(2),k}-\frac{1}{2}\nabla v_i^{(2),k}\cdot\nabla v_j^{(2),k}{\bf n}_l\Big).\nonumber
\end{align}
Noticing that
$v_i^{(1),k}-v_i^{(2),k}=\sigma_k\xi_i^k$
with $\sigma_k=O(\varepsilon^{M(m^*-2)})$, we get
\begin{align}\label{e:PohozaevXi}
\sum_i\varepsilon_t^{(2),k}&\int_\Omega\partial_l H_i^{M,k}(p_t^{(2),k}+\varepsilon_t^{(2),k}y)e^{v_i^{(2),k}}\xi_i^k\nonumber\\
&=\sum_{ij}a^{ij}\int_{\partial\Omega}(\partial_{\bf n}\xi_i^k\partial_l v_j^{(2),k}+\partial_{\bf n}v_i^{(2),k}\partial_l \xi_j^k-\nabla v_i^{(2),k}\cdot\nabla\xi_j^k{\bf n}_l)+O(\varepsilon^{M(m^*-2)-1}).
\end{align}

~~

\noindent{\bf Analysis of RHS of (\ref{e:PohozaevXi})}

~~

Take $\Omega=B_{{\tau/\varepsilon}}(p_t^{(2),k})$ for small $\tau>0$.
For the right hand side of (\ref{e:PohozaevXi}), we get for the $l=1$ case 
\begin{align}
&\partial_{\bf n}\xi_i^k\partial_1 v_j^{(2),k}+\partial_{\bf n}v_i^{(2),k}\partial_1 \xi_j^k-\nabla v_i^{(2),k}\cdot\nabla\xi_j^k{\bf n}_l\nonumber\\
&=\partial_{\bf n}\Big(b_{0t}Z_{0,i}+b_{1t}Z_{1,i}+b_{2t}Z_{2,i} +\psi_{\xi,i}^k+w_{\xi,i}^k\Big)\times\partial_1(V_{tj}^{(2),k}+w_{t,j}^{(2),k})\nonumber\\
&\quad+\partial_{\bf n}(V_{tj}^{(2),k}+w_{t,j}^{(2),k})\times\partial_1\Big(b_{0t}Z_{0,i}+b_{1t}Z_{1,i}+b_{2t}Z_{2,i} +\psi_{\xi,i}^k+w_{\xi,i}^k\Big)\nonumber\\
&\quad-\nabla (V_{tj}^{(2),k}+w_{t,j}^{(2),k})\cdot\nabla \Big(b_{0t}Z_{0,i}+b_{1t}Z_{1,i}+b_{2t}Z_{2,i} +\psi_{\xi,i}^k+w_{\xi,i}^k\Big)\frac{y_1}{r}\nonumber
\end{align}
We first study the main terms. By a direct computation, we get
\begin{align}
\sum_{ij}a^{ij}\int_{\partial\Omega } \partial_{\bf n}Z_{0,i}\partial_1 V_{t,j}^{(2),k}+\partial_{\bf n}V_{t,j}\partial_1Z_{0,i}-(\partial_1 V_{t,j}^{(2),k}\partial_1 Z_{0,i}+\partial_2 V_{t,j}^{(2),k}\partial_2 Z_{0,i})\frac{y_1}{r}d\mathcal{H}^1=0,\nonumber
\end{align}
\begin{align}\label{e:PohozaevRHSZ2}
\sum_{ij}a^{ij}\int_{\partial\Omega} \partial_{\bf n}Z_{2,i}\partial_1 V_{t,j}^{(2),k}+\partial_{\bf n}V_{t,j}\partial_1Z_{2,i}-(\partial_1 V_{t,j}^{(2),k}\partial_1 Z_{2,i}+\partial_2 V_{t,j}^{(2),k}\partial_2 Z_{2,i})\frac{y_1}{r}d\mathcal{H}^1=0
\end{align}
and
\begin{align}\label{e:PohozaevRHSZ11}
&\sum_{ij}a^{ij}\Big[\partial_{\bf n}Z_{1,i}\partial_1 V_{t,j}^{(2),k}+\partial_{\bf n}V_{t,j}\partial_1Z_{1,i}-(\partial_1 V_{t,j}^{(2),k}\partial_1 Z_{1,i}+\partial_2 V_{t,j}^{(2),k}\partial_2 Z_{1,i})\frac{y_1}{r}\Big]\nonumber\\
&=\sum_{ij}a^{ij}\Big\{(V_{t,i}^{(2),k})''(V_{t,j}^{(2),k})'\frac{y_1^2}{r^2}+(V_{t,i}^{(2),k})'\cdot\Big[(V_{t,j}^{(2),k})''\frac{y_1^2}{r^2}+(V_{t,j}^{(2),k})'\Big(\frac{1}{r}-\frac{y_1^2}{r^3}\Big)\Big]
\nonumber\\
&\quad-(V_{t,i}^{(2),k})'\frac{y_1}{r}\cdot\Big[(V_{t,j}^{(2),k})''\frac{y_1^2}{r^2}+(V_{t,j}^{(2),k})'\Big(\frac{1}{r}-\frac{y_1^2}{r^3}\Big)\Big]\cdot\frac{y_1}{r}-(V_{t,i}^{(2),k})'\frac{y_2}{r}\cdot\Big[(V_{t,j}^{(2),k})''\frac{y_1 y_2}{r^2}-(V_{t,j}^{(2),k})'\frac{y_1 y_2}{r^3}\Big]\cdot\frac{y_1}{r}\nonumber\\
&=\sum_{ij}a^{ij}(V_{t,i}^{(2),k})'(V_{t,j}^{(2),k})''\cdot\Big(\frac{2y_1^2}{r^2}-\frac{y_1^4+y_1^2 y_2^2}{r^4}\Big)+\sum_{ij}a^{ij}(V_{t,i}^{(2),k})'(V_{t,j}^{(2),k})'\Big(\frac{1}{r}-\frac{2y_1^2}{r^3}+\frac{y_1^4}{r^5}+\frac{y_1^2 y_2^2}{r^5}\Big).
\end{align}
By symmetry, we get
\begin{align}\label{e:PohozaevRHSZ12}
&\sum_{ij}a^{ij}\int_{\partial\Omega}\Big[\partial_{\bf n}Z_{1,i}\partial_1 V_{t,j}^{(2),k}+\partial_{\bf n}V_{t,j}\partial_1Z_{1,i}-(\partial_1 V_{t,j}^{(2),k}\partial_1 Z_{1,i}+\partial_2 V_{t,j}^{(2),k}\partial_2 Z_{1,i})\frac{y_1}{r}\Big]\nonumber\\
&=\frac{1}{2}\sum_{ij}a^{ij}\int_{\partial\Omega}\Big((V_{t,i}^{(2),k})'(V_{t,j}^{(2),k})''+\frac{1}{r}(V_{t,i}^{(2),k})'(V_{t,j}^{(2),k})'\Big)\nonumber\\
&=\frac{1}{2}\sum_{ij}a^{ij}\int_{\partial\Omega}\Delta V_{i,t}^{(2),k} \partial_{\bf n}V_{t,j}^{(2),k}d\mathcal{H}^1\nonumber\\
&=-\frac{1}{2}\sum_{ij}a^{ij}\int_{\partial\Omega}e^{V_{t,j}^{(2),k}}\partial_{\bf n}V_{t,j}^{(2),k}d\mathcal{H}^1=O(\varepsilon^{m^*}).
\end{align}

Now we analysis the error parts, which reads as follows:
\begin{align}
&\sum_{ij}a^{ij}\int_{\partial\Omega}\Big\{\partial_{\bf n}\psi_{\xi,i}^k\partial_1 V_{t,j}^{(2),k}+\partial_{\bf n}\psi_{\xi,i}^k\partial_1 w_{t,j}^{(2),k}+\partial_{\bf n}w_{\xi,i}^k \partial_1 V_{tj}^{(2),k}+\partial_{\bf n}w_{\xi,i}^k\partial_1 w_{t,j}^{(2),k}\nonumber\\
&\quad\quad+\partial_1\psi_{\xi,i}^k\partial_{\bf n} V_{t,j}^{(2),k}+\partial_1\psi_{\xi,i}^k\partial_{\bf n} w_{t,j}^{(2),k}+\partial_1w_{\xi,i}^k \partial_{\bf n} V_{tj}^{(2),k}+\partial_1w_{\xi,i}^k\partial_{\bf n} w_{t,j}^{(2),k}\nonumber\\
&\quad\quad-\Big[\partial_1\psi_{\xi,i}^k\partial_1 V_{t,j}^{(2),k}+\partial_1\psi_{\xi,i}^k\partial_1 w_{t,j}^{(2),k}+\partial_1w_{\xi,i}^k \partial_1 V_{tj}^{(2),k}+\partial_1 w_{\xi,i}^k\partial_1 w_{t,j}^{(2),k}\nonumber\\
&\quad\quad+\partial_2\psi_{\xi,i}^k\partial_2 V_{t,j}^{(2),k}+\partial_2\psi_{\xi,i}^k\partial_2 w_{t,j}^{(2),k}+\partial_2w_{\xi,i}^k \partial_2 V_{tj}^{(2),k}+\partial_2 w_{\xi,i}^k\partial_2 w_{t,j}^{(2),k}\Big]\frac{y_1}{r}\Big\}d\mathcal{H}^1.\nonumber
\end{align}
{Recall that we assume that $\Omega=B_{{\tau/\varepsilon}}(p_t^{(2),k})$.}
In the above, all the terms formed by the product of the derivatives of $\psi_{\xi,i}^k$ and of $w_{t,j}^{(2),k}$ are of $\varepsilon^{m^*}$-order. All the terms formed by the product of the derivatives of $w_{\xi,i}^k$ and of $V_{t,j}^{(2),k}$ are of $\varepsilon^{m^*}$-order. All the terms formed by the product of the derivatives of $w_{\xi,i}^k$ and of $w_{t,j}^{(2),k}$ are of $\varepsilon^{2m^*-2}$-order. The only problem left are the terms formed by the product of the derivatives of $\psi_{\xi,i}^k$ and of $V_{t,j}^{(2),k}$.
By (\ref{Expansion:Phi}) and a similar computation as in (\ref{e:PohozaevRHSZ2}), (\ref{e:PohozaevRHSZ11}) and (\ref{e:PohozaevRHSZ12}), these terms are bounded by $O(\varepsilon^{m^*})$.
Therefore, we get
\begin{align}
\mbox{LHS of (\ref{e:PohozaevXi})}=O(\varepsilon^{m^*}).\nonumber
\end{align}

~~

\noindent{\bf Analysis of LHS of (\ref{e:PohozaevXi})}

~~

Now we analyze the LHS of (\ref{e:PohozaevXi}).
An easy computation gives
\begin{align}
\mbox{LHS of (\ref{e:PohozaevXi})}&=(\varepsilon_1^{(2),k})^2\Bigg(\sum_i\partial_{11}H_i^{M,k}(0)\int_{ (\varepsilon_t^{(2),k})^{-1}B_t}e^{V_{it}^{(2),k}}V'_i\frac{y_1^2}{r}dy\cdot b_1^k\nonumber\\
&\quad\quad
+\sum_i\partial_{12}H_i^{M,k}(0)\int_{(\varepsilon_t^{(2),k})^{-1}B_t}e^{V_{it}^{(2),k}}V'_i\frac{y_2^2}{r}dy\cdot b_2^k\Bigg)(1+o_\varepsilon(1))\nonumber\\
&=(\varepsilon_1^{(2),k})^2\Bigg(\sum_i\partial_{11}H_i^{M,k}(0)\int_{\mathbb{R}^2}e^{V_{it}^{(2),k}}dy\cdot b_1^k+\sum_i\partial_{12}H_i^{M,k}(0)\int_{\mathbb{R}^2}e^{V_{it}^{(2),k}}dy\cdot b_2^k\Bigg)(1+o_\varepsilon(1))\nonumber
\end{align}
Here, we use the computation that
\begin{align}
\int_{(\varepsilon_t^{(2),k})^{-1}B_t} e^{V_{it}^{(2),k}}V'_i\frac{y_1^2}{r}dy=\pi\int_0^{\tau/\varepsilon_t^{(2),k}}e^{V_{it}^{(2),k}}V'_ir^2dr=-\pi\int_{r=0}^{r={\tau/\varepsilon_t^{(2),k}}}r^2d\Big(e^{V_{it}^{(2),k}}\Big)=2\pi\int_0^{\tau/\varepsilon_t^{(2),k}}e^{V_{it}^{(2),k}}rdr.\nonumber
\end{align}
The non-degeneracy implies that $b^k_1,b^k_2\to0$ as $k\to\infty$.

\noindent{$\Box$}


\subsection{Proof of Theorem \ref{t:MAIN} and Theorem \ref{non-deg}}

In this subsection, we prove Theorem \ref{t:MAIN} and Theorem \ref{non-deg}.

~~

\noindent{\bf Proof of Theorem \ref{t:MAIN}.}
Based on the above estimate, we know that $\|\xi_i^k\|_{L^\infty(M)}\to0$ for any $i=1,\cdots,n$ as $k\to\infty$. This is a contradiction with the definition of $\xi_i^k$. Therefore $\xi_i^k\equiv0$ for $i=1,\cdots,n$ and $k=1,2,\cdots$. This implies that $u_i^{(1),k}=u_i^{(2),k}$ for any $i=1,\cdots,n$ and $k$ large.

\noindent{$\Box$}

~~

\noindent{\bf Proof of Theorem \ref{non-deg}.}
Assume that there exists a sequence of nontrivial Sobolev functions $\phi_i^k\in \mathring{H^1}(M)$ for $i=1,\cdots,n$ and $k=1,2,\cdots$ with
\begin{align}
\Delta\Big(\sum_j a^{ij}\phi_j^k\Big)+\rho_i^k\frac{h_ie^{u_i^k}}{\int_M h_ie^{u_i^k}}\Bigg(\phi_i^k-\int_M\frac{h_i e^{u_i^k}\phi_i^k}{\int_M h_i e^{u_i^k}}\Bigg)=0\mbox{ in }M.\nonumber
\end{align}
Up to the renormalization $\int_M h_i e^{u_i^k}\equiv 1$ for $i=1,\cdots,n$ and $k=1,2,\cdots$, we get
\begin{align}
\Delta\Big(\sum_j a^{ij}\phi_j^k\Big)+\rho_i^k h_ie^{u_i^k}\Big(\phi_i^k-\int_M h_i e^{u_i^k}\phi_i^k\Big)=0\mbox{ in }M.\nonumber
\end{align}
Denoting
$\xi_i^k:=\frac{\phi_i^k-\int_M h_i e^{u_i^k}\phi_i^k}{\big\|\phi_i^k-\int_M h_i e^{u_i^k}\phi_i^k\big\|_{L^\infty(M)}}$
and
repeating the argument for the uniqueness result, we prove that $\xi_i^k\equiv0$. Therefore, we get $\phi_i^k\equiv0$. The non-degeneracy result follows.

\noindent{$\Box$}

\appendix

\section{Pohozaev identity and an invertibility result}\label{App:Solvability}

In order to build a leading term for 1-frequency part, we study a modified invertibility of the linearized operator.
This is a Fredholm theory in the spirit of arguments in \cite{LinYan2013,Huang2019,ChengLiZhang2025}. {The central idea is to \emph{use Pohozaev identity as a criterion for the existence of 1-frequency part.} In this part, for the sake of simplicity, we denote $\tau_k^M$ and $\varepsilon_t^{(l),k}$ as $\tau$ and $\varepsilon$, respectively.}

\subsection{A Fredholm result and solvability of a linear problem}

{In this subsection, we first build a Fredholm result (Proposition \ref{prop:invertibility}). Then, we show that the Pohozaev identity implies the existence of 1-frequency part of the solution (Corollary \ref{coro:1-FreqSolvable}).}
We begin with two weighted Sobolev spaces.
Denote
\begin{align}
\rho_\beta(x)=(1+|x|)^{1+\frac{\beta}{2}}\mbox{ and }\widetilde{\rho}_\beta(x)=\frac{1}{(1+|x|)(\ln(2+|x|))^{1+\frac{\beta}{2}}}\nonumber
\end{align}
for a sufficiently small positive number $\beta$.
\begin{definition}
\begin{align}
\mathbb{X}=\{u\in(L_{loc}(B_{3\tau/\varepsilon}))^n|\|u\|_{\mathbb{X}}<\infty,\,u_i|_{\partial B_{3\tau/\varepsilon}}\equiv 0\}\mbox{ and }\mathbb{Y}=\{u\in(L_{loc}(B_{3\tau/\varepsilon}))^n|\|u\|_{\mathbb{Y}}<\infty\}.\nonumber
\end{align}
Here,
\begin{align}
\|u\|_{\mathbb{X}}^2=\sum_{i=1}^n\Big(\|\Delta u_i\rho_\beta\|_{L^2(B_{3\tau/\varepsilon})}^2+\|u_i\widetilde{\rho}_\beta\|_{L^2(B_{3\tau/\varepsilon})}\Big)\mbox{ and }
\|u\|_{\mathbb{Y}}^2=\sum_{i=1}^n\|u_i\rho_\beta\|^2_{L^2(B_{3\tau/\varepsilon})}.\nonumber
\end{align}
\end{definition}
\begin{definition}
We denote $\mathbb{X}_{odd}$ and $\mathbb{Y}_{odd}$ the subspaces of  $\mathbb{X}$ and $\mathbb{Y}$ consisting of their odd functions, respectively.    
\end{definition}

For the case $\partial_1 h_i^{M,k}(0)\neq 0$, we proceed the following modification on the kernel. We start by defining several functions.
Let us introduce two smooth cut-off functions $\chi$ and $\chi_*$. Here, we assume that $\chi(x)=1$ if $|x|\leq\tau$ and $\chi(x)=0$ if $|x|\geq2\tau$ and $\chi_*(x)=1$ if $\frac{7\tau}{3}\leq|x|\leq\frac{8\tau}{3}$ and $\chi_*(x)=0$ if $|x|\leq 2\tau$ or $|x|\geq 3\tau$. 
Then, we define
\begin{align}
Z_{\varepsilon,0,i}(y)=Z_{0,i}(y)\cdot\chi(\varepsilon y)\nonumber
\end{align}
and
\begin{align}
Z_{\varepsilon,s,i}(y;t)=Z_{s,i}(y)\cdot\chi(\varepsilon y)+t \chi_*(\varepsilon y)\cdot Z_{s,i}(y)\cdot\mbox{sgn}(\partial_s h_i^{M,k}(0))\nonumber
\end{align}
for $t\in\mathbb{R}$, $s=1,2$ and $i=1,\cdots,n$. 
For the case $\partial_1h_i^{M,k}(0)\equiv0$ and $V_i^{M,k}=V_j^{M,k}$ for $i\neq j$, we set $t=0$.

\begin{definition}
Define the weighted spaces as follows.
\begin{align}
E^{odd}_\varepsilon(t_1,t_2)=\Bigg\{u=(u_1,\cdots,,u_n)\subset\mathbb{X}_{odd}\Bigg|\int_{B_{3\tau/\varepsilon}}\sum_{i=1}^n u_iY_{\varepsilon,1,i}(y;t_1)=\int_{B_{3\tau/\varepsilon}}\sum_{i=1}^n u_iY_{\varepsilon,2,i}(y;t_2)=0\Bigg\}\nonumber
\end{align}
and
\begin{align}
F^{odd}_\varepsilon(t_1,t_2)=\Bigg\{u=(u_1,\cdots,,u_n)\subset\mathbb{Y}_{odd}\Bigg|\int_{B_{3\tau/\varepsilon}}\sum_{i=1}^n u_iZ_{\varepsilon,1,i}(y;t_1)=\int_{B_{3\tau/\varepsilon}}\sum_{i=1}^n u_iZ_{\varepsilon,2,i}(y;t_2)=0\Bigg\}.\nonumber
\end{align}
Here, for $s=1,2$
\begin{align}
Y_{\varepsilon,s,i}(y;t):=\Delta \Big(\sum_j a^{ij} Z_{\varepsilon,s,j}(y;t)\Big).\nonumber
\end{align}
\end{definition}

Define the projection $Q_\varepsilon:\mathbb{Y}_{odd}\to F_\varepsilon^{odd}$ as follows.
\begin{align}
Q_\varepsilon(u)=u-\sum_{j=1}^2 r_j\Bigg(e^{V_1^{M,k}}Z_{\varepsilon,j,1}(y;t_j^\varepsilon),\cdots,e^{V_n^{M,k}}Z_{\varepsilon,j,n}(y;t_j^\varepsilon)\Bigg).\nonumber
\end{align}
Here, the numbers $r_1$ and $r_2$ are chosen such that $Q_\varepsilon(u)\in F_\varepsilon^{odd}$.
Now we study the invertibility result for the operator
{
$L_\varepsilon=(L_{1,\varepsilon},\cdots,L_{n,\varepsilon}):\mathbb{X}_{odd}\to\mathbb{Y}_{odd}$
\begin{align}
L_{\varepsilon}(u)=(L_{1,\varepsilon}(u),\cdots,L_{n,\varepsilon}(u))\,\mbox{ and }\,L_{i,\varepsilon}(u)=\Delta \Big(\sum_j a^{ij} u_j\Big)+h_i^{M,k}(0)e^{V_i^{M,k}}u_i\nonumber
\end{align}
for $i=1,\cdots,n$
satisfying
\begin{align}
L_{\varepsilon}(u)=(h_{1,\varepsilon},\cdots,h_{n,\varepsilon})+\sum_{j=1}^2 r_{\varepsilon,j}\Bigg(e^{V_1^{M,k}}Z_{\varepsilon,j,1},\cdots,e^{V_n^{M,k}}Z_{\varepsilon,j,n}\Bigg)\nonumber
\end{align}
with suitable $r_{\varepsilon,1}$ and $r_{\varepsilon,2}$ such that $(h_{1,\varepsilon},\cdots,h_{n,\varepsilon})\in F_\varepsilon^{odd}$.
Our main result is that
\begin{proposition}\label{prop:invertibility}
There exists a positive constant $C$ such that
\begin{align}\label{ineq:Invert}
\|u\|_\mathbb{X}+\|u\|_{L^\infty}\leq C\|Q_\varepsilon L_\varepsilon u\|_\mathbb{Y}
\end{align}
for any $u\in \mathbb{X}_{odd}$. Moreover, $Q_\varepsilon L_\varepsilon$ is a homeomorphism from $E_\varepsilon^{odd}$ to $F_\varepsilon^{odd}$.
\end{proposition}

Based on the above notions, we have the following lemmas.
\begin{lemma}\label{l:Ortho1-Freq}
There exists $t_1^\varepsilon,t_2^\varepsilon\in\mathbb{R}$ such that
\begin{align}
\sum_{i=1}^N\int_{B_{3\tau/\varepsilon}}(y\cdot\nabla h_i^{M,k}(0))e^{V_i^{1,k}}Z_{\varepsilon,s,i}(y,t^\varepsilon_s)=0.\nonumber
\end{align}
for $s=1,2$. In other words, we get
\begin{align}\label{e:Orthogonal01}
(y\cdot\nabla h_1^{M,k}(0)e^{V_1^{M,k}},\cdots,y\cdot\nabla h_n^{M,k}(0)e^{V_n^{M,k}})\in F_\varepsilon^{odd}(t_i^\varepsilon,t_2^\varepsilon)
\end{align}
Moreover, we get $|t^\varepsilon_1|+|t^\varepsilon_2|=O(1)$ independence in $\varepsilon$.
\end{lemma}
\noindent{\bf Proof.}
We only prove the case $s=1$ since the other is similar.
By a direct computation, we get
\begin{align}
0&=\sum_i\int_{B_{3\tau/\varepsilon}}(y\cdot\nabla h_i^{M,k}(0))e^{V_i^{M,k}}Z_{\varepsilon,1,i}(y,t^\varepsilon_s)\nonumber\\
&=\sum_i\int_{B_{3\tau/\varepsilon}}(y\cdot\nabla h_i^{M,k}(0))e^{V_i^{M,k}}\Big(Z_{s,i}(y)\cdot\chi(\varepsilon y)+t_1^\varepsilon \chi_*(\varepsilon y)\cdot Z_{1,i}(y)\cdot\mbox{sgn}(\partial_s h_i^{M,k}(0))\Big)\nonumber\\
&=O(\varepsilon^{m^*-2-\delta})+t_1^\varepsilon\int_{B_{3\tau/\varepsilon}\backslash B_{2\tau/\varepsilon}}\chi_*(\varepsilon y) \sum_{i}y_1 Z_{1,i}(y) e^{V_i^{M,k}}|\partial_1 h_1^{M,k}(0)|\chi_*(\varepsilon y).\nonumber
\end{align}
{Here, we use a Pohozaev identity in the form of \cite[Theorem 1.3]{LinZhang2010}.}
This implies that
$|t_1^\varepsilon|=O(1)$.

\noindent{$\Box$}

An important corollary of Lemma \ref{l:Ortho1-Freq} is the solvability of the equation
\begin{align}\label{e:1-FreqLeadingSolvability}
\Delta\Bigg(\sum_{j=1}^n a^{ij}C^k_{1,j}\Bigg)+h^{M,k}_i(0)e^{V^{1,k}_i}C^k_{1,i}=-\nabla h^{M,k}_i(0)\cdot y e^{V^{1,k}_i}\mbox{ in }B_{3\tau/\varepsilon}
\end{align}
for $i=1,\cdots,n$. To be precise, we get
\begin{corollary}\label{coro:1-FreqSolvable}
For the numbers $t_\varepsilon^1$ and $t_\varepsilon^2$ in Lemma \ref{l:Ortho1-Freq}, there exists a unique $(C^k_{1,1},\cdots,C_{1,n}^k)\in E_\varepsilon^{odd}(t_1^\varepsilon,t_2^\varepsilon)$ satisfies (\ref{e:1-FreqLeadingSolvability}).
\end{corollary}

\begin{remark}\label{r:MultiOrthognal}
Now we are in a position to establish a generalization of the above constructions. In order to establish the solvability of 
\begin{align}
\Delta\Bigg(\sum_{j=1}^n a^{ij}C^k_{1,l,j}\Bigg)+H^{M,k}_i(0)e^{V^{1,k}_i}C^k_{1,l,i}=-\nabla h^{M,k}_i(0)\cdot y e^{V^{1,k}_i}+A_{l,i}\mbox{ in }B_{3\tau/\varepsilon}\nonumber
\end{align}
with $u_{1,l,j}\equiv0$ on $\partial B_{{3\tau/\varepsilon_k}}$
for $l=1,\cdots,L$ and $A_l$'s denotes small perturbations with their frequencies equal to 1, we 
take a further step to modify the functions $Z_{\varepsilon,1,i}$ as follows. Taking $\Omega_l\subset B_{3\tau}\backslash B_{2\tau}$ for $l=1,\cdots,L$ such that for any $l_1\neq l_2$ we get $\Omega_{l_1}\cap \Omega_{l_2}=\emptyset$. Then, we define $L$ smooth nonnegative cut-off functions $\chi_{*,l}$ supported on $\Omega_l$. Taking $t_1,\cdots,t_L\in\mathbb{R}$ and defining
\begin{align}
Z_{\varepsilon,s,i}(y;t_{s,1},\cdots,t_{s,L})=Z_{s,i}(y)\cdot\chi(\varepsilon y)+\sum_{l=1}^L t_{s,l} \chi_{*,L}(\varepsilon y)\cdot Z_{s,i}(y)\cdot\mbox{sgn}(\partial_s h_i^{M,k}(0))\nonumber
\end{align}
$s=1,2$ and $t=1,\cdots,n$. As in the above computation, we determine the numbers $t_{s,1},\cdots,t_{s,l}$ using the Pohozaev identity and the smallness of $A_{i,l}$'s.
Writing $T_s:=(t_{s,1},\cdots,t_{s,L})$ with $s=1,2$ for short, we define the weighted spaces as follows.
\begin{align}
E^{odd}_\varepsilon(T_1,T_2)=\Bigg\{u=(u_1,\cdots,,u_n)\subset\mathbb{X}_{odd}\Bigg|\int_{B_{3\tau/\varepsilon}}\sum_{i=1}^n u_iY_{\varepsilon,1,i}(y;T_1)=\int_{B_{3\tau/\varepsilon}}\sum_{i=1}^n u_iY_{\varepsilon,2,i}(y;T_2)=0\Bigg\}\nonumber
\end{align}
and
\begin{align}
F^{odd}_\varepsilon(T_1,T_2)=\Bigg\{u=(u_1,\cdots,,u_n)\subset\mathbb{Y}_{odd}\Bigg|\int_{B_{3\tau/\varepsilon}}\sum_{i=1}^n u_iZ_{\varepsilon,1,i}(y;T_1)=\int_{B_{3\tau/\varepsilon}}\sum_{i=1}^n u_iZ_{\varepsilon,2,i}(y;T_2)=0\Bigg\}.\nonumber
\end{align}
Here, for $s=1,2$
\begin{align}
Y_{\varepsilon,s,i}(y;T):=\Delta \Big(\sum_j a^{ij} Z_{\varepsilon,s,j}(y;T)\Big)\,\mbox{ with }\,T\in\mathbb{R}^L.\nonumber
\end{align}
An analogue to (\ref{e:Orthogonal01}) follows. To be precise, we get
 for certain $T_1^\varepsilon,T_2^\varepsilon\in \mathbb{R}^L$
\begin{align}
(y\cdot\nabla h_1^{M,k}(0)e^{V_1^{M,k}}+A_{1,l},\cdots,y\cdot\nabla h_n^{M,k}(0)e^{V_n^{M,k}}+A_{n,l})\in F^{odd}_\varepsilon(T^\varepsilon_1,T^\varepsilon_2)\nonumber
\end{align}
holds for any $l=1,\cdots,L$. Analogues to Proposition \ref{prop:invertibility} and Corollary \ref{coro:1-FreqSolvable} for $E_\varepsilon^{odd}(T_1^\varepsilon,T_2^\varepsilon)$ and $F_\varepsilon^{odd}(T_1^\varepsilon,T_2^\varepsilon)$ holds.
\end{remark}

\subsection{Proof of Proposition \ref{prop:invertibility}}

To prove Proposition \ref{prop:invertibility}, we first estimate the decay rate of $r_{1,\varepsilon}$ and $r_{2,\varepsilon}$.
\begin{lemma}\label{l:Coefficientrj}
It holds that
$r_{\varepsilon,j}=O(\|h_\varepsilon\|_\mathbb{Y})+O(\varepsilon^{m^*-2-\delta})\max_i\|u_i\|_{L^\infty}$
for $j=1,2$ and a small positive constant $\delta$.
\end{lemma}
\noindent{\bf Proof.}
Multiplying the $i$-th equation by $Z_{\varepsilon,j,i}$ and integrating over $B_{3\tau/\varepsilon}$, we get
\begin{align}\label{e:CoefficientTOTAL}
{\sum_{i=1}^n}\int_{B_{3\tau/\varepsilon}}L_{i,\varepsilon}(u)Z_{i,\varepsilon}={\sum_{i=1}^n}\int_{B_{3\tau/\varepsilon}}h_{i,\varepsilon}Z_{\varepsilon,j,i}+r_{\varepsilon,j}{\sum_{i=1}^n}\int_{B_{3\tau/\varepsilon}}e^{V_i^{M,k}}Z_{\varepsilon,j,i}^2.
\end{align}
For the left hand side, we get
\begin{align}
\sum_{i=1}^n\int_{B_{3\tau/\varepsilon}}L_{i,\varepsilon}(u)Z_{i,\varepsilon}&=\sum_{i=1}^n\int_{B_{3\tau/\varepsilon}}\Bigg(\Delta\Big(\sum_{j=1}^n a^{ij} u_j\Big)+h_i^{M,k}(0)e^{V_i^{M,k}}u_i\Bigg) Z_{\varepsilon,s,i}\nonumber\\
&=\sum_{i=1}^n\int_{B_{3\tau/\varepsilon}}u_i\Bigg(\Delta \Big(\sum_{j=1}^n a^{ij} Z_{\varepsilon,s,i}\Big)+h_i^{M,k}(0)e^{V_i^{M,k}}Z_{\varepsilon,s,i}\Bigg)\nonumber\\
&=O(\varepsilon^{m^*-2-\delta})\max_{i}\|u_i\|_{L^\infty}.\nonumber
\end{align}
In the second equality, we use divergent theorem and the construction of $Z_{\varepsilon,j,i}$. For the left hand side, we get
\begin{align}
{\sum_{i=1}^n}\int_{B_{3\tau/\varepsilon}}h_{i,\varepsilon}Z_{\varepsilon,j,i}+r_{\varepsilon,j}{\sum_{i=1}^n}\int_{B_{3\tau/\varepsilon}}e^{V_i^{M,k}}Z_{\varepsilon,j,i}^2=O(\|h_\varepsilon\|_\mathbb{Y})+(C+o(1))r_{\varepsilon,j}\nonumber
\end{align}
for certain positive constant $C$.
This implies the result.

\noindent{$\Box$}

~~

\noindent{\bf Proof of Proposition 
\ref{prop:invertibility}.}
We prove the result by contradiction. Assume that there exist sequences $\{\varepsilon_m\}_m$ and $u_m\in E_{\varepsilon_m}^{odd}$ such that
\begin{align}
\|u_m\|_{\mathbb{X}}+\|u_m\|_{L^\infty}=1\mbox{ and }\|QL u_m\|_{\mathbb{Y}}=o(1).\nonumber
\end{align}

~~

\noindent{\bf Step 1. Estimate in the $L_{loc}^\infty$-sense.}

~~

We claim that $\|u_m\|_{L^\infty}=o(1)$. To prove this, we begin by a $L^p$ estimate as follows.

\begin{align}
\|QLu_m\|_{L^p}&\leq\|Lu_m\|_{L^p}=\Bigg(\int_{B_{3\tau/\varepsilon}}|Lu_m|^p\frac{\rho_\beta^p}{\rho_\beta^p}dx\Bigg)^\frac{1}{p}\nonumber\\
&\leq\Bigg(\int_{B_{3\tau/\varepsilon}}|Lu_m|^2\rho_\beta^2\Bigg)^\frac{1}{2}\Big(\int_{\mathbb{R}^2}\frac{dx}{\rho_\beta^\frac{2p}{2-p}}\Big)^\frac{2-p}{2p}\leq C\|Lu_m\|_{\mathbb{Y}}=o(1)\nonumber
\end{align}
for suitable $p$ close to 2 and smaller than 2. Therefore, we see that
\begin{align}
-\Delta\Big(\sum_{j=1}^n a^{ij}u_{m,j}\Big)-h_j^{M,k}(0)e^{V_i^{M,k}}u_{m,i}=o_{L^p}(1)\mbox{ in }\mathbb{R}^2\mbox{ for }i=1,\cdots,n.\nonumber
\end{align}
Here, $o_{L^p}(1)$ is the infinitesimal in $L^p$-norm.
Therefore, $u_{m,i}\to u_{\infty,i}$ in $C_{loc}^2$ with
$(u_{\infty,1},\cdots,u_{\infty,n})\in\mbox{span}\{Z_{0},Z_1,Z_2\}$, we see that $u_{\infty,i}=0$ for any $i=1,\cdots,n$.

~~

\noindent{\bf Step 2. Estimate in $L^\infty$ space.}

~~

We argue by contradiction. Assume that there exists $x_m\in B_{2\tau/\varepsilon}$ such that there exists $i_0\in I$ and a constant $\delta_0>0$ such that
$|u_{m,i_0}(x_m)|\geq\delta_0$.
By Step 1, we see that $|x_m|\to\infty$. Assume that $|x_m|/\varepsilon_m=o_m(1)$.
Rescaling the function as follows.
\begin{align}
\overline{u}_{m,i}(y)=u_{m,i}(|x_k|y)\mbox{ for }i=1,\cdots,n.\nonumber
\end{align}
Then,
\begin{align}
-\Delta_y \overline{u}_{m,i}=|x_m|^2\Big((L_i u_m)(|x_m|y)+\sum_{j=1}^n a_{ij} h^{M,k}_j(0)e^{V_j^{1,k}(|x_k|y)}\overline{u}_{m,j}+r_{\varepsilon,j}e^{V_i^{M,k}}Z_{\varepsilon,j,i}\Big)\mbox{ for }i=1,\cdots,n.\nonumber
\end{align}
Select any $R_1,R_2\in(0,+\infty)$ with $R_1<R_2$, we get
\begin{align}
\| |x_m|^2(L_i u_m)\|_{L^p(B_{R_1,R_2};dy)}
&=|x_m|^{2-\frac{2}{p}}\Bigg(\int_{B_{|x_m|R_1,|x_m|R_2}}|L_i u_m(x)|^pdx\Bigg)^\frac{1}{p}\nonumber\\
&\leq |x_m|^{2-\frac{2}{p}}\Bigg(\int_{B_{|x_m|R_1,|x_m|R_2}}|L_i u_m(x)|^2\rho_\beta^2 dx\Bigg)^\frac{1}{2}\Bigg(\int_{B_{|x_m|R_1,|x_m|R_2}}\rho_\beta^{-\frac{2p}{2-p}}\Bigg)^\frac{2-p}{2p}\nonumber\\
&\leq C |x_m|^{-\frac{\beta}{2}}\|Lu_m\|_{\mathbb{Y}},\nonumber
\end{align}
\begin{align}
||x_m|^2 e^{V_j^{1,k}(|x_m|y)}\overline{u}_{m,j}|=O\Bigg(\frac{|x_m|^{2-m^*}}{|y|^{m^*}}\Bigg)\mbox{ for }R_1<|y|<R_2\nonumber
\end{align}
and
\begin{align}
||x_m|^2r_{\varepsilon,j}e^{V_i^{M,k}(|x_m|y)}Z_{\varepsilon,j,i}(|x_\varepsilon|y)|=O\Bigg(\frac{|x_m|^{2-m^*}}{|y|^{m^*}}\Bigg)\mbox{ for }R_1<|y|<R_2.\nonumber
\end{align}
For the last assertion we use Lemma \ref{l:Coefficientrj}.
Therefore, $\overline{u}_i$ converges to a harmonic function in $B_{R_1,R_2}$. Since $R_1$ and $R_2$ are arbitrarily selected, the harmonic function is defined on $\mathbb{R}^2\backslash\{0\}$. On the other hand, $u_m$ are bounded. Hence, $\overline{u}_{m,i}\to(c_1,\cdots,c_n)$ with $c_i$ is constant for any $i=1,\cdots,n$. In other words, we get
\begin{align}\label{e:wkiAsy}
u_{m,i}(x)=c_i+o_m(1)\mbox{ for }R_1|x_m|\leq|x|\leq R_2|x_m|.
\end{align}
Here, $\max_i|c_i|\geq\frac{\delta_0}{2}>0$.
In the case $\frac{|x_m|}{\varepsilon_m}=C+o_m(1)$ for some positive constant $C$, we obtain a contradiction with the boundary condition of $u_{m,i}$.

On the other hand, $w_k\in E_{odd}$, this is a contradiction with the assumption. Therefore, we know that
\begin{align}\label{ineq:wkLinftyo1}
\|u_m\|_{L^\infty}=o_k(1).
\end{align}

~~

\noindent{\bf Step 3. Estimating the $\mathbb{X}$-norm.}

~~

In this part, we estimate the $\mathbb{X}$-norm of $u_m$.
By a direct computation,
we get

\begin{align}
\|\Delta u_{m,i}\rho_\beta\|_{L^2}\leq \sum_i\|e^{V_i^{M,k}}u_{m,i}\rho_\beta
\|_{L^2}+\|L_i u\rho_\beta\|_{L^2}\leq C\sum_{i}\|u_{m,i}\|_{L^\infty}+{\|QL_\varepsilon u\|_{\mathbb{Y}}+o_\varepsilon(1)}\nonumber
\end{align}
and
\begin{align}
\|u_{m,i}\widetilde{\rho}_\beta\|_{L^2}\leq C\sum_{i}\|u_{m,i}\|_{L^\infty}.\nonumber
\end{align}

This implies that
\begin{align}\label{ineq:wkX}
\|u_m\|_{\mathbb{X}}=o_k(1).
\end{align}

~~

\noindent{\bf Step 4. Completing the proof.}

~~

Combining (\ref{ineq:wkLinftyo1}) and (\ref{ineq:wkX}), we get
$\|u_m\|_{L^\infty}+\|u_m\|_{\mathbb{X}}=o(1)$.
This is a contradiction. {Now we prove (\ref{ineq:Invert})}.

To complete the proof, we only need to check that $Q_\varepsilon L_\varepsilon: E_\varepsilon^{odd}(t_1^\varepsilon,t_2^\varepsilon)\to F_\varepsilon^{odd}(t_1^\varepsilon,t_2^\varepsilon)$ is surjective. On one hand, we find that $L^*: E_\varepsilon^{odd}(t_1^\varepsilon,t_2^\varepsilon)\to F_\varepsilon^{odd}(t_1^\varepsilon,t_2^\varepsilon)$ defined as $$L^*u=\Big(\Delta\Big(\sum_j a^{1j}u_1\Big),\cdots,\Delta\Big(\sum_j a^{nj}u_n\Big)\Big)$$ is an isomorphism. On the other hand, for any $h=(h_1,\cdots,h_n)\in F_\varepsilon^{odd}(t_1^\varepsilon,t_2^\varepsilon)$, we rewrite
\begin{align}
Q_\varepsilon L_\varepsilon u_h =h\nonumber
\end{align}
as
\begin{align}
u_h=(L^*)^{-1}(L^* u_h-L_\varepsilon u_h)+(L^*)^{-1}h.\nonumber
\end{align}
Here, $(L^*)^{-1}(L^*-L_\varepsilon)$ is compact.  Thus the result follows from the Fredholm alternative.

\noindent{$\Box$}

\subsection{A remark on a non-vanishing blowup result}\label{subsectionb:NonvanishBlowup}

In \cite{ChengLiZhang2025}, the authors glue a blowup solution under the assumptions $(\mathcal{H}_1)$,  (\ref{SharpEstimate1}), (\ref{SharpEstimate2}) and (\ref{e:LHitHis}). The assumption (\ref{SharpEstimate1}) does not require 
$$\Big[\nabla \ln h_i(q_{t}^{k})+2\pi m_i^* \nabla_1 G^*(q_t^k,q_t^k)\Big]\rho_i^*\to0\,\,\forall i.$$
{Note that such a requirement is also absent in \cite{ChengLiZhang2025}, in which the authors provide sufficient conditions for the existence of the blowup solutions. Therefore, it is possible to find indices $t_0=1,\cdots,N$ and $i_0=1,\cdots,n$ such that $\Big|\nabla \ln h_{i_0}(q_{t_0}^{k})+2\pi m_{i_0}^* \nabla_1 G^*(q_{t_0}^k,q_{t_0}^k)\Big|\rho_{i_0}^*\geq C> 0$ for an uniform constant $C$. In this viewpoint, we say that Liouville systems have {\em a non-vanishing blowup result}.}

This is drastically different from the case of single equations. In \cite{BartYangZhang20241}, the authors prove strong vanishing results at blowup points, which play a crucial role in their uniqueness results.

~~

\noindent{{\bf Acknowledgment.} The authors express their gratitude to Dr. Xiaopeng Huang for his technical assistant.}

\noindent{{\bf Author Contributions.} The authors are listed in alphabetical order of their surnames, and all authors have made equal contributions.}

\noindent{{\bf Competing Interests.} The authors have no relevant financial or non-financial competing interests.}

\noindent{{\bf Data Availability Statement.} No datasets were generated or analysed during the current study.}

\noindent{{\bf Funding.} Lei Zhang is partially supported by Simon's Foundation grant SFI-MPS-TSM-00013752}

~~

\Vs\Vs

\end{document}